\documentclass[12pt]{article}

\usepackage{epic}
\usepackage{eepic}
\usepackage{epsf} 
\usepackage{ amssymb, amscd}
\usepackage{amsmath}
\usepackage{color}
\usepackage{epsfig}
\usepackage{graphicx}

\numberwithin{equation}{section}

\def\cA            {{\mathcal{A}}}

\def\cC            {{\mathcal{C}}}
\def\cD            {{\mathcal{D}}}
\def\cE            {{\mathcal{E}}}
\def\cF            {{\mathcal{F}}}

\def\cL            {{\mathcal{L}}}
\def\cM            {{\mathcal{M}}}

\def\cO            {{\mathcal{O}}}

\def\cR            {{\mathcal{R}}}

\def\cT            {{\mathcal{T}}}

\def\cX            {{\mathcal{X}}}

\def\cZ            {{\mathcal{Z}}}

\def\bbC           {\mathbb{C}}

\def\bbK            {\mathbb{K}}

\def\bbQ           {\mathbb{Q}}

\def\bbT           {\mathbb{T}}
\def\bbZ           {\mathbb{Z}}

\def\a             {\alpha}

\def\i{{\rm i}}
\def\sdprod{{\times\!\vrule height5pt depth0pt width0.4pt\,}}

\def\bbF      {\mathbb{F}}

\title{Near-group fusion categories and their doubles}
\author{
{\sc David E.\ Evans}\\
 {\footnotesize School of Mathematics, Cardiff University,}\\
 {\footnotesize Senghennydd Road, Cardiff CF24 4AG, Wales, U.K.}\\
 {\footnotesize e-mail: {\tt EvansDE@cf.ac.uk}}\\ \\
 {\sc Terry  Gannon }\\
 {\footnotesize Department of Mathematics, University of Alberta,}\\
{\footnotesize Edmonton, Alberta, Canada T6G 2G1}\\
{\footnotesize e-mail: {\tt tgannon@math.ualberta.ca}} }

\begin{document}
\maketitle

\begin{abstract} 
A near-group fusion category is a fusion category $\cC$ where all but 1 simple objects
are invertible. Examples of these  include the 
Tambara-Yamagami categories and the even sectors of the $D_5^{(1)}$ and $E_6$ subfactors, though there are infinitely many others.
We classify the near-group fusion categories,  and compute their doubles and the modular data
relevant to conformal field theory.
Among other things, we explicitly construct over 40 new finite depth subfactors, with Jones index ranging
from  around 6.85 to around 14.93. We expect all of these doubles to be realised by rational conformal
field theories.
\end{abstract}

{\footnotesize
\tableofcontents
}

\section{Introduction}

Considerable effort in recent years has been directed at the classification of subfactors of
small index.
Subfactors of index $\le 5$ are now all known (see e.g. \cite{MoSn}). 
The classification for index $\le 4$ was established some time ago. 
The Haagerup subfactor \cite{Haag} with Jones index $(5+\sqrt{13})/2\approx 4.30278$, 
the {Asaeda-Haagerup} subfactor \cite{ahaag} with index $(5+\sqrt{17})/2\approx 
 4.56155$, 
 and the 
{extended Haagerup} subfactor \cite{BMPS} with index  $\approx 4.37720$,
arose in Haagerup's classification {\cite{Haag}} of irreducible finite depth 
subfactors of index between 4 and 
$3+\sqrt{3}\approx  4.73205$. A Goodman-de la Harpe-Jones subfactor \cite{GHJ},
coming from the even sectors of the subfactor corresponding to the
$A_{1,10}\subset C_{2,1}$ conformal embedding,
has index $3+\sqrt{3}$. Then comes the Izumi-Xu subfactor $2221$ \cite{iz3} with index $(5+\sqrt{21})/2\approx 4.79$ and principal graph in Figure 1,  coming from the $G_{2,3}\subset
E_{6,1}$ conformal embedding.

The punchline is that, at least for small index, there are unexpectedly few subfactors.
Does this continue with higher index?
Are the aforementioned subfactors exotic, or can we put them into sequences?
In \cite{iz3}, Izumi realised the Haagerup and Izumi-Xu subfactors using endomorphisms in
Cuntz algebras, and suggested that his construction may generalise.
More precisely, to any abelian group $G$ of odd order, Izumi wrote down a nonlinear system of  
equations; any solution to them corresponds to  a subfactor of index $(|G|+2+\sqrt{|G|^2+4})/2$.
He showed the Haagerup subfactor corresponds to $G=\bbZ_3$, and that there also is a
solution for $G=\bbZ_5$. In \cite{EGh} we found solutions for the next several $G$, explained
that the number of these depends on the prime decomposition of $|G|^2+4$, and argued
that the Haagerup subfactor belongs to an infinite sequence of subfactors and so should not
be regarded as exotic.

Izumi in \cite{iz3} also  associated a second nonlinear system of equations to each finite
abelian group; to any solution of this system he constructs a subfactor of index
$(|G|+2+\sqrt{|G|^2+4|G|})/2$ and with principal graph $2^{|G|}1$, i.e. a star with 
one edge of length 1 and $|G|$ edges of length 2 radiating from the central vertex
 (Figure 1 is an example).
Izumi then found solutions for $G=\bbZ_n$ ($n\le 5$) and $\bbZ_2\times\bbZ_2$.
$G=\bbZ_1$ and $\bbZ_2$ correspond to the index $<4$ subfactors $A_4$ and $E_6$, respectively;
his solution for $\bbZ_3$ provides his construction for Izumi-Xu.
An alternate construction of $2221$, involving the conformal embedding $G_{2,3}\subset
E_{6,1}$, is due to Feng Xu as described 
in the appendix to \cite{CMS} (see also \cite{Hn}). As we touch on later in the paper,
there may be a relation between the series containing the Haagerup subfactor, and 
that containing the Izumi-Xu subfactor.

\medskip\epsfysize=1.3in\centerline{ \epsffile{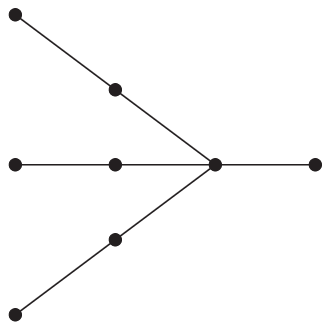}}\medskip
\centerline{Figure 1. The $2^31$ principal graph}
\label{fig1}\medskip

One of our tasks in this paper is to construct several more solutions to Izumi's second family
of equations, strongly suggesting that this family also contains infinitely many subfactors.
But more important, in this paper we study a broad class of systems of endomorphisms,
the \textit{near-group fusion categories},
including the Izumi-Xu series as a special case. We obtain a system of equations, generalising
those of Izumi, providing necessary and sufficient
conditions for their existence. We identify the complete list of solutions to the first several
of these systems, which permits us the construction of over 40 new finite-depth subfactors 
of index $<15$, 

A \textit{fusion category} $\cC$  \cite{ENO} is a $\bbC$-linear semisimple rigid monoidal category 
with finitely many simple objects and finite-dimensional spaces of morphisms, such that the 
endomorphism algebra of the neutral object is $\bbC$. The Grothendieck ring of a fusion
category is called a fusion ring. Perhaps the simplest examples
are associated to a finite group $G$: the objects are  $G$-graded vector spaces $\oplus_g
V_g$, with monoidal product $V_g\otimes V'_h=(V\otimes V')_{gh}$. Its fusion ring
is the group ring $\bbZ G$. We call
such examples \textit{group categories}. The category Mod$(G)$ of finite-dimensional 
$G$-modules is also a fusion category.

We're actually interested in certain concrete realisations of fusion categories, which we call
\textit{fusion $C^*$-categories}: the objects are  endomorphisms (or rather \textit{sectors}, 
i.e. equivalence classes of endomorphisms under the adjoint action of unitaries) on some
infinite factor $M$, the spaces Hom$(\rho,\sigma)$ of morphisms are intertwiners, and
the product is composition. Two fusion C$^*$-categories are equivalent iff they are
equivalent as fusion categories --- all that matters for us is that the factor $M$ exists,
not which one it is.
 Every finite-depth subfactor $N\subset M$ gives rise to two of these, one
corresponding to the principal, or $N$-$N$, sectors and the other to the dual principal, or $M$-$M$, ones. 
For example, given an outer action $\alpha$ of a finite group $G$ on an infinite
factor $N$, we get a  subfactor $N\subset N\sdprod G=M$ coming from the 
crossed product construction: the $N$-$N$ system
realises the group category for $G$, while the $M$-$M$ system realises Mod$(G)$.
Not all fusion categories can be realised as fusion C$^*$-categories
(e.g. the modular tensor categories associated to the so-called nonunitary Virasoro minimal models
are not fusion C$^*$-categories).

Perhaps the simplest nontrivial example of the extension of a fusion category is when the
category $\cC$ has precisely 1 more simple object than the subcategory $\cC_0$, and the latter  corresponds  to a finite abelian group. More precisely,
 simple objects $[g]$ in $\cC_0$ correspond to group elements $g\in G$, with tensor product 
 $[g][h]=[gh]$ corresponding to group multiplication. The simple objects of $\cC$ consist of
 the $[g]$, together with some object we'll denote $[\rho]$. Then $[\rho]$ must be self-conjugate,
 $[g][\rho]=[\rho]=[\rho][g]$, and $[\rho]^2=n'[\rho]+\sum_{g\in G}[g]$ (the multiplicities $n''_g$
 in the second term must be independent of $g$ because of equivariance $[g][\rho]=[\rho]$;
 because $[\rho]$ is its own
 conjugate, the multiplicity of $[1]$ must be 1).
 
 We call these \textit{near-group categories of type $G+n'$}. In this paper we restrict to abelian $G$,
 and we reserve $n$ always for
 the order of $G$. Examples of these have been studied in the literature:
 
 \begin{itemize}
 
 \item{} the Ising model and the module category of the dihedral group $D_4$, which are
 of type $\bbZ_2+0$ and $\bbZ_2\times\bbZ_2+0$, respectively;
 
\item{} more generally, the Tambara-Yamagami systems are by definition those of type $G+0$;

\item{} the $A_4,E_6$ and Izumi-Xu subfactors are of type $G+n$ for $G=\bbZ_1,\bbZ_2,\bbZ_3$ respectively;

\item{} more generally, Izumi's second hypothetical  family would be of type  $G+n$;

\item{} the $D_5^{(1)}$ subfactor and the module category of groups $S_3$ and $A_4$ are of type 
$\bbZ_2+1$, $\bbZ_2+1$, and $\bbZ_3+2$, respectively;

\item{} more generally, the representation category
of the affine group Aff$_1(\bbF_q)\cong \bbF_q\sdprod\bbF_q^\times$ of a finite field 
$\bbF_q$ is of type $\bbZ_{q-1}+(q-2)$.

\end{itemize}

More precisely, Aff$_1(\bbF_q)$ is the group of all affine maps $x\mapsto ax+b$ where $a\in\bbF^*_q$
and $b\in \bbF_q$. It has precisely $q=p^k$ conjugacy classes, with representatives 
$(a,0)$ and (1,1). 
It has precisely $q-1$ 1-dimensional representations, corresponding to the characters of $\bbF^*_q$. 
The remaining 
irrep is thus of dimension $\sqrt{q(q-1)-q-1}=q-1$, and is the nontrivial summand of
the natural permutation representation of  Aff$_1(\bbF_q$) on $\bbF_q$ given by the affine maps: 
$(a,b).x=ax+b$. 

In this paper we classify the near-group C$^*$-categories $G+n'$, in the sense that we obtain
polynomial equations in finitely many variables, whose solutions correspond bijectively
to equivalence classes of the near-group C$^*$-categories.  
Given any near-group C$^*$-category $\cC$ with $n'>0$, we identify a natural subfactor 
$\rho(M)\subset M$ whose even systems 
are both identified with $\cC$. We also work out the principal graph of the closely related
subfactor $\rho(M)\subset M^G$. By contrast, we can realise some but not all $\cC$ with $n'=0$, as the even sectors of  a subfactor.

There is a fundamental dichotomy here: $n'$ either equals $n-1$,
or is a multiple of $n$, where as always $n=|G|$. When $n'<n$, we have a complete classification:

\medskip\noindent\textbf{Fact.} \textit{Let $G$ be any abelian group of order $n$.}

\smallskip\noindent\textbf{(a)} \textit{There are precisely two $C^*$-categories of type
$G+0$.}

\smallskip\noindent\textbf{(b)} \textit{When $n'$ is not a multiple of $n=|G|$, 
the only  $C^*$-categories of type $G+n'$ are Mod(Aff$_1(\bbF_{n+1})$), except for $n=1,2,3,7$
which have 1,2,1,1 additional $C^*$-categories. In all cases here, $n+1$ is a prime power, 
$n'=n-1$, and $G=\bbZ_{n+1}$.}\medskip

This is our Corollary 4 and Proposition 5 respectively, proven below.
Type $G+0$ and type $\bbZ_n+n-1$ fusion categories were classified by Tambara-Yamagami
\cite{TY} and Etingof-Gelaki-Ostrik \cite{EGO}, respectively; we find that for these types, all
fusion categories can be realised as C$^*$-categories. Our proof of (a) is independent
of and much simpler than \cite{TY}. 

\medskip\noindent\textbf{Conjecture 1.} \textit{For every nontrivial cyclic group $G=\bbZ_n$, there 
are at least 2 inequivalent subfactors with principal graph $2^{n}1$whose principal even sectors satisfy the near-group fusions of type $G+n$.}\medskip

We have verified this for $n\le 13$. For those $n$, the complete classification is given in
Table 2 below. 
In the process, we construct dozens of new finite depth subfactors of small index with principal even sectors of near-group type.
This classification for $n=3$ yields a uniqueness proof
(up to complex conjugation) for the principal even sectors of the Izumi-Xu 2221 subfactor; this can be
compared to Han's uniqueness proof \cite{Hn} of the 2221 subfactor. Again, our proof is independent of and both considerably
shorter and simpler than the original one. 
We do not yet feel confident speculating on systems with $n'>n$; the corresponding subfactors
would have principal graph as in Figure 3 below.

Two morals can be drawn from this paper together with our previous one \cite{EGh}. 
One is that there is surely a plethora of undiscovered finite-depth subfactors,
of relatively small index. This is in marked contrast to the observations of e.g.
\cite{MoPe}, who speak of the `little desert' in the interval $5<[M:N]<3+\sqrt{5}$. The situation
here is probably very analogous to the classification of finite groups, which also is very
tame for small orders. The second moral is that, when the fusions are close to that of a group,
a very promising approach to the classification and construction of corresponding systems of endomorphisms, equivalently C$^*$-categories, or the corresponding subfactors, is the Cuntz algebra
method developed in e.g. \cite{iz3} and championed here. This approach  also makes the computation
of the tube algebra and corresponding modular data etc (to be discussed shortly) completely accessible. In contrast,
 the technique of planar algebras is more robust, able to handle subfactors unrelated to groups,
 such as Asaeda-Haagerup
 and the extended Haagerup. But planar algebra techniques applied to e.g. the Haagerup
 fail to see that it (surely) lies in an infinite family. In a few minutes  the
 interested computer can construct several more subfactors of the type described in Conjecture 1,
 using the Cuntz algebra method here, each of which would be a serious challenge for
 the planar algebra method.

The underlying presence of groups here begs the question of $K$-theory realisations of these
fusion rings. For example, the fusion ring of the near-group C$^*$-categories when $n'=n-1$
can be expressed as $K_0^{\mathrm{Aff}_1(\bbF_q)}(1)$. Is there a $K$-theoretic expression in the
other class, i.e. when $n'\in n\bbZ$?

An important class of fusion categories are the \textit{modular tensor categories} \cite{tur}, which are
among other things braided and carry a representation (called \textit{modular data})
of the modular group SL$_2(\bbZ)$ of dimension equal to the rank of the category,
from which e.g. the fusion coefficients can be computed. These
arise from braided systems of endomorphisms on an infinite factor, from representations of
completely rational conformal nets, or from the modules of
a rational vertex operator algebra.
 
There is a standard construction, called the \textit{quantum} or \textit{Drinfeld double}, to
 go from fusion categories (with mild additional properties) to modular tensor categories \cite{m1}. We construct the doubles of our
 C$*$-categories, following the tube algebra approach \cite{iz1}, 
 and in particular explicitly compute its modular data. As with the Haagerup series,
 our formulae are unexpectedly simple. This simplicity also challenges the perceived exoticness
 of these subfactors.
 
 A natural question is,
 are these modular tensor categories realised by conformal nets of factors, or by rational vertex
 operator algebras (VOAs)?
 Ostrik (see Appendix A in \cite{CMS}) shows that the double of Izumi-Xu $2221$
has a VOA interpretation, in fact it is the affine algebra VOA corresponding to
$G_{2,3}\oplus A_{2,1}$. No construction is
known for the large $G$. Curiously, this VOA conformally embeds into that of 
 $E_{6,1}\oplus A_{2,1}$ (which realises the fusions
of the double of $\bbZ_3$), and this was where \cite{EGh} suggests to look for the
VOA associated to the double of the Haagerup. Could there be a relation  between Izumi-Xu $2221$ and
the Haagerup? Other reasons suggest a relationship between $2^91$ and the Haagerup.
We discuss this latter possibility briefly in Subsections 3.5 and 4.4.
 
 More generally we could consider quadratic extensions of a group category. More precisely, 
 let $G$ be a finite group (not necessarily abelian) and suppose $[\rho]^*=[\rho][g_\rho]$
for some $g_\rho\in G$. Let $N$ be any subgroup of $G$: we require $[g][\rho]=[\rho]$ iff
$g\in N$. Then $[\rho][g]=[\rho]$ iff $g\in g_\rho Ng_\rho^{-1}=:N'$. 
The simple objects in this category are $[g]$ for $g\in G$ as well as $[g_i][\rho]$ for 
representatives $g_i$ of cosets $G/N$. Let $\phi$ be any isomorphism 
$G/N\rightarrow G/N'$; we require
$[g][\rho]=[\rho][g']$ iff $g'\in\phi(gN)$. Then $[\rho]^2 =\sum_{g\in N}[g]+\sum_{i}n'_i[g_i][\rho]$.
We require $\phi$ to satisfy $g_\rho^{-1}\phi(\phi(g))g_\rho=g$ for all $g\in G$.
This large class of examples should be accessible to a similar treatment. The near-group
categories correspond to the choice $N=G$ and $g_\rho=1$; the Haagerup-Izumi series
\cite{iz3,EGh}  corresponds to $G=\bbZ_{2n+1}$, $N=1$, $\phi(g)=-g$, $n'_i=1$; in particular,
 the Haagerup subfactor at index $(5+\sqrt{13})/2$ corresponds to $G=\bbZ_3$.
It would be very interesting to extend the analysis in this paper to this larger class. 

Here is a summary of our main results. Theorem 1 associates numerical invariants 
to a near-group C$^*$-category, which according to  Corollary 1 completely characterise the
category. Corollary 2 (and the end of Subsection 2.2) associate to each C$^*$-category two
subfactors and work out their principal graphs. Theorem 2 establishes the fundamental
dichotomy of near-group C$^*$-categories: either $n'=n-1$ or $n'\in n\bbZ$. When $n'=n-1$,
Theorem 3 lists the identities necessarily obeyed by the numerical invariants and shows
they are also sufficient. Theorem 4 does the same when $n|n'$. In Proposition 5 we find all
 near-group C$^¥$-categories with $n'=n-1$; we see that almost all of these are 
known. In Table 2 we list the first several with $n'=n$, and find that almost none of these are
known. In Theorem 5 and Corollary 6 we work out the tube algebra and modular data
for any near-group C$^*$-category with $n'=n-1$. \cite{iz3} had found a very complicated expression
for the modular data when $n'=n$; we notice in Subsection 4.4 that it collapses to cosines.

\medskip\noindent\textit{Note added in proof.} After completing this manuscript, we received
in July 2012 \cite{iz5} from Masaki Izumi, which overlaps somewhat the contents of our paper.
In particular, he also obtained necessary and sufficient conditions for the Cuntz
algebra construction to realise a near-group C$^*$-category of type $G+n'$. On the one hand, 
unlike us, he does not require $G$ to be abelian, and he allows the possibilities of an $H^2$-twist.
On the other hand, unlike us, he does not address principal graphs of associated
subfactors, nor the tube algebra, nor the modular data (simplified or otherwise) for
the doubles, and he does not construct new solutions of the resulting equations and
hence does not construct new subfactors.


\section{The near-group systems}

\subsection{The numerical invariants}

Let $G$ be a finite abelian group (written additively) with order $n=|G|$,  and as usual write 
$\widehat{G}$ for its irreps. Let $M$ be an infinite factor,
$\rho$ a self-conjugate irreducible endomorphism on $M$ with finite statistical  dimension $d_\rho<\infty$, 
and $\alpha$  an outer action of $G$ on $M$. Suppose the following fusion rules hold:
\begin{align}
&&[\alpha_g\rho]=[\rho]=[\rho\alpha_g]\,,\\
&&[\rho^2]=\bigoplus_g[\alpha_g]\oplus n'[\rho]\,,\label{fusrho2}\end{align}
for some $n'\in\bbZ_{\ge 0}$. Then the $d_\rho$ satisfies $d_\rho^2=n'd_\rho+n$ so
\begin{equation}d_\rho=\frac{n'+\sqrt{n'{}^2+4n}}{2}=:\delta\,.\label{delta}\end{equation}
 Let $\cC(G,\alpha,\rho)$ denote the fusion C$^*$-category generated by $\alpha,\rho$.
 We call these, \textit{$C^*$-categories of type $G+n'$}.
 
\medskip\noindent\textbf{Definition 1.} \textit{By a} {pairing} $\langle g,h\rangle$ \textit{on $G$ we mean a complex-valued function on
 $G\times G$ such that for all fixed $g\in G$, both $\langle g,*\rangle,\langle *,g\rangle\in\widehat{G}$.
 By a} {symmetric pairing} \textit{we mean a pairing satisfying $\langle g,h\rangle=\langle h,g\rangle$.
 By a} {nondegenerate pairing} \textit{we mean a pairing for which the characters $\langle g,*\rangle$
 are distinct for all $g$.}\medskip
 
  Note that a nondegenerate pairing is equivalent to a choice of
 group isomorphism $G\rightarrow\widehat{G}$, $g\mapsto\phi_g$, by $\phi_g(h)=\langle g,h\rangle$. The  nondegenerate symmetric pairings for $G=\bbZ_n$ are
  $\langle g,h\rangle=e^{2\pi\i mgh/n}$ for some integer $m$ coprime to $n$.

\medskip\noindent\textbf{Theorem 1.} \textit{Let $G,\alpha,\rho$ be a $C^*$-category of type $G+n'$. 
Suppose in addition
that $H^2(G;\bbT)=1$. Then there are $n+n'$ isometries $S_g,T_z$ ($g\in G$, $z\in\cF$) 
satisfying the
Cuntz relations, such that $\alpha_g\rho=\rho$, $\rho\alpha_g=\mathrm{Ad}(U_g)\,\rho$, for
a unitary representation $U_g$ of $G$ of the form
\begin{equation}
U_g=\sum_h\langle g,h\rangle\,S_hS_h^*+\sum_{z}u_{z,g}T_zT_{gz}^*\,,\label{Ug}\end{equation} 
where  $G$ permutes the $z\in\cF$ and $\langle g,h\rangle$ is a pairing on $G$.
 Moreover, $\alpha_g(S_h)=S_{g+h}$ and $\alpha_g(T_z)=\widetilde{z}(g)T_z$ for 
some $\widetilde{z}\in\widehat{G}$. Finally,
\begin{align}
&&\rho S_g=\left({s}{\delta}^{-1}\sum_h\langle g,h\rangle\, S_h+\sum_{x,z}u_{x,g}a_{gx,z}T_xT_z\right) U_g^*\,,\label{rhoS}\\
&&\rho(T_z)=\sum_{h,x}\overline{\widetilde{x}(h)}\,b_{z,x}S_hT_x^*+\sum_{x,h}\widetilde{x}(h)\,b'_{z,x}T_xS_hS_h^*+\sum_{w,x,y}b''_{z;w,x,y}
T_wT_xT^*_y\,,\label{rhoT}\end{align}
for some sign $s\in\{\pm 1\}$ and complex parameters $a_{y,z},b_{z,x},b'_{z,x},b''_{z;w,x,y}$, for
$w,x,y,z\in\cF$. }

\medskip\noindent\textit{Proof.} 
 Our argument follows in part that of the first theorem of \cite{iz4}. 
Because $[\alpha_g\rho]=[\rho]$, there exists a unitary $W_g\in
U(M)$ for each $g\in G$, satisfying $\alpha_g\rho=\mathrm{Ad}(W_g)\,\rho$. But
\begin{equation}
\mathrm{Ad}(W_{g+h})\,\rho=\alpha_h\alpha_g\rho=\mathrm{Ad}(\alpha_h(W_g)\,W_h)\,\rho\end{equation}
for all $g,h\in G$, so $\alpha_h(W_g)\,W_h=\xi(g,h)\,W_{g+h}$ for some 2-cocycle
$\xi\in Z^2(G;\bbT)$. Because $H^2(G;\bbT)=1$, we can require that $\xi$ be identically 1, by
tensoring  $W_g$ with the appropriate  1-coboundary. Since $G$ is a finite group and
$\alpha$ is outer, the $\alpha$-cocycle $W_g$ is a coboundary, so there exists a unitary
$V\in U(M)$ so that $W_g=\alpha_g(V^*)\,v$ for all $g\in G$. This means $\mathrm{Ad}(\alpha_g(
V))\,\alpha_g(\rho)=\mathrm{Ad}(V)\,\rho$, i.e. $\alpha_g(\mathrm{Ad}(V)\,\rho)=\mathrm{Ad}(V)
\,\rho$. Thus if we replace $\rho$ by $\mathrm{Ad}(V)\,\rho$ we obtain $\alpha_g\rho=\rho$
as endomorphisms, not just as sectors.

This has exhausted most of the freedom in choosing $\rho$. The fusion $[\rho\alpha_g]=[\rho]$ 
means $\rho\alpha_g=\mathrm{Ad}(U_g)\,\rho$ for some unitaries $U_g$; because $H^2(G;\bbT)=1$, we 
can in addition insist that $g\mapsto U_g$ defines a unitary representation of $G$.
Note that we still have a freedom in replacing $U_g$ with $\psi(g)\,U_g$ for any character
$\psi\in\widehat{G}$.

The fusion \eqref{fusrho2} means 
\begin{equation}
\rho^2(x)=\sum_{g\in G}S_g\alpha_g(x)S_g^*+\sum_{z\in\cF}T_z\rho(x)T_z^*\,,\label{endrho2}\end{equation}
where $S_g$ and $\{T_z\}_{z\in\cF}$ are bases of isometries for the intertwiner spaces Hom$(\alpha_g,\rho^2)$ 
and Hom$(\rho,\rho^2)$ respectively (so $\rho^2(x)S_g=S_g\alpha_g(x)$ etc). 
Then \eqref{endrho2} implies $S_g,T_z$ obey the Cuntz relations. 
Since $\mathrm{Ad}(U_g)\,\rho^2=\rho\alpha_g\rho=\rho^2$,
$U_h$ maps Hom$(\alpha_g,\rho^2)$ to itself and Hom$(\rho,\rho^2)$ to itself, i.e. $U_hS_g=
\mu_g(h)\,S_g$ and $U_hT_z=\sum_wu(h)_{z,w}T_w$ for some
$\mu_g(h),u(h)_{z,w}\in\bbC$. Since $U_{h+h'}=U_hU_{h'}$, we have that $\mu_g\in\widehat{G}$ 
for each $g\in G$, and the matrices $u$ define a unitary 
representation on Hom$(\rho,\rho^2)$. This gives us 
\begin{equation}U_g=\sum_h\mu_h(g)S_hS_h^*+\sum_{z,y}u(g)_{z,y}T_zT_{y}^*\,.\label{Ugen}\end{equation} 
Define $U_g'=\overline{\mu_0(g)}\,U_g$. Then $\rho\,\alpha_g=\mathrm{Ad}(U_g')\,\rho$ and
$U_g'$ is still a unitary representation of $G$. For this reason we may assume that
$\mu_0$ is identically 1.

Similarly, $\alpha_h$ maps Hom$(\alpha_g,\rho^2)$ to Hom$(\alpha_{g+h},\rho^2)$ and 
Hom$(\rho,\rho^2)$ to itself,  as $\alpha_h\rho=\rho$.
This means $\alpha_h(S_g)=\psi_{g,h}S_{g+h}$  for some nonzero $\psi_{g,h}\in\bbC$, 
and $\alpha_h$ defines an $n'$-dimensional unitary $G$-representation
on Hom$(\rho,\rho^2)$. Because
$H^2(G;\bbT)=1$ we can choose $\psi_{g,h}$ to be identically 1. Because $G$ is abelian, 
we can diagonalise the $n'$-dimensional representation, i.e. choose our basis $T_z$ so that
 $\alpha_gT_z=\widetilde{z}(g)T_z$ for some $\widetilde{z}\in\widehat{G}$.
 
Because $\rho$ is self-conjugate and $S_0\in\mathrm{Hom}(\mathrm{id},\rho^2)$, the isometry $S_0$ will
satisfy $S_0^*\rho (S_0)=s/\delta$ for a sign $s$. Hence $S_g^*\rho(S_0)=\alpha_g(S_0^*\rho(S_0))
=s/\delta$.

For any $T\in\mathrm{Hom}(\rho,\rho^2)$, define the right and left Frobenius maps $\cR(T)=\sqrt{\delta}T^*\rho(S_0)$
and $\cL(T)=\sqrt{\delta}\rho(T^*)S_0$, as in Section 3.2 of \cite{BEK1}. Then  $\rho^2(x)\cR(T)=T^*\rho^3(x)\rho(S_0)
=T^*\rho(\rho^2(x)S_0)=T^*\rho(S_0x)=\cR(T)\rho(x)$ and $\rho^2(x)\cL(T)=\rho(T^*\rho^2(x))S_0
=\cL(T)\rho(x)$ so both $\cL,\cR$ are conjugate-linear on the space Hom$(\rho,\rho^2)$.
$\cR$ is surjective: $\cR(\cR(T))=\delta\rho(S_0)^*T\rho(S_0)=
\delta\rho(S_0^*\rho S_0)T=sT$. A similar calculation shows that for any $T,T',T''\in\mathrm{Hom}(\rho,\rho^2)$,
$T^*\rho(T')T''\in\mathrm{Hom}(\rho,\rho^2)$: $(T^*\rho(T')T'')\rho(x)=T^*\rho(T'\rho(x))T''=\rho^2(x)T^*\rho(T')T''$.

Since $1=\sum_gS_gS^*_g+\sum_zT_zT^*_z$, we find
\begin{align}
&&\rho(S_0)=\sum S_gS_g^*\rho(S_0)+
\sum T_zT_z^*\rho(S_0)=s\delta^{-1}\sum S_g+\sum T_z\cR(T_z)\nonumber\\&&
={s}{\delta}^{-1}\sum_hS_h+\sum_{z,y}a_{z,y}T_zT_y\label{rhos0}
\end{align}
for some complex numbers $a_{z,y}$, and covariance $\rho(S_g)=\rho(\alpha_g(S_0))=\mathrm{Ad}(U_g)\,\rho(S_0)$ forces \eqref{rhoS}.

We can identify the shape of $\rho(T)$ similarly. Choose some $T_z\in\mathrm{Hom}(\rho,\rho^2)$; then surjectivity of
$\cR$ implies there is some $T_z'\in\mathrm{Hom}(\rho,\rho^2)$ such that $T_z=\cR(T_z')$. We find
\begin{align}
&&\rho(T_z)=\sum_gS_gS_g^*\rho(T_z)+\sqrt{\delta}\sum_{w,g}T_wT^*_w\rho(T'{}^*_z)\rho^2(S_0)S_gS_g^*+\sum_{w,x}T_wT^*_w\rho(T_z)
T_xT^*_x\nonumber\\&&
=\sum_gS_g\alpha_g(\cL(T_z)^*)+\sum_{w,g}T_wT^*_w\alpha_g(\cL(T_z))\alpha_g(S_0)S_g^*
+\sum_{w,x}T_w(T^*_w\rho(T_z)T_x)T^*_x\nonumber\\&&
=\sum_{g,x}b_z(g,x)\,S_gT_x+\sum_{g,x}b'_z(x,h)\,T_xS_gS_g^*
+\sum_{w,x}T_w(T^*_w\rho(T_z)T_x)T^*_x\,.\nonumber\end{align}
$\alpha_g\rho(T_z)=\rho(T_z)$ forces $b_z(h,x)=b_{z,x}\overline{\widetilde{x}(h)}$
and $b'_z(x,h)=b'_{z,x}\widetilde{x}(h)$, which gives \eqref{rhoT}.

All that remains is to show that a basis $T_z$ of Hom$(\rho,\rho^2)$ can be found for which
$u(g)$ in \eqref{Ugen} is a generalised permutation matrix. For each $\phi\in\widehat{G}$ let
 $\cT_\phi$ denote the (possibly empty)  subspace of Hom$(\rho,\rho^2)$ on which $\alpha_g$ acts as $\phi$,
 so Hom$(\rho,\rho^2)=\oplus_\phi\cT_\phi$.
Note that $\alpha_g\rho(S_h)=\rho(S_h)$ implies (among other things) the selection rule:
 $u(g)_{x,y}\ne 0\Rightarrow\widetilde{y}={\widetilde{x}}\mu^g$ for 
$\mu^g\in\widehat{G}$ defined by $\mu^g(k)=\mu_k(g)$. This means there is a pairing
$\langle g,h\rangle$ on $G$ such that $\mu_h(g)=\mu^g(h)=\langle g,h\rangle$. Each
$u(g)$ defines a linear isomorphism from $\cT_\phi$ to $\cT_{\mu^g\phi}$ (with inverse $u(-g)$).  

Define $H$ to be the set of all $h\in G$ such that $\mu^h=1$; then $H$ is a subgroup of $G$.
Choose a set $O$ of orbit representatives  of this $G$-action $\phi\mapsto\mu^g\phi$ on $\widehat{G}$,
and a set $C$ of coset representatives  for $G/H$. Note that  $u$ restricts to a unitary representation of $H$ on each space $\cT_\phi$; for each representative $\phi\in O$
 choose a basis $\cF_{\phi}$ of $\cT_{\phi}$ diagonalising this $H$-representation.
For any $\psi\in \widehat{G}$ there will be a unique choice of representatives $\phi^\psi\in O,k_\psi\in C$ such that
$\mu^{k_\psi}\phi^\psi=\psi$; let the basis $\cF_\psi$ on $\cT_\psi$ be the image under $u(k_\psi)$ of the
basis for $\cT_{\phi^\psi}$. Our basis $\cF=\{T_z\}$ of the intertwiner space Hom$(\rho,\rho^2)$ 
will be the union of these bases $\cF_\psi$ for the subspaces 
$\cT_\psi$. For any basis vector $T_x\in\cF$ and $g\in G$, write $T_x=
u(k_x)T_y$ for some representative $k_x\in C$ and basis vector $T_y\in\cF_\phi$ for $\phi\in O$,
and $g+k_x=h+k'$ for $k'\in C$ and $h\in H$. Then because $G$ is abelian
we have $u(g)_{xx'}=u(g+k_x)_{yx'}=u(h)_{yy}\delta_{T_{x'},u(k')T_y}$. Now write $u_{x,g}=u_{y,h}$
and define  $gx$ to label the basis vector $T_{gx}=u(k')T_z$. Then $gx$ defines a $G$-action
on $\cF$, and $u(h)$ is diagonal $\forall h\in H$ for this basis, and
 the matrix entries $u(g)_{x,y}$ equal $u_{x,g}\delta_{y,gx}$ as desired.
  \qquad\textit{QED to Theorem 1}\medskip

\noindent\textbf{Corollary 1.} \textit{Let $G$ be a finite abelian group with $H^2(G;\bbT)=1$
and choose any $C^*$-category $\cC$ of type $G+n'$. Then
the sign $s$, pairing $\langle*,*\rangle$, and complex numbers
$u_{x,g},a_{x,y},b_{z,x},b'_{z,x},b''_{z;w,x,y}$ form a complete invariant of $\cC$,
up to gauge equivalence and  automorphism of $G$.} \medskip

By \textit{gauge equivalence} we mean equivalence up to a change of basis on Hom$(\rho,\rho^2)$ (in contrast, relative rescaling of the $S_h$ would wreck \eqref{rhos0} so isn't allowed). 
More precisely, for each representative $\phi\in O$ and $\psi\in\widehat{H}$ 
let $\cT_{\phi,\psi}$ denote the subspace of Hom$(\rho,\rho^2)$ spanned by $T_x\in\cF$
with $\widetilde{x}(g)=\phi(g)$ for $g\in G$, and $u_{x,h}=\psi(h)$ for $h\in H$. Then gauge
equivalence amounts to a change-of-basis $P_{\phi,\psi}\in U(\cT_{\phi,\psi})$.
Let $P\in U(\mathrm{Hom}(\rho,\rho^2))$ be the direct sum of  $\|O\|$ copies of each 
$P_{\phi,\psi}$ and define $T^{old}_{x}=\sum_{w}P_{x,w}T_w^{new}$. Then
$a^{new}=P^Ta^{old}P$, $b^{new}={P}^{-1}b^{old}\overline{P}$, $b^{\prime\,new}=P^{-1}
b^{\prime\,old}P$, $b^{\prime\prime\,new}_{z;w,x,y}=\sum_{z',w',x',y'}b^{\prime\prime\,old}_{z';w',x',y'}\,\overline{P}_{z',z}\,P_{w',w}\,P_{x',x}\,\overline{P}_{y',y}$. An automorphism
$\phi$ of $G$ acts by permuting $S_h\mapsto S_{\phi h}$, $T_{x}\mapsto T_{\phi x}$ where
$\widetilde{\phi x}=\phi\widetilde{x}$, $\langle g,h\rangle\mapsto\langle \phi g,\phi h\rangle$,
etc.

The requirement that $H^2(G;\bbT)=1$  is made for 
simplicity; dropping it would introduce more possibilities (e.g. $U_g$ need only be a projective representation). Recall that for an abelian group, $H^2(G;\bbT)=1$ iff $G$ is cyclic.
We will see shortly that either $n'=n-1$ or $n'\in n\bbZ$; when $n'=n-1$ $G$ must be
cyclic and therefore in this case the hypothesis $H^2(G;\bbT)=1$ in Theorem 1 etc
is redundant and can be dropped.

\medskip\noindent\textbf{Corollary 2.} \textit{Suppose $\cC$ is a $C^*$-category of type
$G+n'$ for $n'>0$. Then $\rho(M)\subset M$ is a subfactor with index $d_\rho^2$, whose $M$-$M$ and
$N$-$N$ systems are both type $G+n'$, with (dual)  principal graph consisting of $n$ vertices attached to
a left vertex, $n$ other vertices attached to a right vertex, and the left and right vertices 
attached with $n'$ edges (see  Figure 2 for an example). This has the intermediate
subfactors $\rho(M)\subset M^G\subset M$ (where the $G$-action is given by $\alpha_g$)
and $\rho(M)\subset\rho(M)\sdprod G\subset
M$ (using the $G$-action $\rho\alpha_g\rho^{-1}$).}

\medskip\epsfysize=.8in\centerline{ \epsffile{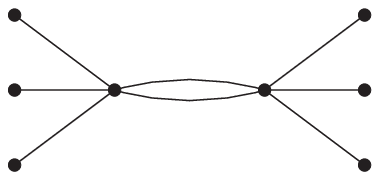}}\medskip

\centerline{Figure 2. The principal graph for $\rho_{3,2}(M)\subset M$}
\label{fig2}\medskip

The subfactor $\rho(M)\subset M$ is self-dual because $\rho=\overline{\rho}$;
the principal graphs for $\rho(M)\subset M^G$ match those of $\rho(M)\sdprod G\subset
M$ because of the basic construction applied to $\rho(M)\subset\rho(M)\sdprod G\subset
M$. At the end of next subsection we
compute the principal graph of  $\rho(M)\sdprod G\subset
M$ in all cases.

$G$ acts on $M$ through the $\alpha_g$; the relation $\alpha_g\cdot\rho= \rho$ says
$\rho(M)\subset M^G$. When $n'=0$, we see $\rho(M)=M^G$ by an index calculation.
The index of $\rho(M)\subset M^G$ is $1+n'n^{-1}\delta$.

Recall the near-group C$^*$-categories listed in the Introduction. For example:

\begin{itemize}

\item{} The Tambara-Yamagami categories \cite{TY}, 
which are of type $G+0$, correspond to $\delta=\sqrt{n}$, $\mu_h(g)
=\langle g,h\rangle$, and $s=\pm 1$.

\item{} The hypothetical Izumi-Xu family of subfactors (Section 5 of \cite{iz3}) correspond to the
parameter choices $s=1$, $\cF=\widehat{G}$, $n'=n$, $\delta=(n+\sqrt{n^2+4n})/2$, $\mu_h(g)=\langle g,h\rangle=\mu^h(g)$, $a_{x,y}
=\sqrt{\delta}^{-1}\delta_{x,\overline{y}}\,a(x)$, $u_{x,g}=1$, $b_{z,x}=c\sqrt{n\delta}^{-1}
\langle z,x\rangle$, $b'_{z,x}=\overline{a(z)c}\sqrt{n}^{-1}\langle z,x\rangle$,
$b''_{z;w,x,y}=\delta_{y,wx}\,a(x)\,b(z\overline{x})\langle z,y\rangle$ for some complex numbers
$c,a(x),b(x)$ where $a(z)^2= \overline{\langle z,z\rangle}$ and $\langle g,h\rangle$
is a nondegenerate symmetric pairing on $G$. We extend this pairing to $\widehat{G}$ 
through the group isomorphism $G\rightarrow\widehat{G}$: 
 $\langle x_g,x_h\rangle:=\langle g,h\rangle$.

\item{} The $D_5^{(1)}$ subfactor (Example 3.2 of \cite{iz0}) has $G=\bbZ_2$, $\cF=\widehat{G}\setminus\{1\}$, $n'=1$, $\delta=2$, $s=1$,
$\mu_h(g)=\mu^h(g)=1$, $u_{x,g}=(-1)^g$, $a_{x,x}=\sqrt{2}^{-1}$,
$b_{x,x}=\overline{a}\sqrt{2}^{-1}$, $b'_{x,x}=a$,
$b''_{x;x,x,x}=0$ for some complex number $a$ satisfying $a^3=1$.


\end{itemize}

\subsection{Generalities}

There is a fundamental bifurcation of the theory of near-group C$^*$-categories:

 \medskip\noindent\textbf{Theorem 2.} \textit{Suppose $\cC$ is a near-group $C^*$-category of type $G+n'$,
and suppose $H^2(G;\bbT)=1$.  Let $\langle *,*\rangle,u_x,s,a,b,b',b''$ be the parameters 
which Corollary 1 associates to $\cC$.} 

\smallskip\noindent\textbf{(a)} \textit{Either $n'=n-1$, or $n'=kn$ for some $k\in\bbZ_{\ge 0}$.}

\smallskip\noindent\textbf{(b)} \textit{Suppose $n'=n-1$. Then $\delta=n$, the pairing $\langle g,
h\rangle$ is identically 1, $gx=x$ $\forall x\in\cF$ and $g\in G$. The assignment $x\mapsto 
\widetilde{x}$ bijectively identifies $\cF$ with $\widehat{G}\setminus \{1\}=:\widehat{G}^*$.
There is a permutation $\sigma$ of $\widehat{G}^*$ such that $u_{x,g}=(\sigma(x))(g)$ for all
$x$ and $g$. Finally, $a_{x,y}=\sqrt{\delta}^{-1}\delta_{y,\overline{x}}$.}

\smallskip\noindent\textbf{(c)} \textit{Suppose $n'=kn$ for $k\in\bbZ_{\ge 0}$. Then the pairing
$\langle g,h\rangle$ is nondegenerate, and for any $x\in \cF$ there is a unique $g_x\in G$
such that $\widetilde{g_xx}=1$. Moreover, $u_{x,g}$ is identically 1, and $a_{x,y}\ne 0$
implies $\widetilde{x}\widetilde{y}=1$.}\medskip

\noindent\textit{Proof.}  Let $\cC=\cC(G,\alpha,\rho)$ be of type $G+n'$, and let $s,\ldots,
b''$ be its numerical invariants, and $\langle *,*\rangle$ its pairing. Define subgroups
$H,H'$ of $G$ by $H=\{h\in G\,|\,\langle h,g\rangle=1\ \forall g\}$ and $H'=\{h'\in G\,|\,\langle g
,h'\rangle=1\ \forall g\}$.  Let $n''=|H|$.  Write $\mu^g(h)=\langle g,h\rangle=\mu_h(g)$ as before.
Note that the orders $|H|$ and $|H'|$ must be equal, since the row-rank of the
matrix $\langle g,h\rangle$ will equal its column-rank.

Let us review some observations contained in the proof of Theorem 1.
Recall the coset representatives $k\in C$ and orbit representatives $\phi\in O$ introduced in
the proof of Theorem 1. We saw there that the phases $u_{x,h}$ restricted to $h\in H$ forms a
representation of $\widehat{H}$, which we'll denote by $\ddot{x}$. Then we found there the 
formula $u_{x,g}=\ddot{x}(g+k_x-k')$ valid for any $g\in G$ and $x\in\cF$, where
$k_x,k'\in C$ satisfy $\mu^{-k}\widetilde{x}\in O$ and $g+k_x-k'\in H$.
 Recall the partition $\cF=\cup\cF_{\phi,\psi}$ where $\phi\in\widehat{G},\psi\in
\widehat{H}$; the $G$-action $x\mapsto gx$ on $\cF$, contained in $U_g$, bijectively relates
$\cF_{\widetilde{x}, \ddot{x}}$ to $\cF_{\phi,\ddot{x}}$ where $\phi\in O$ is the unique
representative with $\phi|_H=\widetilde{x}|_H$. 

Recall the Cuntz algebra  $\cO_{n,n'}$ generated by $S_g$ and the $T_z$.
Being an endomorphism of $\cO_{n,n'}$, $\rho$ preserves the
Cuntz relations. Firstly, $\rho(S_g)^*\rho(S_h)=\delta_{g,h}$ for $g\in H$ is equivalent to 
\begin{equation}\delta_{g,0}=n\delta^{-2}\delta_{g+H,H}+\sum_{x,z}\ddot{x}(g)\,\overline{a_{x,z}}\,a_{gx,z}\,.\label{Cu:SS00}
\end{equation}
Putting $g=0$ gives $\sum_z|a_{x,z}|^2=\delta^{-1}$ for all $x$.
Hitting both sides with $\overline{\psi(g)}$ for any $\psi\in\widehat{H}$ and summing over $g\in H$,
\eqref{Cu:SS00} is equivalent to
\begin{equation}\label{Cu:SS}1-nn''\delta^{-2}\delta^{H}_{1,\psi}=\delta^{-1}n''\ddot{n}_\psi\,,
\end{equation}
where $\ddot{n}_\psi$ denotes the number of $x\in\cF$ with
 $\ddot{x}=\psi$. The  $T_xS_0S_0^*T_y$ and $S_0S_g^*$ coefficients of 
completeness $1=\sum_g\rho S_g\,\rho S_g^*+\sum_z\rho T_z\,\rho T_z^*$ give
unitarity of the matrix $b'$ together with
\begin{equation}
\delta_{g,0}={n}{\delta^{-2}}\delta_{g+H',H'}+\sum_{z,x}|b_{z,x}|^2\,{\widetilde{x}(g)}\,.\label{compleeta}\end{equation}
Putting $g=0$ into \eqref{compleeta} tells us $\sum_z|b_{z,x}|^2=\delta^{-1}$, so 
\eqref{compleeta} becomes
\begin{equation}
\label{compleeeta}1-nn''\delta^{-2}\delta^{H}_{1,\phi}=\delta^{-1}n\widetilde{n}_\phi
\end{equation}
for any $\phi\in\widehat{G}$, where $\widetilde{n}_\phi$ denotes the number of $x\in\cF$ with
 $\widetilde{x}=\phi$.
 
 Now suppose $H\ne 0$. Then there exist $\psi\in\widehat{H}$ and $\phi\in\widehat{G}$
 such that $\psi\ne 1$ and $\phi|_H\ne 1$, so \eqref{Cu:SS},\eqref{compleeeta} say
 $\delta=n''\ddot{n}_\psi=n\widetilde{n}_\phi\in \bbZ$. Then $n=\delta^2-n'\delta$ tells us
 $\delta$ divides $n$, but $\delta=n\widetilde{n}_\phi\ge n$, so $\delta=n$. Hence $n'=n-1$,
 so some $\ddot{n}_\psi$ must vanish, so $\ddot{n}_1=0$, so $1-nn''\delta^{-2}=0$
so $n''=n$, i.e. $H=G$. Since $|H'|=|H|$, we know $H'$ would also equal $G$. This means 
$\langle g,h\rangle$ is identically 1, so $gx=x$ for all $g,x$. We've just proved $\widetilde{n}_\phi
=1-\delta_{\phi,1}=\ddot{n}_\phi$ for all $\phi\in\widehat{G}$; in particular we can (and will)
identify $\cF$ with $\widehat{G}^*$ via $x\mapsto\widetilde{x}$, and then the assignment
$x\mapsto u_x$ corresponds to a permutation $\sigma$ of $\widehat{G}^*$. 

On the other hand, when $H=0$, $\langle g,h\rangle$ will be nondegenerate, and $H'$ also
equals 0. The element $g_x$ is then the unique one with $\mu^{g_x}=\overline{\widetilde{x}}$.
We know from the proof of Theorem 1 that the cardinalities $\widetilde{n}_\phi$ and $\widetilde{n}_{\mu^g\phi}$ 
must be equal for any $\phi\in\widehat{G}$ and $g\in G$, and so in this case they all
equal $\widetilde{n}_1=: k\in\bbZ_{\ge 0}$. Thus $n'=\sum_{\phi\in\widehat{G}}\widetilde{n}_\phi=nk$.
Now, $u_{x,g}$ is uniquely defined by its values at $g\in H$, where it is a 
character, so in this case it is identically 1 for all $g\in G$.  

Return to the general case (i.e. arbitrary $n'$).
The equivariance $\alpha_g\rho=\rho$ yields the selection rule: $a_{x,y}\ne 0$
implies $\widetilde{x}\,\widetilde{y}=1$. When $n'=n-1$, this means $a_{x,y}=a'(x)\delta_{y,\overline{x}}$
for some $a'(x)\in \bbC$. But $\sum_z|a_{x,z}|^2=\delta^{-1}$ then implies $|a'(x)|^2=\delta^{-1}$,
so $a(x):=a'(x)\sqrt{\delta}\in\bbT$.\qquad \textit{QED to Theorem 2}

\medskip For the convenience of the next two subsections, let us run through the identities
which must be satisfied by the numerical invariants $s,a,b,b',b''$, in order that they define
a near-group C$^*$-category of type $n'$. Let $H=G$ respectively 0, and $K=0$ respectively $G$,
for $n'=n-1$ and $n|n'$ respectively. Then we know from Theorem 2 that $\langle *,*\rangle$ 
is symmetric and nondegenerate on $K$. Write $u_{x,g}=\ddot{x}(g)$, where $\ddot{x}\in
\widehat{G}$ equals 1 on $K$. Assume that $a_{x,y}=\sqrt{\delta}^{-1}a_x\delta_{y,\overline{ x}}$ for 
some order-2 permutation $x\mapsto\overline{x}$ of $\cF$ (we already know this when $n'=n-1$,
and will prove it in subsection 2.3 when $n$ divides $n'$).

Define endomorphisms $\rho,\alpha_g,U_g$ on the Cuntz algebra $\cO_{n,n'}$, as in Theorem 1.
It is immediate that $\alpha_g$ defines a well-defined $G$-action on $\cO_{n,n'}$, and
$U_g$ a unitary representation of $G$. In order for $\rho$ to be a well-defined endomorphism
on $\cO_{n,n'}$, we need it to preserve the Cuntz relations
$S_g^*S_h=\delta_{g,h}$, $S_g^*T_z=0=T_z^*S_g$, $T^*_{z}T_{z'}=\delta_{z,z'}$
and $\sum_gS_gS_g^*+\sum_zT_zT^*_z=1$.  $(\rho S_g)^*(\rho S_h)=\delta_{g,h}$
reduces to \eqref{Cu:SS00}. $(\rho S_g)^*(\rho T_z)=0$ (or its adjoint) 
 is equivalent to 
\begin{equation}-sn\delta_{\widetilde{w},\overline{\mu^g}}\,
b_{z,w}=\sqrt{\delta}\sum_{x}\overline{\ddot{x}(g)}\,
\overline{a_{gx}}\,b''_{z;x,\overline{gx},w}\,,\label{CuST}\end{equation}
while the  relation $(\rho T_z)^*(\rho T_{z'})=\delta_{zz'}\,
(\sum S_hS_h^*+\sum T_xT_x^*)$ gives 
\begin{equation}\delta_{z,z'}\delta_{y,y'}=n\delta_{\widetilde{y},\widetilde{y'}}\,\overline{b_{z,y}}\,
b_{z',y'}+\sum_{w,x}\overline{b''_{z;w,x,y}}\,b''_{z';w,x,y'}\label{Cu:TTa}\end{equation}
and the unitarity of $b'$: $\sum_{x}\overline{b'_{z,x}}\,b'_{z',x}=
\delta_{z,z'}$. 
Finally, completeness $1=\sum_g\rho S_g\,\rho S_g^*+\sum_z\rho T_z\,\rho T_z^*$ is equivalent
to \eqref{compleeta} (with $H'=H$), unitarity of $b'$, and
\begin{align}
&&-{\delta}^{-3/2}a_z\sum_g\ddot{x}(g)\,\delta_{gx,\overline{z}}\,\overline{\mu_h(g)}=\sum_{w,y}
\widetilde{y}(h)\,\overline{b_{w,y}}\,{b''_{w;x,z,y}}\,,\label{completb}\\
&&\delta_{x,x'}\delta_{z,z'}=\delta^{-1}a_z\overline{a_{z'}}\sum_g\ddot{x}(g)\,\delta_{gx,\overline{z}}
\,\overline{\ddot{x'}(g)}\,\delta_{gx',\overline{z'}}+\sum_{w,y}b''_{w;x,z,y}\,\overline{b''_{w;x',z',y}}\,.\label{completdg}
\end{align}

To establish $\alpha_g\rho(x)=\rho(x)$ and $\rho(\alpha_g(x))=\mathrm{Ad}U_g(\rho(x))$ 
for all $x$, it suffices to prove both for $x=S_h$ and 
$x=T_z$. The first follows from the bilinearity of $\langle g,h\rangle$, that $\widetilde{gx}
=\mu^g\widetilde{x}$, and that $\widetilde{\overline{ x}}=\overline{\widetilde{x}}$. The second 
follows from the factorisations of $b_z(g,x)$ and $b'_z(x,g)$ given in Theorem 1, and the
selection rule $b''_{z;w,x,y}\ne 0\Rightarrow \widetilde{y}=\widetilde{w}\widetilde{x}$.
The identity $\rho(\alpha_g S_h)=U_g\rho(S_h)U_g^*$ is built into \eqref{rhoS}, while  
$\rho(\alpha_g T_z)=U_g\rho(T_z)U_g^*$ is implied by the covariances
\begin{align}
&&b_{z,gx}={\widetilde{z}(g)}\,{\ddot{x}(g)}\,b_{z,x}\,,
\label{ssrUaaa}\\
&&b'_{z,gx}=\widetilde{z}(g)\,\overline{\ddot{x}(g)}\,b'_{z,x}\,,\label{rhoalpTbb}\\ 
&&b''_{z;gw,x,gy}={\widetilde{z}(g)}\,\overline{\ddot{w}(g)}\,{\ddot{y}(g)}\,b''_{z;w,x,y}\,.
\label{rhoalpTcc}\end{align}

All that remains is to consider are the fusion rules. 
We require  $S_g$ to be in the intertwiner space Hom$(\alpha_g,\rho^2)$. We will
follow as much as we can the proof of Lemma 5.1(a) in \cite{iz3}. Thanks to Lemma 2.2 in \cite{iz3}, $S_g\in\mathrm{Hom}(\alpha_g,\rho^2)$ iff
$S_g^*\rho^2(x)S_g=\alpha_gx$ for all generators $x$. Hitting with $\alpha$
we see that this is true iff $S_0^*\rho^2(x)S_0=x$ for all generators $x$.
 Because $S_0^*\rho(U_g)=S_g^*$, the calculation in the
middle of p.625 of \cite{iz3} still goes through and $S_0^*\rho^2(S_h)S_0=S_h$ will follow
once we know it for $h=0$. So we have learned that 
$S_g\in\mathrm{Hom}(\alpha_g,\rho^2)$ iff both  $S_0^*\rho^2(S_0)S_0=S_0$ and $S_0^*\rho^2(T_z)S_0=T_z$ for all $z$.
Those two identities  are equivalent to
\begin{align}
&&{n'}\sqrt{\delta}^{-1}=\sum_{w,x}a_{x}\,b_{x,w}\,b'_{\overline{ x},w}\,,\label{Sintertwa}\\
&&\delta_{z,w}={sn}{\delta}^{-1}\sum_{x}b_{z,x}\,\overline{b_{x,w}}+
{n}{\delta}^{-3/2}\sum_{x}b'_{z,x}\,b_{x,\overline{ w}}\,a_{w}\nonumber\\&&
+\sum_{x,y,x',y',z'}b''_{z;x',x,y}\,b_{x',y'}\,b''_{x;y',w,z'}\,\overline{b_{y,z'}}\,,\label{Sintertwb}\end{align}
respectively, where we use \eqref{ssrUaaa} to simplify \eqref{Sintertwb}.

Clearly, if  $T_z$ is in the intertwiner space Hom$(\rho,\rho^2)$, then
 $\rho(y^*\rho(x))T_z=\rho(y^*)T_z\rho(x)$ for all generators $x,y$.
Conversely, if $\rho(y^*\rho(x))T_z=\rho(y^*)T_z\rho(x)$ for all generators $x,y$,
then the calculation
\begin{align}
\rho^2(x)T_z&=&\sum_g\rho(S_gS_g^*)\rho^2(x)T_z+\sum_w\rho(T_wT_w^*)\rho^2(x)T_z
\nonumber\\&=&\sum_g\rho(S_g)\rho(S_g^*)T_z\rho(x)+\sum_w\rho(T_w)\rho(T_w^*)T_z\rho(x)=T_z\rho(x)
\nonumber\end{align}
shows that $T_z^*\rho^2(x)T_z=\rho(x)$ for all generators $x$, and Lemma 2.2 of
\cite{iz3} then would imply $T_z\in\mathrm{Hom}(\rho,\rho^2)$. 
We compute
\begin{equation}
T_w^*\rho(U_g)=\sum_{x,y,z}V_{wxyz}(g)T_xT_y^*T_z^*
+\sum_{h,y}W_{why}(g)S_hS_h^*T_y^*+\sum_{h,z}X_{wzh}(g)T_zS_h^*\,,\nonumber\end{equation}
where 
\begin{align}
&&V_{wxyz}(g)=\delta^{-1}a_x\,\overline{a_y}\,\sum_k\langle k,g\rangle\,\ddot{w}(k)\,\overline{\ddot{z}(k)}\,\delta_{kw,\overline{x}}\,\delta_{kz,\overline{y}}
+\sum_{z',y'}\ddot{z'}(g)\,b''_{z';w,x,y'}\,\overline{b''_{gz';z,y,y'}}\,,\nonumber\\
&&W_{why}(g)=\overline{\widetilde{y}(h)}\,\widetilde{w}(h)\,\sum_z\ddot{z}(g)\,b'_{z,w}\,\overline{b'_{gz,y}}\,,\nonumber\\
&&X_{wzh}(g)={\delta}^{-3/2}a_z\sum_k\langle k,g-h\rangle\,{\delta_{kw,\overline{z}}\,\ddot{w}(k)}+\sum_{x,y}
\ddot{x}(g)\,b''_{x;w,z,y}\,\widetilde{y}(h)\,\overline{b_{gx,y}}\,.\nonumber\end{align}

The identity
$\rho(S_h^*\rho(S_g))T_w=\rho(S_h^*)T_w\rho(S_g)$ gives
\begin{align}
&& \langle g,h\rangle\,\overline{W_{wky}(g)}= s\,\delta\,
\overline{\langle g,k\rangle}\sum_z\ddot{z}(g)\,\overline{a_{hw,z}}\,a_{gz,hy}\,,\label{rSrSTb}\\
\label{rSrSTa}
&&s\langle g,h\rangle\overline{V_{wxyz'}(g)}=\delta_{w,gz'}\,\delta_{y,gx}\,
\ddot{z}(h)\,\ddot{\overline{hw}}(g)\,\overline{\ddot{x}(g)}\,\overline{\ddot{w}(h)}\,
{a_{\overline{(h-g)w}}}\,\overline{a_{hw}}\,,\\
&&X_{wzh}(g)=0\,.\label{X0}\end{align}
Using 
\eqref{X0}, the identity $\rho(T_x^*\rho(S_g))T_w=\rho(T_x^*)T_w\rho(S_g)$ gives
\begin{align}
&&\overline{\widetilde{w}(k)}\,\langle g,k-h\rangle\,\overline{b'_{x,w}}=s\,\sqrt{\delta}\,a_{gx}\,
 \ddot{x}(g)\,\sum_{y}
b_{\overline{g x},y}\,\overline{\widetilde{y}(k)}\,\overline{W_{why}(g)}\,,\label{rTrSTa}\\
&&\overline{b'_{x,w}}\,\overline{\widetilde{w}(h)}\,\langle h,g\rangle\,\overline{\ddot{z}(g)}\,\delta_{y,gz}=s\sqrt{\delta}\, \ddot{x}(g)\,
\sum_{x'}a_{gx}\,b_{\overline{g x},x'}\,\overline{\widetilde{x'}(h)}\,\overline{V_{wzyx'}(g)}\,,
\label{rTrSTb}\\
&&s\,\overline{\langle g,h\rangle}\,\overline{b''_{x;w,\overline{gx'},w'}}\,\ddot{\overline{gx'}}(g)\,a_{x'}=\ddot{x}(g)\,
a_{gx}\,\sum_yb''_{\overline{gx};w',x',y}\,\overline{W_{why}(g)}\,,\label{rTrSTc}\\
&&s\,\delta_{y,gy'}\,\overline{\ddot{y'}(g)}\,\overline{b''_{x;w,\overline{g x'},w'}}\,\ddot{\overline{ x'}}(g)\,a_{x'}=\ddot{x}(g)\,a_{gx}\sum_{z'}
b''_{\overline{g x};w',x',z'}\,\overline{V_{wy'yz'}(g)}\,.\label{rTrSTd}
\end{align}
$\rho(S_g^*\rho(T_z))T_x=\rho(S_g^*)T_x\rho(T_z)$ and $\rho(T_x^*\rho(T_z))T_w=\rho(T_x^*)T_w\rho(T_z)$  now simplify to
\begin{align}
&&\sqrt{\delta}\,\sum_{w}\overline{\widetilde{w}(g)}\,b_{z,w}\,\overline{b'_{w,x}}
=\overline{\ddot{x}(g)}\,\overline{a_{gx}}\,b'_{z,-g\overline{ x}}\,,\label{rSrTTa}\\
&&\sqrt{\delta}\,\sum_{x'}\overline{\widetilde{x'}(g)}\,b_{z,x'}\,\overline{b''_{x';x,w,y}}=
\ddot{y}(g)\,\overline{\ddot{x}(g)}\,\overline{a_{gx}}\,b''_{z;\overline{gx},gy,w}\,,\label{rSrTTb}\\
&&\overline{b'_{x,w}}
\,b_{z,y}-{\delta}^{-3/2}\,\overline{a_y}\,b'_{z,x}\sum_k\widetilde{x}(k)\,\overline{\ddot{w}(k)}
\,\delta_{ky,\overline{ w}}
= \sum_{w',z',x'}b''_{z;x,w',z'}\,b_{w',x'}\,\overline{b''_{z';w,y,x'}}\,,\label{rTrTTa}\\
&&\sum_{y',w',x'}\overline{b'_{y',z}}\,b''_{y';x,w',x'}\,b'_{w',y}\,\overline{b'_{x',w}}
=\overline{b''_{x;w,z,y}}\,,\label{rTrTTb}\\
&&\sum_{y'}\overline{b''_{x;w,y',x'}}\,b''_{z;y',y,z'}-\delta^{-1}\,a_y\,\overline{a_{z'}}\,b'_{z,x}\sum_h\widetilde{x}(h)\,\ddot{x'}(h)\,\overline{\ddot{w}(h)}\,\delta_{hy,\overline{ x'}}\,\delta_{hz',\overline{ w}}\qquad
\nonumber\\&&=\sum_{w',y',w''}b''_{z;x,w',w''}\,b''_{w';x',y,y'}\,\overline{b''_{w'';w,z',y'}}\,,\label{rTrTTc}
\end{align}
where \eqref{rTrTTa} was simplified using the selection rules $a_{kw,y}\ne 0\Rightarrow \mu_k
\widetilde{w}\widetilde{y}=1$ and $b''_{z';w,y,x'}\ne 0\Rightarrow \widetilde{x'}=\widetilde{w}\widetilde{y}$, both obtained earlier. \eqref{rTrTTb} was simplified using unitarity of $b'$ and
the selection rule for $b''$. 

 Let $M$ be
the weak closure of the Cuntz algebra $\mathcal{O}_{n,n'}$ in the GNS representation of a KMS state (as in \cite{iz0}, Remark 4.8). Then the endomorphisms $\alpha_g$ and $\rho$ extend
to $M$ and obey the same fusions as sectors.
Note that the $\alpha_g$ are outer because if they were implemented by a unitary, it would
have to commute with $\rho$ since $\alpha_g\rho=\rho$, but as $\rho$ is irreducible 
only the scalars can commute with it.

We have proved:\medskip

\noindent\textbf{Proposition 1.} \textit{Fix an abelian group $G=H\times K$, where $H=G$ or $H=0$,
and a symmetric pairing $\langle g,h\rangle$ nondegenerate on $K$. Let $a_{x,y}=\sqrt{\delta}^{-1}a_x\delta_{y,\overline{ x}}$ for some permutation $x\mapsto\overline{x}$ of the index set $\cF$ with $\overline{\overline{ x}}=x$ and 
$\widetilde{\overline{ x}}=\overline{\widetilde{x}}$, and suppose $b''$ obeys the selection rule
$b''_{z;w,x,y}\ne 0\Rightarrow \widetilde{y}=\widetilde{w}\widetilde{x}$.
Suppose the quantities $s,a,b,b',b''$ satisfy the equations \eqref{Cu:SS00}, unitarity of $b'$, \eqref{compleeta}, \eqref{CuST}-\eqref{rTrTTc}.
Then $\alpha_g,\rho$ defined as in Theorem 1 yield a near-group $C^*$-category of
type $G+n'$ for $n'=\|\cF\|$.}\medskip

We can now identify the principal graph of the subfactor $\rho(M)\sdprod G\subset M$ of index
$d_\rho^2/n=1+n'n^{-1}d_\rho$, introduced at the end of Section 2.1. Write the inclusions as $\iota:\rho(M)\subset \rho(M)\sdprod G$
and $j:\rho(M)\sdprod G\rightarrow M$. Then as $M$-$M$ sectors the canonical endomorphism 
$j\overline{j}$ is a subsector of the canonical endomorphism
$j\iota\overline{\iota}\overline{j}=\rho^2$, i.e. of $\sum_g[\alpha_g]+n'[\rho]$, which contains
$[\alpha_0]$.

 Consider first $n'=n-1$; then $d_\rho=n$ is the desired index, so the only possibility for
 $[j\overline{j}]$ is $\sum_g[\alpha_g]$. We see that the principal graph matches that
 of the orbifold $M^G\subset M$. This means this subfactor is isomorphic to $M^G\subset M$,
 up to a 3-cocycle. Could this 3-cocycle be related to the 3-cocycle appearing at the end of
 Subsection 4.2?
 
 Now consider $n'$ a multiple of $n$. When $n'=0$ there is nothing to say: the subfactor
 has index 1 so is trivial. When $n'>0$, there is only one possibility for $[j\overline{j}]$,
 namely $[\alpha+0]+n'n^{-1}[\rho]$. We recover the graph as in Figure 3, i.e. the $2^n1$
 graph but with the  $n$ valence-2 vertices attached to the central vertex with $n'n^{-1}$
 edges.
This generalises the paragraph concluding Section 5 of  \cite{iz3}, which in the
$G+n$ case considered there associates 
a subfactor $M\supset\rho(M)\sdprod G$ 
with canonical endomorphism  $[\alpha_0]\oplus[\rho]$, 
index $\delta+1$ and principal graph $2^n1$.

\medskip\epsfysize=1.3in\centerline{ \epsffile{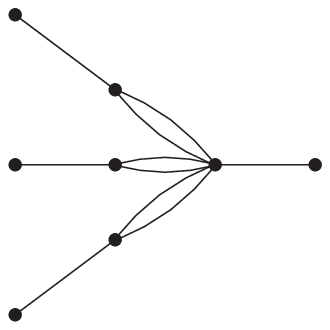}}\medskip

\centerline{Figure 3. Principal graph for the intermediate subfactor for type $\bbZ_3+6$}
\label{fig3}\medskip

\subsection{First class: Near-group categories with $n'=n-1$}

Theorem 2 says that there are two classes of near-group C$^*$-categories: $n'=n-1$ and $n'\in n\bbZ$. In this subsection
we focus on the former, and identify a complete set of 
 relations satisfied by the numerical invariants $s,...,b''$ of Corollary 1. 
We know from Corollary 2 that these systems are always realised by the even part of
a subfactor.

\medskip\noindent\textbf{Theorem 3(a)} \textit{Let $G$ be an abelian group of order $n$. Put
$\delta=n$ and $\cF=\widehat{G}^*:=\widehat{G}\setminus\{1\}$. Let $\sigma$ be
 a permutation of $\widehat{G}^*$ satisfying 
\begin{equation}
\sigma(\overline{a})=\overline{\sigma^{-1}a}\,,\ \sigma^3=id\,, \mathrm{and}\
\sigma(\sigma a\,\overline{\sigma b})=\sigma^2a\,\sigma(b\overline{a})\,,\label{sigma}\end{equation}
for all $a\ne b\in G$.  Put $\widetilde{x}=x$, $u_{x,g}=(\sigma\,x)(g)$, 
$a_{x,y}=\sqrt{n}^{-1}\delta_{y,\overline{x}}\,a(x)$, $b_{x,y}=\sqrt{n}^{-1}\delta_{y,\overline{\sigma\,x}}\,
b(x)$, $b'_{x,y}=s\,\delta_{y,\sigma^2\,x}\,\overline{b(\overline{x})}\,a(x)$, and $b''_{z;w,x,y}=\delta_{z,\sigma w\,\overline{\sigma y}}
\,\delta_{y,wx}\,b''(w,x)$ for quantities $a(x)\in\{1,s\}$, $b(x),b''(x,y)\in\bbT$ 
(provided $xy\ne 1$). Suppose these parameters satisfy $a(x)=
a(\sigma x)=sa(\overline{x})$, $b(\overline{\sigma x})=sb(x)$, $b(x)\,b(\sigma x)\,b(\sigma^2 x)=s\,a(x)$ and
\begin{align}
&&b''(x,{y})=s\,a(y)\,a(\sigma(x{y})\,\overline{\sigma x})\,\overline{b''(x{y},\overline{y})}
\qquad\forall xy\ne 1\,,\label{bpplcbpp}\\
&&b''({x},y)=s\,a(x)\,b(\sigma{x}\,\overline{\sigma({x}y)})\,\overline{b''(\overline{x},xy)}\qquad\forall xy\ne 1\,,\label{bpprcbpp}\\
&&b''({x},y)=s\,b(\overline{x})\,\overline{b(y)}\,b(x{y})\,b(\overline{\sigma^2 x}\,\overline{\sigma y})\,\overline{b''(\overline{\sigma^2 x},\overline{\sigma y})}\qquad\forall xy\ne 1\,,
\label{bpppibpp}\\
&&b''(\sigma w\,\overline{\sigma x},\sigma x\,\overline{\sigma(xy)})\,b''(x,y)\,b''(w,x\overline{w})=b''(w,xy\overline{w})\,b''(x\overline{w},y)\,,\label{bppbppbpp}\end{align}
where the last equation requires $w\ne xy$, $xy\ne 1$, and $w\ne x$. Then $\alpha_g$,
$U_g$, and $\rho$ defined as in Theorem 1 constitute a near-group $C^*$-category
of type $G+(n-1)$.}

\smallskip\noindent\textbf{3(b)} \textit{Conversely, let $\cC$ be a near-group $C^*$-category of type
$G+(n-1)$, and assume $H^2(G;\bbT)=1$. Then $\cC$ is $C^*$-tensor equivalent to one in
part 3(a).}
\medskip

\noindent\textit{Proof.} We'll prove part (b) first. Theorem 2(b) tells us $\delta=n$, we can 
identify $\cF$ with $\widehat{G}^*$ through $x\mapsto\widetilde{x}$, $a_{x,y}=\sqrt{n}^{-1}
\,a(x)\,\delta_{y,\overline{x}}$, and there is a permutation $\sigma$ of $\widehat{G}^*$
obtaining $u_{x,g}$ as $(\sigma\,x)(g)$.

Recall the Cuntz algebra  $\cO_{n,n'}$ generated by the $S_g$ and $T_z$. 
Select representatives $z\in\cR$ of each $\bbZ_2$-orbit $\{z,\overline{z}\}$ in $\cF$; then
by rescaling the $T_z$ appropriately we can fix the values of $a_{z,\overline{z}}$ to be
1 for $z\in\cR$. Now, if  $T_w$ is in the intertwiner space Hom$(\rho,\rho^2)$, then
 $\rho(y^*\rho(x))T_w=\rho(y^*)T_w\rho(x)$ for all generators $x,y$. In particular,
the $T_yS_kS_k^*$ coefficient of $\rho(S_h^*\rho(S_g))T_w=\rho(S_h^*)T_w\rho(S_g)$ reads
\begin{equation}
s\delta^{-1}\sum_z\overline{\ddot{z}(g)}\,\overline{b'_{z,w}}\,b'_{gz,y}=\overline{\ddot{w}(h)}\ddot{y}(h)\sum_z\ddot{z}(g)\,\overline{a_{hw,z}}\,a_{gz,hy}\,.\label{rSrSTb}
\end{equation}
Putting $g=h=0$ in \eqref{rSrSTb} gives $\overline{a_{x}}\,a_{\overline{x}}=s$, and hence
$a_{\overline{z}}=s$ for all $z\in\cR$. 

We get
selection rules for $b,b',b''$ through the equivariance $\alpha_g\rho=\rho$ and
the identity $\rho(\alpha_g T_z)=U_g\rho(T_z)U_g^*$, namely:
$b''_{z;w,x,y}\ne 0$ implies $y={w}\,{x}$ and $z=\sigma{w}\,\overline{\sigma{y}}$;
$b_{z,x}\ne 0$ implies ${z}\,{\sigma{x}}=1$; and $b'_{z,x}\ne 0$ implies
 ${z}=\sigma(x)$. Therefore we can write 
$b_{x,y}=\sqrt{n}^{-1}\delta_{y,\overline{\sigma\,x}}\,
b(x)$, $b'_{x,y}=\delta_{y,\sigma\,x}\,b'(x)$, and $b''_{z;w,x,y}=\delta_{z,\sigma w\,\overline{\sigma y}}
\,\delta_{y,wx}\,b''(w,x)$ for quantities $b(x),b'(x),b''(w,x)\in\bbC$. Note that
$b''(w,x)=0$ when $wx=1$ because 1 is a forbidden value for $y=wx$. \eqref{compleeta} forces
$b_x\in\bbT$. $g=h=0$ in \eqref{rSrSTb} forces unitarity of the matrix $b'$. The $S_0S_0S_0^*$
coefficient
of the identity $\rho(T_x^*\rho(S_0))T_w=\rho(T_x^*)T_w\rho(S_0)$ (again coming from
$T_w\in\mathrm{Hom}(\rho,\rho^2)$) gives $b'(x)=sa(x)\overline{b(\overline{x})}$.
The $T_xT_yT_{y'}^*T_{x'}^*$ coefficient of 
$1=\sum_g\rho S_g\,\rho S_g^*+\sum_w\rho T_w\,\rho T_w^*$ collapses now to
\begin{equation}
\delta_{x,x'}\delta_{y,y'}=\delta_{y,\overline{x}}\,\delta_{y',\overline{x'}}\,\delta_{x,x'}
+b''(x,y)\,\overline{b''(x',y')}\,\sum_{w,z}\delta_{w,\sigma x\overline{\sigma y}}\,\delta_{z,xy}\,\delta_{w,\sigma x'\overline{\sigma y'}}\,\delta_{z,x'y'}\,.\label{completdgg}
\end{equation}
Choose any $x,y\in\widehat{G}^*$ with $xy\ne 1$; we claim that the only solution $x',y'\in
\widehat{G}^*$ with $xy=x'y'$ and $\sigma x\,\overline{\sigma y}=\sigma x'\,\overline{\sigma y'}$
is $x=x'$ and $y=y'$: otherwise \eqref{completdgg} would force $b''(x,y)\,\overline{b''(x',y')}=0$,
which contradicts \eqref{completdgg} with $x=x'$ and $y=y'$. Thus each $b''(x,y)\in \bbT$
(provided $xy\ne 1$).



 The $S_0S^*_0$ respectively $T_yT^*_w$ coefficients of 
$\rho(S_g^*\rho(T_z))T_x=\rho(S_g^*)T_x\rho(T_z)$, together with $b'=s\delta\overline{ab}$, 
 gives 
$\overline{\sigma^2z}=\sigma\overline{z}$ and $b(z){b(\overline{z})}b(\sigma z)
 =a(z)a(\sigma z)a(\sigma^2 z)$, respectively $\sigma(\sigma a\overline{\sigma b})=\sigma^2a\,
 \sigma(b\overline{a})$ and \eqref{bpprcbpp}. Taking $\sigma$ of the complex conjugate of
 $\sigma(\sigma a\overline{\sigma b})=\sigma^2a\,
 \sigma(b\overline{a})$, we obtain $\sigma^3=1$; iterating \eqref{bpprcbpp} twice
 gives $b(y)=s\,b(\overline{\sigma y})$. The $T_yT_{z'}S_0S_0^*$ coefficient of
  $\rho(T_x^*\rho(S_0))T_w=\rho(T_x^*)T_w\rho(S_0)$ recovers \eqref{bpplcbpp}.
The $T_yS_0S_0^*$ and $T_yT_{z'}S_0S_0^*$ coefficients of 
$\rho(T_x^*\rho(T_z))T_w=\rho(T_x^*)T_w\rho(T_z)$ give \eqref{bpppibpp} and
\eqref{bppbppbpp}.

Now that we know $\sigma$ has order dividing 3, we know we can choose the $\bbZ_2$-orbit
representatives $\cR$ so that $a_{z}$ is constant on $\sigma$-orbits. Indeed, if $\sigma^ix=\overline{x}$  
 for some $i$ and $x$, then $\overline{\sigma^{-i}x}=\sigma^{-i}x$ and $s$ must equal 1,
 so there is nothing to do; when there is no such $i,x$, there is no obstruction to putting all
 of $x,\sigma x,\sigma^2 x$ in $\cR$.

Conversely, suppose $s,a(x),b(x),b''(x,y)$, and $\sigma$ are as in Theorem 3(a). 
We need to verify the conditions of Proposition 1 are satisfied.
For this purpose note that: (i) $\sum_xx(g)=n\delta_{g,0}-1$ since $x$ runs over $\widehat{G}^*$;
(ii) that given a pair $w,x\in\widehat{G}^*$, there will be $y,z\in\widehat{G}^*$ with
$b''_{z;w,x,y}\ne 0$ (namely $y=wx$ and $z=\sigma w\,\overline{\sigma y}$) iff $x\ne\overline{w}$; 
and (iii) that given a pair $y,z\in\widehat{G}^*$, there will be $w,x\in\widehat{G}^*$ with
$b''_{z;w,x,y}\ne 0$ (namely $w=\sigma^2((\sigma y)z)$ and $x=y\overline{w}$) iff $\overline{z}\ne
\sigma y$. 

The permutation in Proposition 1 is the usual complex-conjugation $\overline{x}$ of characters. One easily computes $W_{why}(g)
=\delta_{w,y}\,(\sigma^2w)(g)$ and $V_{wxyz}(g)=\delta_{w,z}\delta_{x,y}(\sigma^2{w})(g)\,
{(\sigma x)(g)}$. In the last term of \eqref{Sintertwb},
$x,y,y',z'$ are determined from $x',z,w$, and is nonzero precisely when $z=w\ne x'$.
That the second term in $X_{wzh}(g)$ vanishes, follows because $\sigma^2(w)\,\sigma(z)
=\overline{wz}$ implies $wz=1$; this can be seen directly from \eqref{sigma} but is
trivial once we have Proposition 2 below. Both \eqref{rTrSTc},\eqref{rTrSTd} follow from 
\eqref{bpplcbpp}, \eqref{rSrTTb} comes from \eqref{bpprcbpp}, and both 
\eqref{rTrTTa},\eqref{rTrTTb} follow from \eqref{bpppibpp}. \eqref{rTrTTc} follows from
\eqref{bppbppbpp} and (when $y=\overline{x'}$) \eqref{bpplcbpp},\eqref{bpprcbpp}.
\textit{QED to Theorem 3}\medskip
%
%

When such  a permutation $\sigma$ exists, we get a solution by taking $s=1$, $a(x)=b(x)=b''(x,y)
=1$. This solution corresponds to Mod(Aff$_1(\bbF_q)$), as we explain at the end of
Subsection 4.3.
It is possible to classify all solutions $\sigma$ to \eqref{sigma} ---  they are essentially unique
when they exist. The key observation in the following proof  (the relation to finite fields)
is due to Siehler \cite{Sie}. (Incidentally, an implicit unwritten hypothesis throughout
\cite{Sie} is that $G$ is abelian.)

\medskip\noindent\textbf{Proposition 2.} \textit{Let $G$ be a finite abelian group which possesses a 
solution $\sigma$ to \eqref{sigma}. Then $G\cong\bbZ_{q-1}$ for some prime power $q=p^k$.
Moreover, if $\sigma'$ is any other solution to \eqref{sigma}, then $\sigma'=\alpha\sigma\alpha^{-1}$
for some group automorphism $\alpha\in\mathrm{Aut}(G)$. Conversely,  any $G=\bbZ_{q-1}$
for $q=p^k$ has exactly $|\mathrm{Aut}(G)|=\phi(q-1)$ solutions $\sigma$ to \eqref{sigma}.}

\medskip\noindent\textit{Proof.} Suppose $G$ has a solution $\sigma$. For convenience in the
following proof,  write $G$ multiplicatively. Then \cite{Sie} explains how to give $G\cup\{0\}$
the structure of a field $\bbF$: the multiplicative structure of $\bbF$ is the multiplication in
$G$, supplemented by $0x=x0=0$; let $-1$ be the unique element in $G$ of order
gcd$(2,1+|G|)$ and write $-x=-1x$; addition in $\bbF$ is defined by $x+y=(\sigma(-x^{-1}y))^{-1}x$
when $x,y\in G$ and $x\ne -y$, supplemented by $0+x=x+0=x$ and $x+(-x)=0$. This means
$G$ is the multiplicative group of the finite field $\bbF$, and thus is isomorphic to $\bbZ_{q-1}$
for some power $q$ of a  prime.

Call this field $\bbF_\sigma$. Suppose there is a second solution $\sigma'$. Let $\alpha$ be
the field isomorphism $\bbF_\sigma\rightarrow\bbF_{\sigma'}$. Then $\alpha$ restricts to a
group isomorphism from $\bbF^\times_\sigma=G$ to $\bbF^\times_{\sigma'}=G$, i.e.
$\alpha\in\mathrm{Aut}(G)$. Conversely, given any solution $\sigma$ to \eqref{sigma} and
$\alpha\in\mathrm{Aut}(G)$, we get a new additive structure on $\bbF_\sigma$ given by
$x+\!'\,y=((\alpha\sigma\alpha^{-1})(-x^{-1}y))^{-1}x$ etc, corresponding to a solution
$\alpha\sigma\alpha^{-1}$ to \eqref{sigma}. This is a bijectiion, since $\alpha$ can be recovered
from $\sigma'=\alpha\sigma\alpha^{-1}$.

Conversely,  any $G=\bbZ_{q-1}$ with $q$ a power of a prime, can be regarded as the
multiplicative group $\bbF^\times$ of a finite field with $q$ elements (so $0\in G$ corresponds
to $1\in\bbF^\times$). Then $\sigma(x)=(1-x)^{-1}$ works.\qquad\textit{QED to Proposition 2}\medskip 


\noindent\textbf{Corollary 3.} \textit{Consider any equivalence class of $C^*$-categories of type 
$G+n-1$.  Fix any finite field $\bbF_{n+1}$, identify $G=\bbF^\times_{n+1}$ and define 
$\sigma' x=1/(1-x)$. Fix any assignment of signs $a'(x)\in\{1,s\}$ such that $a'(x)=a'(\sigma x)
=sa'(\overline{x})$.  There is a set $\cL$ of functions $f:
\widehat{G}^*\times\widehat{G}^*\rightarrow\bbZ$, defined in the proof below, such that any $C^*$-category is equivalent to one with:}
\begin{itemize}

\item{} \textit{$a(x)=a'(x)$ and $\sigma x=\sigma'x$ for all $x$;}

\item{}  \textit{$b(x)=sa(x)$ for all $x\ne -1$; in addition, $b(-1)=sa(-1)$ unless 
$n+1$ is a power of 3 in which case $b(-1)$ must be a third root of unity $\omega$;}

\item{} \textit{$\prod_{(x,y)}b''(x,y)^{f(x,y)}=1$ for all $f\in\cL$.} 

\end{itemize}

\noindent\textit{Conversely, any two $C^*$-categories with numerical invariants satisfying these constraints,
and with identical $b(x)$ and $b''(x,y)$, will be equivalent.  Finally,
 $s=1$ unless $n+1=q$ is a power of 2.}\medskip

\noindent\textit{Proof.} Note that gauge equivalence by a diagonal matrix $P$ with entries in
$\bbT$,  permits us to change $a(x)^{new}=P_xP_{\overline{x}}\,a(x)^{old}$, $b^{new}(x)
=\overline{P_x}\overline{P_{\overline{\sigma x}}}\,b^{old}(x)$, and $b^{\prime\prime new}(w,x)
=\overline{P_{\sigma w\,\overline{\sigma(wx)}}}P_wP_x\overline{P_{wx}}b^{\prime\prime old}(w,x)$.
First note that, for any given $x$, we can change both signs $a(x)$ and $a(\overline{x})$ (and 
leave all other $a(y)$ unchanged) by taking
$P_x=P_{\overline{x}}=\mathrm{i}$ and all other $P_y=1$. Now
 choose any $x$ with $\sigma x\ne x$.
Without loss of generality assume both $x\ne -1$ and $\sigma^2(x)\ne -1$. Then
$\overline{x}\ne x$ and $\overline{\sigma^2 x}\ne \sigma^2x$, so take 
$P_x=sa(x)b^{old}(x)=\overline{P_{\overline{x}}}$, $P_{\sigma x}=1=P_{\overline{\sigma x}}$,
$P_{\sigma^2 x}=sa(x)\overline{b^{old}(\sigma x)}=\overline{P_{\overline{\sigma^2x}}}$
and all other $P_y=1$. This gives $b^{new}(x)=  sa(x)=b^{new}(\sigma x)$. Then 
$b(x)b(\sigma x)b(\sigma^2x)=sa(x)$ forces $b(\sigma^2x)=sa(x)$, and $b(\overline{\sigma y})
=sb(y)$ forces $b(\sigma^i\overline{x})=a(x)=a(\overline{x})$.

$-1\in\widehat{G}^*$ precisely when $G$ has even order, i.e. precisely when
$\bbF_{n+1}$ has odd characteristic. In this case, $a(-1)=sa(\overline{-1})=sa(-1)$ so $s=+1$.
This forces all $a(x)=1$, and by the previous paragraph we know $b(x)=1$ unless $\sigma x=x$.
When $\sigma x=x$, the relation $b(x)b(\sigma x)b(\sigma^2 x)=sa(x)$ says $b(x)$ will
be a third root of 1.

Thanks to Proposition 2, we can find a finite field $\bbF_q$ with $G=\bbF^\times_q\cong\bbZ_n$
for which $\sigma x=1/(1-x)$. Suppose $\sigma x=x$. Then $x^2-x+1=0$ in $\bbF_q$, so
$x^3=-1$. If 3 does not divide $n$, the only solution to $x^3=-1$ is $x=-1$, but $\sigma(-1)=-1$
iff $q$ is a power of 3, i.e. $n\equiv 2$ (mod 3). When 3 divides $n$, there are thus exactly
2 fixed-points of $\sigma$. Let $f$ be one of these. Then we calculate $b''(f,f)=
\overline{b''(\overline{f},\overline{f})}$ from \eqref{bpplcbpp}, and $b''(f,f)=\overline{b(\overline{f})}
\overline{b''(\overline{f},\overline{f})}$ from \eqref{bpprcbpp}, so $b(\overline{f})=b(f)=1$.

By Corollary 1,
two C$^*$-categories are equivalent,   if they have identical numerical invariants
 $s,a,b,b',b''$ modulo gauge equivalence and
automorphism of $G$. Fixing $\sigma$ fixes the automorphism of $G$. Suppose we also
fix $a$ and $b$. The remaining gauge freedom are the quantities $P_z=P_{\sigma z}=\overline{P_{\overline{z}}}\in \bbT$. Note that $\overline{z}=\sigma^i z$ iff $\sigma^{-i}(z)=-1$ (the only
order-2 element in $\bbF_q^\times$).  Together, $\sigma$ and complex-conjugation  form a group isomorphic to $S_3$; the even-length orbits in $\widehat{G}^*$ are
precisely those which don't contain $-1$.  Let $\ell$ be the number of  even-length orbits; the
point is that for each of these,  the gauge phase $P_z\in\bbT$ is arbitrary. Associate to each pair
$(x,y)\in\widehat{G}^*\times\widehat{G}^*$, $xy\ne 1$, a vector $\vec{v}(x,y)\in\bbZ^\ell$
such that the gauge action is $b''(x,y)\mapsto b''(x,y)\prod_z P_z^{\vec{v}(x,y)_z}$. The $\bbZ$-span
of these $\vec{v}(x,y)$ is a lattice (of dimension $\le \ell$), and thus will have a $\bbZ$-basis;
$\cL$ consists of those linear combinations $\sum_{(x,y)}f(x,y)\vec{v}(x,y)$, $f(x,y)\in\bbZ$,
 corresponding to such a basis.
 \qquad \textit{QED to Corollary 3}\medskip

We will show in Proposition 5 below that $s=a=b=b'=b''=1$ except for $n=1,2,3,7$. One can see 
from the data in
Subsection 3.4 that for all $n<32$, the basis defining $\cL$ consists of exactly $\ell$ of the
$\vec{v}(x,y)$.

 Let $\cC(b,b'')$ denote the equivalence class corresponding to the constraints of Corollary 3
 (we suppress the choice of field $\bbF$ and subset $\cL$, though these are implicit). 
Note that there is an obvious product structure on the collection of equivalence classes of
type $G+n-1$ C$^*$-categories: $\cC(b_1,b''_1)*\cC(b_2,b''_2)=\cC(b_1b_2,b''_1b''_2)$ where  
$(b_1b_2)(x)=b_1(x)b_2(x)$ and $(b''_1b''_2)(x,y)=b_1''(x,y)b_2''(x,y)$. 
(Of course $s=s_1s_2$, $a(x)=a_1(x)a_2(x)$ and  $b'(x)=b_1'(x)b_2'(x)$).
Then $s,a,b,b',b''$ will obviously also satisfy the conditions of Theorem 3(a) and Corollary 3,
and thus uniquely determine an equivalence class of  type $G+n-1$ C$^*$-categories.
Likewise, the identity is $s=a=b=b'=b''=1$, and the inverse is complex-conjugation.
We find this abelian group structure very useful in Section 4.2  (though we will find in Section 3.4 
that this  group is usually trivial!). 

In any case, this group of equivalence classes should be closely related to the
set $H^2((\bbZ_p^k,\bbZ_{p^k-1});\bbT)/\!\!\sim$ defined in Chapter III of  \cite{IK}, as
suggested by their Theorem IX.8 (see also \cite{HS}). Those equivalence classes parametrise 
deformations of Kac algebras (equivalently, depth-2 subfactors) possessing what
we would call near-group fusions of type $\bbZ_{p^k-1}+p^k-2$.  We return to this briefly in
Section 4.2.

Note that in all cases in Corollary 3, $b(\sigma x)=b(x)=sb(\overline{x})$. 
In the following we fix a finite field $\bbF_q$, identify the labels  $x\in\widehat{G}^*$ with 
the entries of $\bbF_q\setminus\{0,1\}$, and choose $\sigma x=(1-x)^{-1}$, so that
(when $\bbF_q$ has odd characteristic, i.e. $n$ is even) $\sigma(-1)=1/2$ and $\sigma(2)=-1$.
Note from the proof that $\sigma$ will have exactly 0,1,2 fixed-points respectively,
for $n\equiv 1,2,0$ (mod 3).

By Corollary 2,  the principal graph for $\rho(M)\subset M$ when $n=2$ is $D_5^{(1)}$,
the McKay graph for binary $S_3$. This suggests an alternate construction of the subfactor
$\rho(M)\subset M$, at least when $b=b''=1$. Construct a central extension BAff$_1(\bbF_q)$
of Aff$_1(\bbF_q)$ by $\bbZ_2$ ---  it will have
precisely 2 $n$-dimensional irreps (one of which, denoted $\rho'$, is faithful) and $2n$ 1-dimensional irreps, and its 
McKay graph (which consists of  a node for each irrep and $m$ edges connecting node $i$ and $j$
if $k$ is the multiplicity of irrep $j$ in the tensor product of irrep $i$ with $\rho'$)
is the desired principal graph.
The irreps of BAff$_1(\bbF_q)$ are separated into
even and odd ones, depending on whether or not the centre is in the kernel; the even vertices
are precisely the irreps of Aff$_1(\bbF_q)$. To construct the  subfactor, start
 with the index-$n^2$ subfactor
\begin{equation}
\bbC I\otimes M_{n\times n}\otimes M_{n\times n}\otimes\cdots\subset M_{n\times n}\otimes M_{n\times n}\otimes M_{n\times n}\otimes\cdots\ ,\end{equation}
identifying $\mathrm{BAff}_1(\bbF_q)$ with its image $K\subset M_{n\times n}$ using the
faithful irrep $\rho'$, and then take  fixed-points:
\begin{equation}
\bbC I\otimes M_{n\times n}^K\otimes M_{n\times n}^K\otimes\cdots\subset M_{n\times n}^K\otimes M_{n\times n}^K\otimes M_{n\times n}^K\otimes\cdots\,.\end{equation}
The principal graph is constructed as in \cite{GHJ}; for this purpose it is important that
$\rho'$ is self-conjugate, as explained in \cite{W1}.

\subsection{The remaining class: $n'$ a multiple of $n$}

In this section we identify a complete set of 
 relations satisfied by the numerical invariants $s,...,b''$ of Corollary 1, for the remaining class
 of near-group C$^*$-categories, namely where the fusion coefficient multiplicity $n'$ is a multiple
 of the order $n$ of the abelian group $G$.  Corollary 2 tells us that when $n'$ is a positive
 multiple of $n$, the system will be realised by the even parts of a subfactor.
There seems no reason to expect all solutions with $n'=0$ to arise from subfactors.

 The main result in Theorem 4 is part (b).
As usual, we identify the basis
$\{T_z\}=\cF$ of Hom$(\rho,\rho^2)$ found in Theorem 1, with the set of labels $\{z\}$.
For a given symmetric pairing $\langle,\rangle$ for $G$, fix a function $\epsilon_{\langle,\rangle}:
G\rightarrow\bbT$ satisfying $\epsilon_{\langle,\rangle}(-g)=\epsilon_{\langle,\rangle}(g)$
and $\epsilon_{\langle,\rangle}(g+h)=\overline{\langle g,h\rangle}\epsilon_{\langle,\rangle}(g)\epsilon_{\langle,\rangle}(h)$. An immediate consequence is that $\epsilon_{\langle,\rangle}
(g)^2=\overline{\langle g,g\rangle}$.
For example, if $n=|G|$ is odd, the unique such function is $\epsilon_{\langle,\rangle}
(g)={\langle g,g\rangle}^{(n-1)/2}$, while if $G=\bbZ_{2k}$ and $\langle g,h\rangle=\exp(m\pi\i gh/k)$ one of the two is $\epsilon_{\langle,\rangle}(g)=\exp(-m\pi\i\,{g^2/(2k))}$.


 \medskip\noindent\textbf{Theorem 4(a)} \textit{Let $G$ be a finite abelian group.
Let $\langle k,k'\rangle=\langle k',k\rangle\in\bbT$ be a nondegenerate 
symmetric pairing  on $G$. 
For each $\psi\in \widehat{G}$ let $\mathcal{F}_1^{\psi}$ be a (possibly empty)  parameter set
and define $\cF_\phi^{\phi\psi}=\cF_1^\psi$ for all $\phi,\psi\in\widehat{G}$; then $\cF$ is the set
of all triples $x=(\widetilde{x},\dot{x},\check{x})$ where $\widetilde{x},\dot{x}\in\widehat{G},\check{x}\in\cF^{\dot{x}}_{\widetilde{x}}$.  Let $n'=\|\cF\|=n\sum_\psi\|\cF_1^\psi\|$ and define $\delta$ by \eqref{delta}.
Let $S_g,T_z$ be standard Cuntz 
generators, for $g\in G$ and $z\in\cF$. Define $\alpha_g$ and $U_g$ as in 
Theorem 1, where $u_{x,g}=1$ and $gx=(\mu^g\widetilde{x},\mu^{-g}\dot{x},\check{x})$ for $\mu^g(k)=\langle g,k\rangle$.  
Define $\rho(S_g)$ and $\rho(T_z)$ by \eqref{rhoS},\eqref{rhoT} where
 $a_{x,y}=\sqrt{\delta}^{-1}a_x\delta_{y,\overline{x}}$ and $b'_{x,y}=s\sqrt{\delta} \overline{a_{x}}
\,\overline{b_{\overline{x},y}}$, for $a_x=s_x\,\overline{\dot{x}(g_x)}\,
\overline{\epsilon_{\langle,\rangle}(g_x)}$, where $g_x\in G$ is as in Theorem 2(c), 
 for some signs $s_x\in\{1,s\}$ and some permutation $x\mapsto\overline{ x}$ of $\cF$ as in
 Proposition 1.  
Then $\cC(G,\alpha,\rho)$ is a $C^*$-category of type $G+n'$, provided the following equations are satisfied:
$s_{gx}=s_x$, $s_{\overline{ x}}=ss_x$, 
 $\dot{\overline{x}}=\overline{\dot{x}}$, 
$b_{y,x}=sb_{x,y}$,} 
\begin{align}
&&b_{z,gx}={\widetilde{z}(g)}\,b_{z,x}\,;
\label{ssrUaaa}\\
&&\overline{\widetilde{z}(g)}\,b''_{z;gw,x,gy}=b''_{z;w,x,y}=
\overline{\dot{w}(g)}\,\overline{\epsilon_{\langle,\rangle}(g)}\,b''_{-gz;gw,-gx,y}\,;\label{rhoalpTcc}\\
&&\delta\sum_yb_{x,y}\,\overline{b_{z,y}}=\delta_{x,z}=s\sum_{y,w}b'_{x,y}\,b'_{y,w}\,
{b'_{w,z}}
\,;\label{aunitar}\\
&&\delta^{-1}\,b''_{z;w',x',y}=s\sqrt{\delta}^{-1}\,a_{w'}\sum_{w}b_{z,w}\,\overline{b''_{w;\overline{w'},y,x'}}=\delta^{-1}\,a_{x'}\,\overline{a_{z}}
\,\overline{b''_{\overline{z};y,\overline{ x'},w'}}\nonumber\\&&
=\sum_{w,x}b_{y,w}\,\overline{b_{w',x}}\,\overline{b''_{x';x,z,w}}
\,;\label{rSrTTb}\\
&&n^{-1}\sum_{y}b''_{z;y,x,y}=-\delta^{-1}\,\delta_{z,x}\delta_{\widetilde{x},1}=\delta_{\widetilde{x},1}
\delta_{\widetilde{z},\widetilde{x}}\sum_yb''_{y;x,y,z}\,;\label{ssrUb}\\
&&\sum_{w,y}{b''_{z;w,x,y}}\,\overline{b''_{z';w,x',y}}=\delta_{z,z'}\delta_{x,x'}-n\delta^{-1}\delta_{\widetilde{x},\widetilde{x'}}\,{\delta_{z,x}}\,\delta_{z',x'}\,;\label{Cu:TT}\\
&&\sum_{w,y}b''_{w;y,z,x}\,\overline{b''_{w;y,z',x'}}=\delta_{x,x'}\delta_{z,z'}-\delta^{-1}
\delta_{g_xx,g_zz}\,\delta_{g_xx',g_zz'}\,;\label{completddgg}\\
&&\delta_{\overline{\widetilde{y}},\widetilde{x'}}\sum_{y'}\overline{b''_{x;w,y',x'}}\,b''_{z;y',y,z'}-\delta^{-1}b'_{z,x}\,a_{x'}\,\delta_{y,\overline{ x'}}\,\overline{a_{w}}\,\delta_{z',\overline{w}}\qquad
\nonumber\\&&=\delta_{\overline{\widetilde{y}},\widetilde{x'}}\sum_{w',y',w''}b''_{z;x,w',w''}\,b''_{w';x',y,y'}\,\overline{b''_{w'';w,z',y'}}\,;\label{rTrTTcc}
\end{align}
\textit{as well as the selection rules $b_{x,y}\ne0\Rightarrow\dot{x}\dot{y}=\overline{\widetilde{x}
\widetilde{y}}$, and} 
$b''_{z;w,x,y}\ne 0\Rightarrow\widetilde{y}=\widetilde{w}\,\widetilde{x}$. 

\smallskip\noindent\textbf{(b)} \textit{Conversely, let $\cC$ be a $C^*$-category of type $G+n'$ for
$n'\in n\bbZ$
and suppose $H^2(G;\bbT)=1$. Then there exist quantities $\langle *,*\rangle,s,a,b,b',b''$
satisfying the above equations and relations, such that the corresponding 
$\cC(G,\alpha,\rho)$ is tensor-equivalent
 to $\cC$.}\medskip

\noindent\textit{Proof.} We will prove part (b) first.
Recall the Cuntz algebra  $\cO_{n,n'}$ generated by the isometries $S_g$ and $T_z$. The desired
selection rule for $b''$ follows from the equivariance $\alpha_g\rho=\rho$.
The identity $\rho(\alpha_g T_z)=U_g\rho(T_z)U_g^*$ implies the
 covariances for $b,b',b''$, namely the left-sides of \eqref{ssrUaaa} and \eqref{rhoalpTcc},
 along with $b'_{z,gx}=\widetilde{z}(g)\,b'_{z,x}$.The $S_0S_0^*$-coefficient of the Cuntz
 relation  $\delta_{z,z'} =\rho T^*_z\,\rho T_{z'}$ forces the unitarity of $b'$.

The $T_yS_kS_k^*$ coefficient of the identity $\rho(S_h^*\rho(S_g))T_w=\rho(S_h^*)T_w\rho(S_g)$ 
(which holds because $T_w$ lies in the intertwiner space $\mathrm{Hom}(\rho,\rho^2)$) reads
\begin{equation}
\langle g,h\rangle s\delta^{-1}\sum_z\overline{b'_{z,w}}\,b'_{gz,y}=\sum_z\overline{a_{hw,z}}\,a_{gz,hy}\,.\label{rSrSTb}
\end{equation}
Putting $g=h=0$ in \eqref{rSrSTb} gives $\sum_y\overline{a_{x,y}}\,a_{y,z}=s\delta^{-1}\delta_{x,z}$.
Applying the triangle inequality to this, and comparing with $\sum_y|a_{x,y}|^2=\delta^{-1}$, yields
${a_{y,x}}=sa_{x,y}$. 

For each $\phi\in\widehat{G}$, recall the subspace $\cF_{\phi}$ consisting of all $T\in\mathrm{Hom}
(\rho,\rho^2)$ with $\alpha_g(T)=\phi(g)\,T$ $\forall g\in G$. 
Using the gauge freedom discussed after Corollary 1, we can now simplify the form of $a$.
We need this unitary change-of-basis $P$ to commute with each $U(g)$, i.e. to satisfy
$P(\cF_{\phi})= \cF_{\phi}$ and $P_{kx,ky}=P_{x,y}$ for all $k\in G$.
For any $g\in G$ define $A(g)_{x,y}=\delta\sum_za_{x,gz}\,\overline{a_{y,z}}$. Then 
\eqref{rSrSTb} tells us for any $g,k\in G$ the covariance 
$A(g)_{x,y}=\overline{\langle g,k\rangle}\,A(g)_{kx,ky}$ and
the selection rule $A(g)_{x,y}=0$ unless both $\widetilde{y}=\mu^g\widetilde{x}$.
 The latter tells us each $A(g)$ is a map from $\cF_{\phi}$ to $\cF_{\mu^g\phi}$. We verify that $A(0)=I$,
$A(g)\,A(k)=A(g+k)$ and $A(-g)=A(g)^*$ so $g\mapsto A(g)$ defines a unitary representation
of $G$. Moreover, covariance yields $A(g)U_k=\langle g,k\rangle U_kA(g)$. Hence
the set of unitary operators $A(g)U_{-g}$ for all $g\in G$ commute and so can be simultaneously
diagonalised; choose $P$ so that 
all $A(g)U_{-g}$ are diagonal on $ \cF_{1}$, and $P$ is defined on arbitrary
$ \cF_{\phi}$ by requiring $P_{kx,ky}=P_{x,y}$. 
This means $A(g)_{x,y}=\epsilon_x(g)\delta_{y,gx}$ for all $x,y\in\cF$, for some $\epsilon_x(g)\in\bbT$. Hence the properties of $A$ reduce to $\epsilon_{kx}(g)=\langle g,k
\rangle\epsilon_x(g)$ and $\epsilon_x(g)\epsilon_{gx}(h)=\epsilon_x(g+h)$. We see that
$\psi_x:=\epsilon_x\epsilon_{\langle,\rangle}$ lies in $\widehat{G}$. The triangle inequality
applied to $A(g)_{x,gx}$, together with $\delta\sum_z|a_{xz}|^2=1$, gives us the covariance
$a_{gx,-gy}=\overline{\epsilon_x(g)}\,a_{x,y}$. We need to refine this $P$ further.

Define an equivalence relation on $\cF$ by $x\sim x'$ iff there exists a sequence $x=x_0,
x_1,
\ldots,x_m=x'$ and $y_1,\ldots,y_m$ in $\cF$ such that the entries $a_{x_{i-1},y_i}$ and $a_{x_i,y_i}$ are
nonzero for all $1\le i\le m$. Let $\cX_x$ denote the equivalence class containing $x$
and write $\cT_x$ for span$_{w\in\cX_x}\{T_w\}$.
Then whenever $a_{x,y}\ne 0$, $a$ restricts to the indecomposable blocks $\cT_x\rightarrow
\cT_y$, where it is unitary. An induction argument (the base step of which was done in the previous paragraph)
verifies that any $w\in\cX_x$ has $\widetilde{w}=\widetilde{x}$.
Moreover, the invertibility of $a$ says $a_{x,y}\ne 0$ implies the cardinalities $\|\cX_x\|=\|\cX_y\|$
are equal.

Choose any $x\in \cF_1$, and suppose $a_{x,y}\ne 0$. Consider first the case where $\cX_x$
and  $\cX_y$ are disjoint, and fix some bijection $\pi:\cX_y\rightarrow\cX_x$. 
Define a unitary $u$ on $\cT_x+\cT_y$ to be the identity on $\cT_x$ and to be $a|_{\cT_x}\circ
\pi$ on $\cT_y$, and replace $a$ on $\cT_x+\cT_y$ with $u^T au$. Otherwise
we have the case $\cX_x=\cX_y$, so we can make use of some facts from linear algebra (see section 4.4 of
\cite{HoJo}) which say that: (i) when a complex matrix $B$ is both symmetric and normal, then
there exists a real orthogonal matrix $Q$ and a diagonal matrix $D$ such that $B=QDQ^T$;
(ii) when a complex matrix $B$ is both skew-symmetric and normal, then there exists
a real orthogonal matrix $Q$ such that $Q^TBQ=0\oplus\cdots\oplus 0\oplus\bigoplus_j\left(\begin{matrix}0&z_j\\ -z_j&0\end{matrix}\right)$ for $z_j\in\bbC^\times$. Our matrix $B$ here
(namely $a|_{\cT_x}$) is in fact unitary, so both the $z_j$ and the
diagonal entries of $D$ lie in $\bbT$ and we can adjust $Q$ by the square-roots of those
numbers and maintain unitarity. The result is a matrix $a$ in the form described in the statement
of Theorem 4,
where we  write $\pi x=\overline{x}$ and  $\dot{x}(g)=\overline{\psi_x(g)}$, and decompose 
each $\cF_\phi$ into $\oplus_\psi\cF_\phi^\psi$ where $\dot{x}=\psi$ for $x\in\cF_\phi^{\psi_x}$.
This means we write $x\in\cF$ as
a triple $(\widetilde{x},\dot{x},\check{x})$ where $\check{x}\in\cF_{\widetilde{x}}^{\dot{x}}$; then  for 
any $g,x$, $\epsilon_x(g)=\overline{\dot{x}(g)}\,\overline{\epsilon_{\langle,\rangle}(g)}$
  and $g(\widetilde{x},\dot{x},\check{x})=(\mu^g\widetilde{x},\mu^{-g}\dot{x},\check{x})$
  as desired. 

Because $S_h\in \mathrm{Hom}(\alpha_h,\rho^2)$, we have 
\begin{equation}S_g^*\rho^2(S_0)S_g=S_g=S_0^*\rho^2(S_g)S_0=S_0^*\rho(U_g\rho(S_0)U_g^*)
S_0=(\rho(U_g)^*S_0)^*\rho^2(S_0)(\rho(U_g)^*S_0)\,.\nonumber\end{equation}
But the intertwiner space Hom$(\alpha_g,\rho^2)$ is one-dimensional, and thus
\begin{equation}
\beta(g)\,S_g^*=S_0^*\rho(U_g)\label{beta}\end{equation} for some scalars $\beta(g)$ 
with $|\beta(g)|=1$.
Because $\rho(U_{g+h})=\rho(U_g)\,\rho(U_h)$, we have from \eqref{beta} and 
$\alpha_h$-covariance that $\beta\in\widehat{G}$. The $S_h^*$-coefficient of \eqref{beta} reads
\begin{equation}
\beta(g)\,\delta_{g,h}=\frac{n}{\delta^2}\,\delta_{\mu^g,\mu_h}+\sum_{x,z}b_{z,x}\,\overline{b_{gz,x}}\,\widetilde{x}(h)\,.\label{ssrUa}\end{equation}
The triangle inequality applied to \eqref{ssrUa} with $g=h$ forces  $\beta=1=\overline{\mu^g}
\mu_g$ on $G$ (i.e. the  pairing $\langle *,*\rangle$ is symmetric).

Putting $g=h=0$ in \eqref{rTrSTa} gives $b'_{z,w}=s\delta\overline{a_{z}}\,
\overline{b_{\overline{z},w}}$. Hence the unitarity of $b'$ implies the left-side of \eqref{aunitar}.
$\rho(S_g^*\rho(T_z))T_x=\rho(S_g^*)T_x\rho(T_z)$ holds because $T_x\in\mathrm{Hom}
(\rho,\rho^2)$; its $S_0S^*_0$ and $T_yT^*_w$ coefficients,
together with $b'=s\delta\overline{ab}$, 
 gives the right-side of \eqref{aunitar} and the  left-side of \eqref{rSrTTb}.
Combining the left-side of \eqref{rSrTTb} with $a_{y,x}=sa_{x,y}$ and the unitarity of 
$\sqrt{\delta}b$,  gives 
\begin{equation} \sum_{z'}b_{z,z'}\,\overline{b''_{z';x,y,w}}=s \sum_{z'}b_{z',z}\,\overline{b''_{z';x,y,w}}\,.\label{bskew}\end{equation}
Choose $g\in G$ so that $\widetilde{w}\mu^g=1$; replacing $y$ with $\overline{gx}$ and
applying $\sum_xa_{gx}$ to \eqref{bskew}, this simplifies by \eqref{CuST} to $\sum_{z'}
b_{z,z'}\overline{b_{z',w}}=s\sum_{z'}b_{z',z}\overline{b_{z',w}}$, i.e. $b_{y,x}=sb_{x,y}$.
Substituting the value for $A(g)$ and our expression for $b'$ into \eqref{rSrSTb} with $h=0$,
and applying the triangle inequality and $\delta\sum_x|b_{x,y}|^2=1$, yields
$b_{gw,gy}=\epsilon_w(g)\epsilon_y(g)b_{w,y}$; comparing with covariance \eqref{ssrUaaa}
yields the selection rule for $b$ given in the theorem.

Because $T_w\in\mathrm{Hom}(\rho,\rho^2)$, we get
 $\rho(T_x^*\rho(S_g))T_w=\rho(T_x^*)T_w\rho(S_g)$. Using \eqref{rSrSTb} and our formula for
 $A(g)$, the  $T_{y}T_{z'}S_0S_0^*$ coefficient of this identity gives
\begin{equation}\overline{b''_{x;w,-g z',y}}\,a_{z'}=a_{gx}\,
b''_{\overline{g x};y,\overline{z'},-gw}\,\overline{\langle g,g\rangle}\,\epsilon_w(g)\,,\label{rTrSTc}\end{equation}
which simplifies to the second equality of \eqref{rSrTTb} (for $g=0$) and the right-side of 
\eqref{rhoalpTcc}. The left-side of \eqref{ssrUb} comes from \eqref{CuST}, while its right-side comes from
the $T_w^*T_z^*$ coefficient of  \eqref{beta}.
\eqref{Cu:TT}  follows by multiplying \eqref{Cu:TTa} by $\sum_{z,z'}b_{z,u}\overline{b_{z',v}}$ and using
the left-side of \eqref{rSrTTb} twice. \eqref{completddgg} is a simplification of \eqref{completdg}.
Equations \eqref{rTrTTb},\eqref{rTrTTc} simplify to
 the right-side of \eqref{rSrTTb} (after applying
the other parts of \eqref{rSrTTb}) and \eqref{rTrTTcc}, respectively.
\eqref{Cu:TT} follows directly from \eqref{Cu:TTa} and \eqref{rSrTTb}.
This concludes the proof of part (b).

To prove part (a), we need to verify the conditions of Proposition 1, given the equations
listed in Theorem 4. Unitarity of $b'$ follows from that of $\sqrt{\delta}a$ and $\sqrt{\delta}b$.
The covariances \eqref{ssrUaaa},\eqref{rhoalpTcc}, $a_{gx}=
\epsilon_{\langle,\rangle}(g)\dot{x}(g)a_x$ and $b'_{gz,hw}=\langle g,h\rangle{\widetilde{z}(h)}
\widetilde{w}(g)\epsilon_z(g)b'_{z,w}$ (covariance for $b'$ follows from that for $a$ and $b$)
are used repeatedly to simplify the expressions of Section 2.2. For example, these covariances
immediately give $W_{why}(g)={\langle g,h\rangle}\dot{y}(g)\epsilon_{\langle,\rangle}(g)\delta_{w,gy}$.

\eqref{CuST} follows from multiplying the left-side of \eqref{ssrUb} by $\sum_z\overline{b_{z',z}}$ 
and replacing $x$ with $gx$ and $y$ with $gy$. \eqref{completb} involves the right-side of
\eqref{ssrUb}. \eqref{rTrSTa} involves the selection rule $b'_{x,w}\ne 0\Rightarrow \dot{x}
\widetilde{x}=\dot{w}\widetilde{w}$. Verifying \eqref{rTrSTb} requires the selection rule for
$b_{\overline{x},w}$. \eqref{rTrTTb} involves all three identities in \eqref{rSrTTb}. To obtain
\eqref{rTrTTc} from \eqref{rTrTTcc}, use three times the left-side of \eqref{rhoalpTcc} with
$g=g_x'g_y$, as well as the selection rule for $b''$.

Equation \eqref{Cu:TT} and the left-side of \eqref{ssrUb} tells us that the inverse of the
matrix $b_{(w,y),(z,x)}:=b''_{z;w,x,y}$ is $b'_{(w,y),(z,x)}:=\overline{b''_{z;w,x,y}}-\overline{\dot{w}(g_x)}
\epsilon_{\langle,\rangle}(g_x)\delta_{w,g_xy}\delta_{z,x}$. Hence we obtain another form
of \eqref{Cu:TT}:
\begin{equation}\label{Cu:TTaa}
\sum_{z,x}b''_{z;w,x,y}\,\overline{b''_{z;w',x,y'}}=\delta_{w,w'}\delta_{y,y'}-\delta^{-1}\delta_{g_ww',
g_yy'}\delta_{g_ww,g_yy}\,.\end{equation}
 \eqref{Cu:TTa} arises from $\sum_{z,z'}\overline{b_{y',z}}b_{w',z'}$ applied to \eqref{Cu:TT},
 while \eqref{completdg} follows from \eqref{Cu:TTaa}. \eqref{Sintertwb} comes from
 \eqref{Cu:TT} and the formula for $\delta$. To see \eqref{rTrTTa}, hit it with $\sum_{w,y}\overline{
 b_{y',w}b_{z'',y}}$ and use \eqref{rSrTTb}, \eqref{Cu:TTa}, and both sides of  \eqref{aunitar}.
 From \eqref{completdg} we compute $V_{wxyz}(g)=\dot{z}(g)\epsilon_{\langle,\rangle}(g)
 \delta_{{w},gz}\delta_{g{x},y}$. \qquad \textit{QED to Theorem 4}\medskip





\section{Explicit classifications}

Recall that $G$ is the abelian group formed by the group-like simple objects, so $n=|G|$ is the number of $S$'s in the Cuntz algebra $\cO_{n,n'}$ of Section 2.
The fusion coefficient $n'=N_{\rho\rho}^\rho=\|\cF\|$ will be the number of $T$'s. We know
that either $n'=n-1$ (`first class') or $n'$ is a multiple of $n$ (second class).
In this subsection we explicitly solve our equations for small $n$ or $n'$.
But first we address the question of the direct relation of near-group systems to character rings
of groups $K$, a question begged by the examples in the Introduction.

\subsection{Which finite group module categories are near-group?}

We have seen several examples of finite groups $K$ whose module categories are
near-group. For example, 
the module categories Mod$(D_4)$ and Mod$(Q_8)$ for the order-8 dihedral and quaternion
groups, are both of type $\bbZ_2\times\bbZ_2+0$, and those examples apparently motivated Tambara-Yamagami to study their class of categories \cite{TY}. Similarly,
the even sectors of the $D_5^{(1)}$  subfactor satisfies the $S_3$ fusions and, more 
generally, the module categories for the affine groups Aff$_1(\bbF_q)$ are of type $\bbZ_{q-1}+q-2$.
The complete list of groups whose representation rings possess the near-group property has been rediscovered several times,
but perhaps originated with Seitz: 

\medskip\noindent\textbf{Proposition 3.} \cite{Sei} \textit{The complete list of all  finite groups $K$ 
whose module category Mod($K$) is a near-group category of type $G+n'$, for some
abelian group $G$ and some $n'\in\bbZ_{\ge 0}$, is:}

\smallskip\noindent\textbf{(a)} \textit{$|K|=2^k$ for $k$
odd, its centre is order 2, and $G/Z\cong \bbZ_2\times\cdots\bbZ_2$.
In this case, $G=K/Z(K)$, $d_\rho=\delta=2^{(k-1)/2}$, and $n'=0$.}

\smallskip\noindent\textbf{(b)} \textit{$
K\cong \mathrm{Aff}_1(\bbF_q)$ for some finite field $\bbF_q$. In this case, $G\cong
\bbZ_{q-1}$, $d_\rho=q-1$, and $n'=q-2$.}

\medskip The groups in part (a) are called \textit{extraspecial 2-groups}; there are
precisely 2 of them for each odd $k>1$. We will see next subsection that most C$^*$-categories 
of type $G+0$ are not Mod($K$) for some $K$. In contrast,
Proposition 5 below says that all but 5 C$^*$-categories of type $G+n-1$
will be Mod$(K)$ for $K$ in part (b). 

\subsection{The type $G+0$ classification}

As a special case of Theorem 4, we recover the Tambara-Yamagami classification \cite{TY}:

\medskip\noindent\textbf{Corollary 4.} \textit{The equivalence classes of $C^*$-categories of type $G+0$  are in one-to-one correspondence
with either choice of sign $s$ and any choice of nondegenerate symmetric pairing $\langle,\rangle$
on $G$, up to
automorphism of $G$.}\medskip

\noindent\textit{Proof.}  Because $n'=0$, the parameters $a,b,b',b''$ must be dropped from all equations in Proposition 1.
All that remains is the sign $s$ and the symmetric pairing $\langle*,*\rangle$, which will
be nondegenerate for $G$. \textit{QED to Corollary 4}\medskip

The proof in \cite{TY} is independent and much longer, involving a detailed study of the pentagon
equations in the category. It is worth remarking that \cite{TY} prove that $G$ must be abelian
(whereas we assume it).

These usually don't seem to be realised by a subfactor. 
As explained after Corollary 2, for both choices of signs the subfactors $\rho_\pm(M)\subset
M$ are equivalent to the $M^{G}\subset M$ subfactor.

\subsection{The near-group categories for the trivial group $G$}

It is generally believed that there are a
finite number of fusion categories of each rank, so in particular one would expect that for each finite $G$,
there are only finitely many near-group C$^*$-categories of type $G+n'$ for arbitrary $n'$. In fact, we
are led to expect that only finitely many cyclic groups $G=\bbZ_n$, when $n+1$ is not a prime power, 
will have C$^*$-categories of type $G+n'$ for $n'>0$.

Nevertheless, until we can bound $n'$ given a $G$, it seems to be nontrivial to
classify all near-group C$^*$-categories whose group-like objects form a given group $G$.
The only example we can fully work out is $G=\{0\}$ (although we expect the tube algebra
analysis for $n'>n$ should yield classifications for other groups of small orders).

\medskip\noindent\textbf{Proposition 4.} \textit{Up to equivalence, there are precisely 3 near-group $C^*$-categories of
type $\{0\}+n'$: namely, two of type $\{0\}+0$, and one of type $\{0\}+1$.}

\medskip The possibility $n'=0$ is Tambara-Yamagami and so is covered by Corollary 4,
while $n'=1$ is most easily handled using Corollary 5 below. 
The reason  there can be no examples with $n'>1$ is that such a solution
would yield a fusion category with rank 2 and a fusion coefficient $=n'\ge 2$, and no such
fusion category can exist \cite{Ost}. 

There are precisely 4 rank 2 fusion categories (2 of type $\{0\}+0$ and 2 of type $\{0\}+1$). The 
 one which is not realised as a fusion C$^*$-category is known as the Yang-Lee model,
 corresponding to one of the nonunitary Virasoro minimal models. The \textit{nonunitary} minimal
 models can never be realised as C$^*$-categories, so it is no surprise that Yang-Lee
 is missing from Proposition 4.

\subsection{The type $G+n-1$ C$^*$-categories}

By Corollary 3, the only parameters we need to identify are a sign $s$ when $n+1$ is
a power of 2, a third root of unity $\omega$ when $n+1$ is a power of 3, and the
$b''_{x,y}$ when $xy\ne 1$.
The complete classification for $n<32$ is collected in Table 1; the only value for the entries
$n\ne 1,2,3,7$ is the permutation $\sigma$. Recall $n+1$ must be a 
power of a prime and $G$ must be cyclic. In Table 1 we identify $\cF$ with the subset
$G\setminus\{0\}$.  $\omega$ there is any third root of 1.
In the $\sigma$ column of the table we use cycle notation, writing A for 10, B for 11, etc.
For $G=\bbZ_7$, $a(x)=1,1,s,1,s,s$ for $x=1,2,\ldots,6$ respectively, and $b''$ is given by
\begin{equation}
b''(i,j)=\left(\begin{matrix}1&s&s&1&s&*\\ 1&1&s&s&*&s\\ s&s&1&*&1&1\\ s&1&*&1&s&s\\ 
s&*&1&s&1&1\\ *&s&1&s&1&1\end{matrix}\right)\,.\end{equation}
It is elementary to verify from Corollary 3 that each entry in Table 1 yields an inequivalent
solution to the equations of Theorem 3.

$$\begin{array}{c||c|c|c|c|c|c} G
& \rm{\#} & \sigma&s & a&b&b''
\\ \hline\hline
\bbZ_1&2&${\scriptsize(1)}$&\pm 1&\emptyset&\emptyset&\emptyset\\ \hline
\bbZ_2& 3&${\scriptsize(1)}$&1&1&\omega&\emptyset\\ \hline
\bbZ_3& 2&${\scriptsize(1)}$&\pm 1&1,s&sa&b''_{xx}=a_x\\ \hline
\bbZ_4&1& ${\scriptsize(123)}$&1&1&1&1\\ \hline
\bbZ_6& 1&${\scriptsize(234)}$&1&1&1&1\\ \hline
\bbZ_7 & 2&${ \scriptsize(142)(365)}$& \pm 1&11s1ss&sa&\mathrm{(see\ above)}\\ \hline
\bbZ_8& 1 & ${\scriptsize(165)(273)}$& 1&1&1& 1\\ \hline
\bbZ_{10} & 1&${\scriptsize(159)(276)(843)}$ & 1& 1&1&1\\ \hline
\bbZ_{12} & 1& ${\scriptsize(16B)(378)(459)}$& 1& 1&1&1\\ \hline
\bbZ_{15} & 1& ${\scriptsize(1B3)(276)(4EC)(8D9)}$& 1& 1&1&1\\ \hline
\bbZ_{16}& 1 & ${\scriptsize(1AD)(2E8)(36F)(4B9)(5C7)}$& 1& 1&1&1\\ \hline
\bbZ_{18} & 1& ${\scriptsize(19H)(2EB)(4G7)(5CA)(6D8)}$& 1& 1&1&1\\ \hline
\bbZ_{22} & 1& ${\scriptsize(1EI)(236)(48L)(5FD)(7H9)(ACB)(GJK)}$& 1& 1&1&1\\ \hline
\bbZ_{24} & 1& ${\scriptsize(129)(3EJ)(5AL)(6CI)(7GD)(8HB)(FMN)}$& 1& 1&1&1\\ \hline
\bbZ_{26} & 1& ${\scriptsize(139)(2LG)(4FK)(5OA)(6BM)(7EI)(8CJ)(HNP)}$& 1& 1&1&1\\ \hline
\bbZ_{28} & 1& ${\scriptsize(1ER)(293)(4FN)(5DO)(6KG)(7HI)(8MC)(ABL)(JQP)}$& 1&1&1&1\\ \hline
\bbZ_{30} &1 &${\scriptsize(1LN)(23A)(4JM)(6OF)(79T)(8BQ)(CHG)(DIE)(KRS)}$ & 1&1&1&1\\ \hline
\bbZ_{31} & 1& ${\scriptsize(1HD)(23Q)(46L)(5ST)(7F9)(8CB)(APR)(EUI)(GOM)(JNK)}$& 1&1&1&1\\ \hline
\end{array}$$

\centerline{\textbf{Table 1.} The C$^*$-categories of type $G+n'$ for $|G|=n'+1\le 31$}\medskip

\noindent\textbf{Proposition 5.} \textit{If a $C^*$-category is of type $G+n'$ with $n'\not\in n\bbZ$,
then $G=\bbZ_n$, $n+1=:q$ is a prime power, and $n'=n-1$. There is precisely one $C^*$-category 
of type $\bbZ_n+n-1$, namely Mod(Aff$_1(\bbF_q)$), except for $n=1,2,3,7$ which have
precisely 1,2,1,1 additional $C^*$-categories, collected in Table 1.}\medskip

\noindent\textit{Proof.} We know $n'=n-1$ from Theorem 2, and $G\cong\bbZ_{n}$ where $n=q-1$ 
for some prime power $q=p^k$, by Proposition 2. Corollary 7.4 of \cite{EGO} tells us that 
the only fusion categories of type $\bbZ_n+n-1$ are 
Mod(Aff$_1(\bbF_q)$), except for $n=1,2,3,7$ where there are
precisely 1,2,1,1 additional fusion categories. There always is at least one C$^*$-category
of type $\bbZ_n+n-1$, namely the one corresponding to the solution $b=b''=1$, so it suffices
to find 1,2,1,1 additional solutions the the equations of Theorem 3, when $n=1,2,3,7$ respectively.
These are collected in Table 1. \textit{QED}\medskip

 Combining \cite{EGO} and Proposition 5,
 we find that each near-group fusion category with $n'=n-1$ and $G=\bbZ_n$ has a 
 C$^*$-category structure, i.e. a system of endomorphisms (unique up to
equivalence). To our knowledge, this gives the first construction of the extra C$^*$-categories
for $\bbZ_3$ and $\bbZ_7$. As pointed out
in Section 2.3, the collection of (equivalence classes of) type $G+n-1$ C$^*$-categories
for a given $G$ will form an abelian group. We find that this group is always
trivial, except for $n=1,2,3,7$ when it is $\bbZ_2,\bbZ_3,\bbZ_2,\bbZ_2$ respectively.

Recall the discussion of  the deformation parameters $H^2((\bbZ_p^k,\bbZ_{p^k-1});\bbT)/\!\sim$ 
near the end of Section 2.3. We expect the group $H^2((\bbZ_p^k,\bbZ_{p^k-1});\bbT)$ is trivial for 
all $n$.
It should be possible to verify this from the results of \cite{IK}, at least  
 for $n\le 4$ and $n=6$. Certainly it says there is a unique depth-2 subfactor with
principal even fusions of type $\bbZ_n+n-1$ for those $n$, namely $M^H\subset M\sdprod \bbZ_{n-1}$
for $H=\bbZ_p^k$ (see \cite{BiHa} for a complete analysis of these subfactors). 
 It is easy to compute $H^2((\bbZ_p^k,\bbZ_{p^k-1});\bbT)$ 
directly from the definition, for $n=1,2$, and we find indeed that it is trivial.
In particular, this means that only the $s=1$ solution at $n=1$, and only $\omega=1$ at $n=2$,
are realised by depth-2 subfactors. The triviality of this for $n=5$ follows from 
uniqueness results for subfactors of index 5.

\subsection{At least as many $T$'s as $S$'s}




Consider now type $G+n'$, where $n'\ge n$ (and therefore must be a multiple of $n$). 
We don't know of any examples where $n'>n$. A natural approach to bounding $n'$,
given $n$, would be carrying through the tube algebra analysis for $n'>n$. After all,
we find in Section 4.2 below that this strategy is effective in pruning the possibilities
for type $G+n-1$. Likewise, Ostrik's analysis \cite{Ost}, which eliminated $n'>n$ for $G=\{0\}$,
investigated the modular data for the double. We will pursue this thought in future work.

\medskip\noindent\textbf{Corollary 5.} \textit{Any solution to the equations of Proposition 1, for arbitrary
$G$ and $n'=n$, has $s=1$, $\cF=\widehat{G}$, $H=0$, $u_{x,g}=1$, and is of the form 
$a_{x,y}=\sqrt{\delta}^{-1}\delta_{x,\overline{y}}\,a(x)$, $b_{z,x}=c\sqrt{n\delta}^{-1}
\langle z,x\rangle$, $b'_{z,x}=\overline{a(z)c}\sqrt{n}^{-1}\langle z,x\rangle$,
$b''_{z;w,x,y}=\delta_{y,wx}\,a(x)\,b(z\overline{x})\langle z,y\rangle$ for some complex numbers
$c,a(x),b(x)$ satisfying}
\begin{align}a(1)=1\,,\ a(x)=a(\overline{x})\,,\ a(xy)\,\langle x,y\rangle=a(x)\,a(y)\,,\  \sum_x{a}(x)=
\sqrt{n}\,c^{-3}\,,&&\label{5.1}\\
b(0)=-1/\delta\,,\ \sum_y\overline{\langle x,y\rangle}\,{b}(y)=\sqrt{n}\,c\,\overline{b(x)}\,,\ 
a(x)\,b(\overline{x})=\overline{b(x)}\,,&&\label{5.2}\\ \sum_{x}b(xy)\,\overline{b(x)}=\delta_{y,1}-{\delta}^{-1}
\, ,\ \sum_{x}b(xy)\,b(xz)\,\overline{b(x)}=\overline{\langle y,z\rangle}\,b(y)\,b(z)-\frac{c}{\delta\sqrt{n}}\,.&&\label{5.4}\end{align}
\textit{Conversely, any $a(x),b(x),c$ satisfying \eqref{5.1}-\eqref{5.4} yields a solution to the equations of Proposition 1
in this way. Two $C^*$-categories $\cC_1,\cC_2$ of type $G+n$ are equivalent iff $c_1=c_2$
and there is a $\phi\in\mathrm{Aut}(G)$ such that $\langle g,h\rangle_2=\langle \phi g,\phi h\rangle_1$, $a_2(x)=a_1(\phi x)$, and $b_2(x)=b_1(\phi x)$.}

\medskip\noindent\textit{Proof.} Through the pairing $\langle*,*
\rangle$ we may identify $\widehat{G}$, and hence $\cF$, with the group $G$. More precisely, let $x_0$ denote the $x\in\cF$ with $\widetilde{x_0}=1$;
then $x_g=gx_0$ has $\widetilde{gx_0}(h)=\mu^g(h)=\langle g,h\rangle$.
So the action of $G$ on $\cF$ corresponds to addition in $G$: $hx_g=x_{h+g}$.

 Equivariance \eqref{ssrUaaa} forces $b_{x_g,x_h}=c
\sqrt{\delta}^{-1}\langle g,h\rangle$ for some $c$ with $|c|=1$. Write $a_{x_g,x_h}=\sqrt{\delta}^{-1}
a(g)\delta_{g,-h}$ for some numbers $a(g)$ with $|a(g)|=1$. Then $a(-g)=sa(g)$, so we get
$s=+1$ by looking at $g=0$. We get $b'_{x_g,x_h}=\overline{c}\,\langle g,h\rangle\overline{a(g)}$
and covariance \eqref{rhoalpTcc} says $b''_{x_g;x_h,x_k,x_l}=\delta_{l,h+k}\langle g,h\rangle
b''_{g,k}$ for some numbers $b''_{g,k}$. The left-side of \eqref{rSrTTb} now reduces to
\begin{equation}
c\sum_g\langle g,h-k\rangle\,\overline{b''_{g,l}}=\overline{a(k)}\,\overline{\langle h,k\rangle}\,b''_{h,k+l}\end{equation}
and hence 
\begin{equation}b''_{h,k}=c\,a(k)\,\langle h,k\rangle\sum_g\langle g,h-k\rangle\,\overline{b''_{g,0}}
=a(k)\,\langle h,k\rangle\,b''_{h-k,0}\,.\end{equation}
Writing $b''_{g,0}=b(g)$, it is now easy to verify that all equations of \cite{iz3} are recovered.
\textit{QED to Corollary 5}\medskip

This Corollary says that the case $n'=n$ reduces to the generalisation of $E_6$ introduced
in Section 5 of \cite{iz3}.

\medskip
\noindent\textbf{Proposition 6.} \textit{There are  (up to equivalence) precisely 1, 2, 2, 2, 3, 4, 2, 8, 2, 4, 4, 4, 4 systems for $G=\bbZ_n$, $1\le n\le 13$, respectively. There is 1 solution each for
$G=\bbZ_2\times\bbZ_2$ and  $G=\bbZ_3\times\bbZ_3$, 2 solutions for $\bbZ_2\times \bbZ_6$, 4 solutions for $\bbZ_2\times\bbZ_4$, and no solutions for  $\bbZ_2\times
 \bbZ_2\times\bbZ_2$.}\medskip

{\scriptsize $$\begin{array}{c||c|c|c|c|c|c|c} G
& c&\langle,\rangle & a&Q'&(j_1,\ldots,j_{\lfloor n/2\rfloor})
\\ \hline\hline
\bbZ_1&1&1&1&2&\emptyset\\ \hline
\bbZ_2&\xi_{24}^{17}& 1&1&(-1,-5)&( 0.78539816)\\ \hline
&\xi_{24}^{7}& 1&-1&(1,5)&(-0.78539816)\\ \hline
\bbZ_3&\xi_{12}^{-1}& 1&1&1&(-2.8484536)\\ \hline
&\xi_{12}&-1& 1&-1&(2.8484536)\\ \hline
\bbZ_4&\xi_8^{3}&1&1& (3,-3)&( -0.60623837,-1.5707963)\\ \hline
&\xi_8^{-3}&-1& 1&(-3,3)&(-1.39163653,1.57079632)\\ \hline
\bbZ_2\times\bbZ_2&1& 1&1&&(-2.356194490,2.356194490,0)\\ \hline
\bbZ_5 &-1& 1&1&(1,1)&(-1.256637,1.256637)\\ \hline
&\xi_3&2 &1&2&(-1.0071249,0.3425266)\\ \hline
&\xi_{3}^{-1}&2 &1&-2&(-0.3425266,-2.263762)\\ \hline
\bbZ_6&\xi_{24}^{-1}& 1 & 1&(-1,-19)&( -2.9552611, -0.055354168, 0.78539816)\\ \hline
&\xi_{24}^{-5}& -1 &1 &(3,13)&( -2.915033694, 1.5909100, -2.3561944)\\ \hline
&\xi^{5}_{24}& 1 & -1&(-3,-13)&( 2.9150336,-1.590910, 2.3561944)\\ \hline
&\xi_{24}& -1 &-1 &(1,19)&(2.9552611, 0.055354168, -0.78539816)\\ \hline
\bbZ_{7} &\i& 1&1 &-2&(-1.05169,1.7936250,-0.3143315)\\ \hline
&-\i& -1&1 &2&( 1.0516925,-1.793625,0.31433)\\ \hline
\bbZ_{8} &\xi_{24}& 1&1&(-5,-7)& (-0.87227636,2.7042615,-2.9767963,3.1415926)\\ \hline
&\xi_{24}& 1&-1& (-5,17)&(-2.9767963,-1.1334651, -0.87227635,3.1415926)\\ \hline
&\xi_{24}^{-1}& -1&1&(5,7) &( 0.87227636,-2.7042615,2.9767963,3.1415926)\\ \hline
&\xi_{24}^{-1}& -1&-1&(5,-17) &( 2.9767963,1.1334651,0.87227635,3.1415926)\\ \hline
&\xi_{24}^{7}& 3&1&(1,11)&( 2.4640490,-3.0755747,-0.49188699,0)\\ \hline
&\xi_{24}^{7}& 3&-1& (1,-13)&(-0.49188700,1.5047784,2.4640490,0)\\ \hline
&\xi_{24}^{17}& -3&1&(-1,-11) &(-2.4640490,3.0755747,0.49188699,0)\\ \hline
&\xi_{24}^{17}& -3&-1& (-1,13)&(0.49188700, -1.5047784,-2.4640490,0)\\ \hline
\bbZ_{2}\times\bbZ_4&\xi_{12}^5 & 1&(1,1)&& (0.7853981,1.77783,-2.497219,1.570796,-0.7853981)\\ \hline
&\xi_{12}^{-5} & 1&(-1,1)& &(-0.785398,0.9924406,1.42977,-1.57079,0.7853981)\\ \hline
&\xi_{12}^{-5} & 1&(1,-1)& &(0.785398, 1.42977,-1.777838,-1.570796,-0.785398)\\ \hline
&\xi_{12}^5 & -1&(-1,-1)& &(-0.7853981,-2.497219,-0.9924406,1.5707963,0.7853981)\\ \hline
\bbZ_{9}&\xi_3^{-1}& 1 & 1&-2&( -2.69568,1.367012,1.41882,-2.38374)\\ \hline
&\xi_3& -1 & 1&2&(2.695680,-1.3670127,-1.418824,2.383744)\\ \hline
\bbZ_3\times\bbZ_3&1&-1&1&2&(2.9557793,-1.2330109,-2.2802084,0)\\ \hline
\bbZ_{10} &\xi_8^{-1}& 1&1 &&(-2.3665026,-3.0894639,3.077894,0.00650245,0.785398)\\ \hline 
&\xi_8^{3}& 7&1 &&(1.7756309,-0.6115079,-1.030618,2.8686859,-2.3561944)\\ \hline 
&\xi_8& 1&-1 &&(-3.077894,2.519776,-1.424024,3.089463,-0.78539816)\\ \hline 
&\xi_8^{-3}& 7&-1 &&(-1.3447773,-2.868685,-1.7756309,0.64512913,2.3561944)\\ \hline
\bbZ_{11} &\xi_{12}^7& 1& 1&2&( 1.9464713,2.0140743,-1.7487929,0.3352432,-0.1427077)\\ \hline
 &\xi_{12}^{-1}& 1& 1&1&(0.53877136,-2.8317431,0.2827610,0.46457259,2.5063157)\\ \hline
 &\xi_{12}& -1& 1&-1&(-2.8884206, 2.3090448,0.85395967,2.1781685,-1.4920749)\\ \hline
 &\xi_{12}^{5}& -1& 1&-2&(-1.4807206,0.87167704,-1.1775942,2.0488391,2.1420869)\\ \hline
\bbZ_{12} &\xi_{24}^7& 1&1 &&(3.0822445,-0.3494640,-3.0450322,-0.7241984,-0.38234715,1.570796)\\ \hline
 &\xi_{24}^{-5}& 1&-1 &&(-0.6247574,-3.044463,-2.3415376,0.4718634,0.99777419,1.5707963)\\ \hline
 &\xi_{24}^{-7}& -1&1 &&(-3.0822445,0.34946402,3.0450322,0.7241984,0.3823471,-1.570796)\\ \hline
&\xi_{24}^{5}& -1&-1 &&(0.6247574,3.0444636,2.3415376,-0.47186343,-0.99777419,-1.5707963)\\ \hline
\bbZ_{2}\times\bbZ_6&\xi_3^{-1} & -1&(1,1)& &( -2.35619,0.0611997,2.469129,0.89833,-1.88433,0.785398,1.57079)\\ \hline
&\xi_3 & 1& (-1,1)&&(2.356194,-0.0611997,-2.4691295,-0.8983332,1.884331,-0.785398,-1.570796)\\ \hline
\bbZ_{13} &-1& 1& 1&3&(-2.4521656,1.9847836,0.42579608, 1.4322079,-1.4550587,1.1404478)\\ \hline
&-1& 1& 1&3&(1.4550587,1.3924399,-1.9847836,-1.2761619,0.44776608,-1.4322079)\\ \hline
&1& 2& 1&1&(-2.4805730,3.0305492, 0.28372451,-0.04125417, 0.44928247,2.9410534)\\ \hline
&1& 2& 1&1&( -0.44928247,2.2170122,-3.0305492,-1.892166, -2.9638949,0.041254182)\\ \hline
\end{array}$$}

\centerline{\textbf{Table 2.} The C$^*$-categories of type $G+n$ for $|G|\le 13$}\medskip

All systems for $G=\bbZ_n$ ($n\le 4$) and $\bbZ_2\times\bbZ_2$, and the first one for $\bbZ_5$ were 
constructed in Section 5 and Appendix A of \cite{iz3}; the rest are to our knowledge new.
The column $Q'$ will be relevant to subsection 4.4, where it will be explained.

In the table, we write $\xi_k$ for $\exp(2\pi\i/k)$.
If $\langle,\rangle$ is a nondegenerate symmetric pairing for a cyclic group $G=\bbZ_n$,
then it equals $\langle g,h\rangle=\exp(2\pi\i m gh/n)$ for some integer $m$ coprime to $n$.
When $G=\bbZ_{n'}\times\bbZ_{n''}$ in Table 2, any pairing appearing there is of the
form $\langle (g',g''),(h',h'')\rangle=\exp(2\pi\i g'h'/n')\,\exp(2\pi\i m g''h''/n'')$ for some $m\in\bbZ$
coprime to $n''$. This $m$ is the entry appearing in the $\langle,\rangle$ column. 

Note that if $a_i$ are both solutions of \eqref{5.1} for fixed group $G$ and pairing $\langle,\rangle$, then $a_2=\psi\,a_1$ for some $\psi\in \widehat{G}$ with $\psi^2=1$. Thus for $G$ of
odd order, the unique $a$ is $a(g)=\langle g,g\rangle^{(n-1)/2}$. For $G=\bbZ_n$ or $G=\bbZ_{n'}
\times \bbZ_{n''}$, for $n,n',n''$ even,
with pairing given by $m$ as above, then $a(g)=s_1^g\exp(-\pi m g^2/n)$ or
$a(g',g'')=s_1^{g'}s_2^{g''}$ respectively, for some $s_1,s_2\in\{\pm 1\}$. These signs
$s_1$ or $(s_1,s_2)$ grace the fourth column of the Table.

For the group $G=\bbZ_n$, the quantities $b(g)$ are recovered through  Table 2 from the formula $b(g)=e^{\i j_g}/\sqrt{n}$,
for $0<g\le n/2$, $b(0)=-1/\delta$, and $b(-g)=\overline{a(g)}\,\overline{b(g)}$.
For the noncyclic group $\bbZ_{n'}\times\bbZ_{n''}$, these parameters $j$ are taken in 
order $(1,0),\ldots,(\lfloor n'/2\rfloor,0),(0,1),(1,1),\ldots,(\lfloor n'/2\rfloor,\lfloor n''/2\rfloor)$.

The table lists representatives of equivalence classes of systems. Using Corollary 5, 
it is easy to determine when numerical invariants determine equivalent
systems. Note that taking the complex conjugate of numerical invariants (i.e. the conjugate of
$c$, and the negatives of the $\langle,\rangle$ and $j_i$ columns) will yield another
(possibly equivalent) solution to the equations of Corollary 5. For example, consider the
first entry for $G=\bbZ_5$: although the complex conjugate solution  has different $j$'s, it is equivalent as 
the $j$'s are permuted back to each other through the $\bbZ_5$ automorphism $-1\in\bbZ_5^\times$.
On the other hand, complex conjugation interchanges the second and third systems
for $\bbZ_5$; these two are inequivalent because they have different $c$'s.

Apart from $G=\bbZ_2$, the solutions in the Table turn out to be precisely the solutions
to the linear equations \eqref{5.2} together with $|b(g)|=1/\sqrt{n}$ for $g\ne 1$ (which is a
consequence of the equations of Corollary 5). 
These values for $j_g$ are floating point; to improve their accuracy arbitrarily is trivial
using mathematics packages like Maple (where you would change `Digits' to say 200,
and use `fsolve' with the provided seed values). The $b(g)$'s are in fact all algebraic,
but providing the exact algebraic expressions (though possible) would not be very enlightening.

We illustrate our method of establishing Proposition 6 and Table 2, with the hardest case, namely
$n=13$. Up to automorphism of $\bbZ_{13}$, there are two possible pairings: $\langle g,h\rangle
=\xi_{13}^{mgh}$ for $m=1$ or 2; choose $m=1$ (this will also yield the solutions for $m=2$).
By \eqref{5.1} there are 3 possible values of $c$; choose $c=-1$ for now. To find a complete list of candidate solutions,  
first solve the linear equations \eqref{5.2} (breaking each $b(g)$
into its real and imaginary parts). This determines each of the 26 variables Re$\,b(g)$, Im$\,b(g)$, 
up to 4 real parameters, which we can take to be the real and imaginary parts of $b(11)$
and $b(12)$. The norms $|b(1)|^2=\cdots=|b(4)|^2=1/n$ yield four independent quadratic 
identities obeyed by those parameters; by Bezout's Theorem they can have at most
 $2^k$ (complex)  solutions, and as always  here
this upper bound is realised (i.e. no zeros have multiplicities). We `solved' these equations
using floating point: our approximate solutions are
{\scriptsize\begin{align} && {\scriptsize(\mathrm{Re}\,b(12),\mathrm{Im}\,b(12),
\mathrm{Re}\,b(11),\mathrm{Im}\,b(11))\in\{( -.0277709, .9996143,  .7988155,   -.601576), ( -.022834, -.999739, .585511,-.810664),}\nonumber\\
&&{(.1154793,.9933099,.1774122,.9841366),
 (  -.2122594, .9772133,   -.0413084, .9991464),(  -.3498383, .9368100,   .9107102, -.4130459), }\nonumber\\
&&{  ( .623311,  .7819739, .7557187,  -.6548963),
 ( -.7923374,  .6100831,  -.1096719,  .9939678),
( .8581882,  .5133351,  .8230236,  -.5680071),}\nonumber\\
&&{( -.9560999,  .2930407, .00007002229,  .9999999),
( .976436,-.2158069,   .6983522,  -.7157542), (.9014165, -.4329529, .5249502,   -.8511329),}
\nonumber\\ &&{(  -.8964165,  -.4432124, -.1923449,  .9813274),
( -.7716105,  -.6360951,-.402263,  .9155241),
(  .2614237,   -.9652241, -.3345549, .9423762),}\nonumber\\
&&{ ( .1595449, -.9871906,  -.5115023,  .8592819),
(.08134136,- .9966863,.4166086,- .9090859)}\}\,.\nonumber\end{align}}
Of these, four  also satisfy the remaining norms $|b(5)|^2=|b(6)|^2=1$, and indeed all equations
\eqref{5.4}. Two of these are related to the other two by the automorphism $-1\in\bbZ_{13}^\times$; 
the resulting  two inequivalent
solutions are collected in
Table 2. Incidentally, had we chosen either of the two remaining possibilities for $c$, we would
have had 3 parameters from the linear system, 8 solutions to 3 norm equations, but none of these
8 would be solutions to all of the remaining 3 independent norm equations; thus those
other values of $c$ don't yield solutions (floating point calculations suffice, thanks to Bezout).

We still need to verify that each of these are indeed  solutions, i.e. that they satisfy \eqref{5.1}-\eqref{5.4}
exactly. For this purpose, define $\sqrt{13}b^{(1)}(g)$ to be the roots of
{\scriptsize \begin{align}
&&P(X)=X^{12}+\frac{13+\sqrt{13}-\sqrt{17}-\sqrt{221}}{2}X^{11}+
\frac{53+8\sqrt{13}-13\sqrt{17}-7\sqrt{221}}{2}X^{10}+
\frac{637+67\sqrt{13}-59\sqrt{17}-43\sqrt{221}}{2}X^9\nonumber\\&&+
\left(\frac{7531}{4}+\frac{1905\sqrt{13}}{4}-416\sqrt{17}-\frac{253\sqrt{221}}{2}+
\frac{t}{628}\left(\frac{57473}{29}+\frac{15935\sqrt{13}}{29}-\frac{16341\sqrt{17}}{34}-
\frac{4531\sqrt{221}}{34}\right)\right)X^8+\nonumber\\&&
\frac{169}{8}\left(-16809-2587\sqrt{13}+2249\sqrt{17}+1127\sqrt{221}+t
\left(\frac{3812952}{4553}+\frac{1057407\sqrt{13}}{4553}-\frac{15721329\sqrt{17}}{77401}
-\frac{4359784\sqrt{221}}{77401}\right)\right)X^7\nonumber\\&&
+\frac{1}{2}\left(- 25283 -7179\sqrt{13}+6279\sqrt{17}+1701\sqrt{221}\right)X^6+\nonumber\\&&
\frac{169}{8}\left(-16809-2587\sqrt{13}+2249\sqrt{17}+1127\sqrt{221}-t
\left(\frac{3812952}{4553}+\frac{1057407\sqrt{13}}{4553}-\frac{15721329\sqrt{17}}{77401}
-\frac{4359784\sqrt{221}}{77401}\right)\right)X^5\nonumber\\&&
+\left(\frac{7531}{4}+\frac{1905\sqrt{13}}{4}-416\sqrt{17}-\frac{253\sqrt{221}}{2}-
\frac{t}{628}\left(\frac{57473}{29}+\frac{15935\sqrt{13}}{29}-\frac{16341\sqrt{17}}{34}-\frac{4531\sqrt{221}}{34}\right)\right)X^4
\nonumber\\&&+\frac{637+67\sqrt{13}-59\sqrt{17}-43\sqrt{221}}{2}X^3+\frac{53+8\sqrt{13}-13\sqrt{17}-7\sqrt{221}}{2}X^2+\frac{13+\sqrt{13}-\sqrt{17}-\sqrt{221}}{2}X+1
\nonumber\end{align}}
for $t=\i\sqrt{75090+2\sqrt{13}}$. The roots $\sqrt{13}b(g)$ are thus algebraic integers, and 
necessarily have modulus 1. 
 The base field (i.e. the one generated over $\bbQ$ by the
coefficients) of $P(X)$ is clearly $\bbK=\bbQ[\sqrt{13},\sqrt{17},t]$, which has Galois
group (over $\bbQ$) $D_4$ generated by $t\mapsto s\i\sqrt{75090+s'2\sqrt{13}}$ and
$\sqrt{17}\mapsto s''\sqrt{17}$, for all signs $s,s',s''$. The splitting field $\bbF$ of $P(X)$
is a quadratic extension of $\bbQ[\xi_{13},\sqrt{17},t]$, with Galois group $\bbZ_{12}$ over $\bbK$ given by $b(j)\mapsto b(\ell j)$, $\xi_{13}\mapsto \xi_{13}^{\ell^2}$, for any
$\ell\in\bbZ_{13}^\times \cong\bbZ_{12}$. (Use resultants to prove the products of roots of $P(X)$
include 13th roots of unity, and hence that $13b(g)b(-g)=\overline{a(g)}$; the first prime $p>17$
with $\sqrt{13},\sqrt{17},\sqrt{-75090-2\sqrt{13}}\in\bbZ_p$ is $p=101$, and modulo 101
\begin{equation}P(X)\equiv  
19 + 74 X+55 X^2 + 62 X^3  + 74 X^4  + 77 X^5  + 39 X^6  + 86 X^7  + 8 X^8  + 47 X^9  + 73 X^{10}  + 61 X^{11}   + X^{12}\nonumber\end{equation}
is irreducible in $\bbZ_{101}[X]$, hence Gal$(\bbF/\bbK)$ contains an element of
order 12; etc.) The 7 nontrivial automorphisms in Gal$(\bbK/\bbQ)$, lifted to $\bbF$, map $b(g)$ to the three other solutions for $n=13$ in Table 2, together with 4 analogous `shadow' solutions
corresponding to $\delta'=(13-\sqrt{221})/2$, which can also be estimated numerically. To show that
$b(g)$ (as well as the other 3 candidates in Table 2 for $n=13$) indeed satisfy the remaining
identities in \eqref{5.1}-\eqref{5.4}, it suffices to replace each $\overline{b(g)}$ with
$13/b(g)$, multiply by $\delta$ and an appropriate power of $\sqrt{n}$ to guarantee that
the equations are manifestly algebraic integers, and then evaluate the equations numerically 
for all 8 choices of $b$ (the 4 from Table 2 and the 4 shadows). We used 200 digits of
accuracy --- far more than necessary but trivial using Maple --- and found that the equations
held to accuracy $10^{-190}$ or so. The errors will therefore be algebraic integers, and from the
above we know that all of their Galois associates will have 
modulus $<<1$. This means the errors must vanish identically, and we are done.
(Incidentally, the polynomial $P(X)$ was found working backwards: the numerical analysis
of the previous paragraph suggested its existence and basic properties.)

The 1 in Proposition 6 for $G=\bbZ_1$ corresponds to (i.e. is also implied by) the uniqueness of the $A_3$ subfactor, and the 2 systems
for $G=\bbZ_2$ correspond to the two versions of the $E_6$ subfactor.
Note that our classification for uniqueness (up to complex conjugation) for $G=\bbZ_3$  corresponds 
to the uniqueness for even sectors of the Izumi-Xu 2221 subfactor. The uniqueness of the Izumi-Xu subfactor was
first shown
in the thesis of Han \cite{Hn}. His proof is independent of ours: it involved planar algebras,
and was quite complicated.

Note the numerology $n\delta_{Hn}=\delta_{IXn^2}$, where
$(n+\sqrt{n^2+4})/2$ is the dimension of the nongrouplike simple objects in the Haagerup
system for $\bbZ_n$ \cite{iz3}, and $\delta_{IXn^2}$ is the dimension of $\rho$ for
near-group C$^*$-categories of type $\bbZ_{n^2}+n^2$. This suggests comparing the subfactor
$\rho_{Hn}(M)\subset M$ with the subfactor $\rho_{IXn^2}(M)\subset M^{\bbZ_{n^2}}$,
as they have the same index, namely $\delta_{Hn}^2$. For $n=3$, the principal
graph of the former is Figure (4') in Lemma 3.10 of \cite{GrSn}, while the principal graph
of the latter is the completely different $2^91$, so the connection (if indeed there is one) is not simply this.


\section{Tube algebras and modular data}

\subsection{The tube algebras of near-group systems}

In this section we compute the tube algebras, for any solution to the equations of Proposition 1.
Our notation will be as in \cite{iz1}. 
We can assume $\cF\ne\emptyset$, i.e. $n'\ne 0$, as the tube algebra for the Tambara-Yamagami
systems was computed in Section 3 of \cite{iz3}.

Let $\Delta=\{\alpha_g,\rho\}_{g\in G}$ be a 
finite system of endomorphisms, as in Section  2. The tube algebra 
Tube $\Delta$ is a finite-dimensional C$^*$-algebra, defined as a vector space by
\begin{equation}\mathrm{Tube}\,\Delta=\oplus_{\xi,\eta,\zeta\in\Delta}\mathrm{Hom}(\xi\cdot\zeta,
\zeta\cdot\eta)\,.\end{equation}
Given an element $X$ of Tube$\,\Delta$, we write $(\xi\zeta|X|\zeta\eta)$ for the restriction to
Hom$(\xi\zeta|X|\zeta\eta)$, 
since the same operator may belong to two distinct intertwiner spaces. 
For readability we will often write $g$ for  $\alpha_g$.
In our case the intertwiner spaces are Hom$(\alpha_{g+h},\alpha_g\alpha_h)=\bbC 1$,
Hom$(\rho,\alpha_g\rho)=\bbC 1$, Hom$(\rho,\rho\alpha_g)=\bbC U_g$, Hom$(\alpha_g,
\rho^2)=\bbC S_g$ and Hom$(\rho,\rho^2)=\mathrm{span}_{z\in\cF}\{T_z\}$. 
Denote the elements
$A_{g,h}=(gh|1|hg)$, $B_{g,h}=(g\rho|U_h|\rho h)$, $C_{g,z}=(g\rho|T_z|\rho\rho)$,
$D_{g,z}=(\rho\rho|U_gT_z^*|\rho g)$, $E_k=(\rho k|U_k^*|k\rho)$, $E'_k=(\rho\rho|S_kS_k^*|
\rho\rho)$, $E_{wz}''=(\rho\rho|T_wT_z^*|\rho\rho)$. Then the
vector space structure of the tube algebra is:
\begin{equation}
\mathrm{Tube}\,\Delta=\bigoplus_{g,h}\cA_{g,h}\oplus\bigoplus_g\cA_{g,\rho}\oplus\bigoplus_g
\cA_{\rho,g}\oplus\cA_{\rho,\rho}\,,\end{equation}
where $\cA_{g,g}=\bbC B_{g,g}\oplus\mathrm{span}_{k\in G}A_{g,k}$,
$\cA_{g,h}=\bbC B_{g,h}$ (for $g\ne h$), $\cA_{g,\rho}=\mathrm{span}_{z\in\cF}
C_{g,z}$, $\cA_{\rho,g}=\mathrm{span}_{z}D_{g,z}$,
$\cA_{\rho,\rho}=\mathrm{span}_{k}E_k\oplus\mathrm{span}_{k}E'_k\oplus\mathrm{span}_{w,z}E''_{w,z}$.

The C$^*$-algebra structure of Tube $\Delta$ is as follows: multiplication is given by
\begin{equation}
(\xi\zeta|X|\zeta\eta)(\xi'\zeta'|Y|\zeta'\eta')=\delta_{\eta,\xi'}\sum_{\nu\prec\zeta\zeta'}
\sum_i(\xi\nu|T(\nu,i)^*\rho_\zeta(Y)X\rho_\xi(T(\nu,i)|\nu\eta'))\,,\label{mult}\end{equation}
where $\rho_\rho=\rho$ and $\rho_g=\alpha_g$, and: when $\zeta=g$ and $\zeta'=h$,
then the unique $\nu$ is $g+h$ and the unique  $T(\nu,i)$ is $1$;
when $\zeta=g$ and $\zeta'=\rho$, the unique $\nu$ is $\rho$ and the unique $T(\nu,i)$ is
1; when $\zeta=\rho$ and $\zeta'=g$, the unique $\nu=\rho$ and the unique $T(\nu,i)$ is
$U_g$; and when $\zeta=\zeta'=\rho$, then $\nu$ runs over all $g\in G$, with $T(g,i)=S_g$,
as well as $\nu=\rho$, with $T(g,i)$ running over all $T_z$. Moreover, the adjoint is
\begin{equation}
(\xi\zeta|X|\zeta\eta)^*=d_{\zeta}(\eta\overline{\zeta}|\rho_{\overline{\zeta}}(\rho_\xi(\overline{R}_\zeta^*)X^*)R_\zeta|\overline{\zeta}\xi)\,,\label{star}\end{equation}
where $d_g=1$, $d_\rho=\delta$, $R_g=1=\overline{R}_g$, and $R_\rho=S_0$, $\overline{R}_\rho
=sS_0$.

Let $\sigma$ be a finite sum of sectors in $\Delta$. A half-braiding for $\sigma$ is a choice of
unitary operator $\cE_\sigma(\xi)\in\mathrm{Hom}(\sigma\rho_\xi,\rho_\xi\sigma)$ for each $\xi\in\Delta$, such that for every $X\in\mathrm{Hom}(\rho_\zeta,\rho_\xi\rho_\eta)$,
\begin{equation} X\cE_\sigma(\zeta)=\rho_\xi(\cE_\sigma(\eta))\,\cE_\sigma(\xi)\,\sigma(X)\,.
\end{equation}
For our systems this reduces to
\begin{align}
&&\cE_\sigma(g+h)=\alpha_g(\cE_\sigma(h))\,\cE_\sigma(g)\,,\label{HBa}\\
&&\cE_\sigma(\rho)=\alpha_g(\cE_\sigma(\rho))\,\cE_\sigma(g)\,,\label{HBb}\\
&&U_g\cE_\sigma(\rho)=\rho(\cE_\sigma(g))\,\cE_\sigma(\rho)\sigma(U_g)\,,\label{HBc}\\
&&S_g\cE_\sigma(g)=\rho(\cE_\sigma(\rho))\,\cE_\sigma(\rho)\,\sigma(S_g)\,,\label{HBd}\\
&&T_z\cE_\sigma(\rho)=\rho(\cE_\sigma(\rho))\,\cE_\sigma(\rho)\sigma(T_z)\,,\label{HBe}\end{align}
for all $g,h\in G$, $z\in \cF$. There may be more than 1 half-braiding associated to a given
$\sigma$; in that case we denote them by $\cE_\sigma^j$. As we will see shortly, 
knowing the half-braiding
is equivalent to knowing the \textit{matrix units} of the corresponding simple summand
of Tube$\,\Delta$ (the matrix units $e_{i,j}$ of a matrix algebra isomorphic to $M_{k\times k}$
are a basis satisfying $e_{i,j}e_{m,l}=\delta_{j,m}e_{i,l}$  --- e.g. the standard basis of
$M_{k\times k}$). If we decompose $\sigma=\sum_gk_g\alpha_g+k_\rho\rho$ into a sum of irreducibles,
then each half-braiding $\cE^j_\sigma$ will correspond to a distinct matrix subalgebra of
Tube$\,\Delta$ isomorphic to $M_{k\times k}$, where $k=\sum k_g+k_\rho$. 

The dual principal graph for the Longo-Rehren inclusion of $\Delta$ can be read off from
the collection of half-braidings as follows. On the bottom are the simple sectors of $\Delta$;
on the top row are the (inequivalent) half-braidings $\cE^j_\sigma$. If we write $\sigma=
\sum k_g\alpha_g+k_\rho\rho$, connect $\cE^j_\sigma$ to $\alpha_g$ with $k_g$
edges, and to $\rho$ with $k_\rho$ edges. This forgetful map $\cE^j_\sigma\mapsto\sigma$
is called \textit{alpha-induction} and plays a central role in much of the theory.
See Figure 4 for an example, which follows from the tube algebra analysis of Subsection 4.2.
(The principal graph for the Longo-Rehren inclusion is much simpler: just  $\Delta\times\Delta$ on
the bottom and $\Delta$ on the top, with edge multiplicities given by fusion multiplicities.)

\medskip\epsfysize=2.2in\centerline{ \epsffile{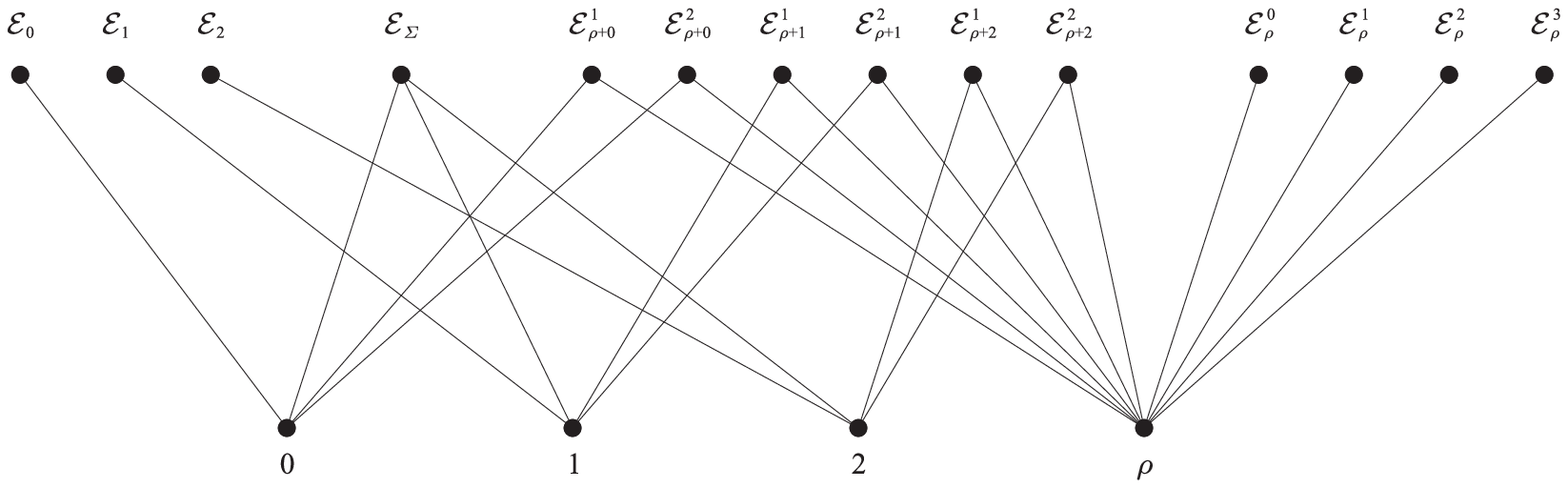}}\medskip

\centerline{Figure 4. The dual principal graph for the double of $\Delta(\bbZ_3+2)$}
\label{fig4}\medskip

The point is that the centre of the tube algebra is nondegenerately braided.
A nondegenerately braided system comes with \textit{modular data}:

\medskip\noindent\textbf{Definition 2.} \textit{Unitary matrices $S=(S_{a,b})_{a,b\in\Phi},T=
(T_{a,b})_{a,b\in\Phi}$ are called} modular data
\textit{if $S$ is symmetric, $T$ is diagonal and of finite order, the assignment $\left({0\atop 1}{-1\atop
0}\right)\mapsto S$, $\left({1\atop 0}{1\atop 1}\right)\mapsto T$ generates a representation of
SL$_2(\bbZ)$, and there is some index $0\in
\Phi$ such that $S_{a,0}^2>0$ $\forall a\in\Phi$ and the quantities 
\begin{equation} N_{a,b,c}:=\sum_{d\in\Phi}\frac{S_{a,d}S_{b,d}S_{c,d}}{S_{0,d}}\label{verl}\end{equation}
are nonnegative integers.}\medskip

The $a\in\Phi$ are called \textit{primaries}; $0\in\Phi$ is called the \textit{identity}; the
$N_{a,b,c}$ are called \textit{fusion coefficients} and \eqref{verl} is called \textit{Verlinde's
formula}. $S^2$ will necessarily be a permutation matrix, called \textit{charge-conjugation} $C$.
In the case of
modular data associated to nondegenerately braided fusion categories,  the primaries $\Phi$
label the simple objects, 
the quantities $N_{a,b}^c:=N_{a,b,Cc}$ are the structure constants of the fusion ring.
Nondegenerately braided systems of endomorphisms always satisfy the
stronger inequality $S_{a,0}>0$, so we will assume this for now. The 
\textit{quantum-dimension} $S_{a,0}/S_{0,0}$ in this case is the statistical dimension $d_a=\sqrt{[M:
a( M)]}$ 
and $1/S_{0,0}=\sqrt{\sum_ad_a^2}$ is the \textit{global dimension}. When $S_{a,0}=S_{0,0}$ then $a$ has an
inverse in $\Phi$ and is called a \textit{simple-current}.

In the case of the centre of Tube$\,\Delta$, the primaries are in one-to-one
correspondence with the simple summands of Tube$\,\Delta$, or equivalently with the
half-braidings defined above. These matrices can be computed once we know the 
\textit{matrix entries} $\cE_\sigma^j(\xi)_{(\eta,\alpha),(\eta',\alpha')}\in\mathrm{Hom}(\rho_\eta
\cdot\rho_\xi,\rho_\xi\cdot\rho_{\eta'})$ for each irreducible $\xi\in\Delta$, as $\eta,\eta'$
run over all simple objects (with multiplicities) in $\sigma$. In fact, the diagonal entries
$(\eta',\alpha')=(\eta,\alpha)$ suffice to determine $S,T$. 
These matrix entries can be computed from either the half-braidings, or from the
diagonal matrix units $e(\sigma^j)_{i,i}$, as follows. Let $W_\sigma(\eta,\alpha)$ be an orthonormal
basis of Hom$(\rho_\eta,\sigma)$; then we have
\begin{align}
&&\cE_\sigma^j(\xi)_{(\eta,\alpha),(\eta',\alpha')}=\rho_\xi(W_\sigma(\eta',\alpha')^*)\,\cE_\sigma^j(\xi)\,
W_\sigma(\eta,\alpha)\,,\\ &&e(\sigma^j)_{(\eta,\alpha),(\eta',\alpha')}=\frac{d_\sigma}{\lambda
\sqrt{d_\eta\,d_{\eta'}}}\sum_\xi d_\xi\,(\eta\xi|\cE_\sigma^j(\xi)_{(\eta,\alpha),(\eta',\alpha')}|
\xi\eta')\,,\end{align}
where 
\begin{equation}\lambda=n+\delta^2=2n+n'\delta\label{globdim}\end{equation}
 is the {global dimension}. 
 The entries of the diagonal unitary matrix $T$ and symmetric unitary matrix $S$ are determined from the matrix entries $\cE_\sigma^j(\xi)_{(
 \eta,\alpha),(\eta,\alpha)}$ through:
 \begin{align}
&& T_{\sigma^j,\sigma^j}=d_\xi\,\phi_\xi(\cE_\sigma^j(\xi)_{(\xi,\alpha),(\xi,\alpha)})\,,\label{Tdef}\\
&&S_{ \sigma^j,\sigma'{}^{j'}}=\frac{d_\sigma}{\lambda}\sum_{(\xi,\alpha)}d_\xi\,\phi_\xi(
\cE_{\sigma'}^{j'}(\eta)_{(\xi,\alpha),(\xi,\alpha)}^*\,\cE_\sigma^j(\xi)_{(\eta,\alpha'),(\eta,\alpha')}^*)
\,,\label{Sdef}\end{align}
where $\phi_\xi$ is the \textit{standard left inverse} of the endomorphism $\rho_\xi$, defined by
$\phi_\xi(x)=R_{\rho_\xi}^*\overline{\rho_\xi}(x)\,R_{\rho_\xi}$. In \eqref{Tdef},
$\xi$ can be any irreducible in $\sigma$, and in \eqref{Sdef} the sum is over all $\xi$ (counting
multiplicities) in $\sigma'$ while $\eta$ is any (fixed) irreducible in $\sigma$.

\subsection{The first class near-group C$^*$-categories: $n'=n-1$}

As before, we will sometimes write $a_x,b_x,b''_{x,y}$ for $a(x),b(x),b''(x,y)$.
We will first show that, as a C$^*$-algebra,
\begin{equation}
\mathrm{Tube}\,\Delta\cong\bbC^{2n+1}\oplus M_{n\times n}\oplus \left(M_{2\times 2}\right)^{n^2-n}\,,\label{tube1}
\end{equation}
unless $s=-1$ and $n=7$, in which case
\begin{equation}
\mathrm{Tube}\,\Delta\cong\bbC^{7}\oplus M_{7\times 7}\oplus \left(M_{2\times 2}\right)^{44}
\,.\label{tube2}\end{equation}

It will prove  useful to know the formulae
\begin{align}
&&\rho U_g=\sum_hS_h\,S_{g+h}^*+\sum_w\sigma^2w(g)\,T_w\,U_g\,T_w^*\,,\label{newa}\\
&&b''(x,y)=a_x\,a_y\,a_{xy}\,\overline{b(\sigma^2x\,\sigma y)}\,\overline{b''(\overline{y},\overline{x})}
\,,\label{newb}\\ &&b''(\sigma^2(-\sigma x),\sigma(-\sigma x))=b''(x,\sigma^2x)\,.\label{newc}
\end{align}
\eqref{newa} is implicit in Section 2.3. \eqref{newb} follows from the sequence
\eqref{bpplcbpp},\eqref{bpppibpp},\eqref{bpplcbpp}. \eqref{newc} is trivial when $n$ is odd
(where $-1=1$), so it suffices to take $s=a=b=1$; then apply the sequence
\eqref{bpplcbpp},\eqref{bpppibpp},\eqref{bpprcbpp} to the left-side.
Note also that we always have $a_{-x}=a_x$ and $b_{-x}=b_x$ (when both sides are
defined).

\medskip\noindent\textbf{Theorem 5.} \textit{Consider any type $G+n-1$ $C^*$-category, i.e.
let $s,a,b,b''$ be any solution of the equations of Theorem 3.} 

\smallskip\noindent\textbf{(a)} \textit{There is precisely one half-braiding 
$\cE_{\alpha_g}$ for any $g\in G$: $\cE_{\alpha_g}(h)=1$ and $\cE_{\alpha_g}(\rho)=(-1)^{n'g} U_g$. We get the diagonal  matrix entries $\cE_g(h)_{g,g}=1$, $\cE_g(\rho)_{g,g}=(-1)^{n'g}U_g$. 
These half-braidings correspond to central projections} 
\begin{equation}
\pi_g=n^{-1}(n+1)^{-1}\sum_hA_{g,h}+(n+1)^{-1}(-1)^{gn'}B_{g,g}\,.\end{equation}

\noindent\textbf{(b)} \textit{There is precisely one half-braiding 
$\cE_{\sum_g\alpha_g}$. It corresponds to the $n\times n$ matrix algebra $\cA(\sum\alpha)$, 
spanned by $B_{g,h}$ ($g\ne h$) and} 
\begin{equation}\pi'_{g}:=(n+1)^{-1}\sum_hA_{g,h}+(n+1)^{-1}(-1)^{gn'+1}B_{g,g}\qquad\forall g\,.
\end{equation}
\textit{The corresponding matrix entries are $\cE_{\Sigma}(g)_{h,h}=1$, $\cE_\Sigma(\rho)_{g,g}=n^{-1}(-1)^{gn'+1}U_g$.}
\smallskip

\noindent\textbf{(c)} \textit{There are precisely $n-1$ half-braidings 
$\cE^w_{\rho+\alpha_g}$ for any  $g\in G$, one for each $2\times 2$ matrix
algebra $\cA(g,w)=\mathrm{Span}\{p_{g,w},C_{g,w},D_{g,\sigma w},D_{g,\sigma w}C_{g,w}\}$.
The corresponding matrix entries are $\cE^w_{\rho+g}(h)_{g,g}=\overline{w(h)}$, $\cE^w_{\rho+g}
(h)_{\rho,\rho}=(-1)^{n'h}\overline{w(g+h)}\,U_h^*$, and}
\begin{align}&&\cE^w_{\rho+g}(\rho)_{\rho,\rho}=
\sum_h(-1)^{n'h}\overline{w(h)}\,S_hS_h^*\qquad\qquad\qquad\qquad\nonumber
\\&&+
n^{-2}\overline{a_w}\,\overline{b_w}\sum_xb''_{\sigma^2w\,\overline{\sigma^2(wx)},x}\,
b''_{\sigma w,x\sigma^2w\,\overline{\sigma^2(wx)}}\,\sigma^2(wx)(g)\,\sigma(w\,\overline{\sigma(wx)})(g)\,T_xT^*_{-\overline{w}x\,\overline{\sigma^2(wx)}}\,.\nonumber\end{align}

\noindent\textbf{(d)} \textit{When $s=\omega=1$ or $n\le 3$, there are precisely $n+1$ half-braidings 
$\cE^\psi_{\rho}$, naturally parametrised by the characters $\psi\in\widehat{\bbF_{n+1}^+}$.
The matrix entries are $\cE_\rho^\psi(g)_{\rho,\rho}=(-1)^{gn'}U_g^*$ and}
\begin{equation}
\cE^\psi_\rho(\rho)_{\rho,\rho}=\zeta_1\,\psi(1)\,\sum_k(-1)^{kn'}S_kS_k^*+
\sum_x\zeta_x\psi(\sigma x)\,T_xT^*_{\overline{\sigma^2x}}\,,\end{equation}
\textit{where $\zeta_1,\zeta_x$ is any particular solution to the equations \eqref{rhoproj}
(when $b''$ and $b$ are identically 1, take $\zeta$ identically 1).}\smallskip

\noindent\textbf{(e)} \textit{When $s=-1$ and $n>3$, then $n=7$ and  there are precisely 
2 half-braidings for $\sigma=2\rho$, with matrix entries $\cE_{2\rho}^{s_1}(g)_{(\rho,s_2),(\rho,s_2)}=U_g^*$ and
\begin{equation}
\cE_{2\rho}^{s_1}(\rho)_{(\rho,s_2),(\rho,s_2)}=\i s_1\sum_g S_gS_g^*+\i s_2 T_1T_5^*+
s_1s_2\i T_6T_4^*\,,\end{equation}
}
where $s_1,s_2\in\{\pm 1\}$.\smallskip

\noindent\textbf{(f)} \textit{There are no $C^*$-categories of type $G+n-1$ with $\omega\ne 1$ and $n>2$.}

\medskip\noindent\textit{Proof.} From \eqref{mult} we obtain the formulae \begin{align}
&&A_{g,h}\,A_{k,l}=\delta_{g,k}A_{g,h+l}\,,\ A_{g,h}\,B_{k,l}=\delta_{g,k}
B_{k,l}\,\ 
B_{g,h}B_{k,l}=\delta_{h,k}\delta_{g,l}\sum_mA_{g,m}+\delta_{h,k}c_{ghl}B_{g,l}\,,\nonumber\\
&&A_{g,h}C_{k,z}=\delta_{g,k}z(h)\,C_{k,z}\,,\ B_{g,h}\,C_{k,z}=0\,,\ 
C_{g,w}D_{k,z}=
\delta_{g,k}\delta_{\sigma w,z}a_zb_z\sum_h\overline{w(h)}A_{g,h}\,,\nonumber\\
&&D_{g,z}C_{k,w}=\delta_{g,k}\delta_{z,\sigma w}n^{-1}a_wb_w\sum_h
(-1)^{n'h}\overline{w(h+g)}E_h+\delta_{g,k}\delta_{z,\sigma w}{b^{-2}_{w}}\sum_h(-1)^{n'h}\overline{w(h)}\,E'_h
\nonumber\\&&+\delta_{g,k}\delta_{z,\sigma w}\sum_xb''_{\sigma^2w\,\overline{\sigma^2(wx)},x}\,
b''_{\sigma w,x\sigma^2w\,\overline{\sigma^2(wx)}}\,\sigma^2(wx)(g)\,\sigma(w\,\overline{\sigma(wx)})(g)\,E''_{x,-\overline{w}x\,\overline{\sigma^2(wx)}}\,,\nonumber\\
&&C_{g,z}E_k=z(k)(-1)^{n'k}C_{g,z}\,,\ C_{g,z}E'_k=s(-1)^{n'k}z(k-g)C_{g,z}\,,\nonumber\\&&
C_{g,z}E''_{w,x}=C_{g,z}\delta_{x,\sigma (z\overline{w})}x(g)\sigma^2w(g)b''(x\sigma^2w,z)
\overline{b''(z,x\sigma^2w)}\,,\nonumber
\end{align}
where $c_{ghl}:=\sum_zz(g)\,\sigma z(h)\,\sigma^2z(l)$. Using $(-1)^{n'h}=z(h)\sigma z(h)
\sigma^2z(h)$, we get $c_{ghh}=c_{hgh}=c_{hhg}=(-1)^{n'h}(n\delta_{g,h}-1)$.
Write $e=\sum_h (-1)^{n'h}E_h$, $e'=\sum_h(-1)^{n'h}E'_h$, $e''_x=E''_{x,\overline{\sigma^2x}}$.
From \eqref{star} we obtain 
\begin{align}&&A_{g,h}^*=A_{g,-h}\,,\ B_{g,h}^*=B_{h,g}\,,\ 
C_{g,z}^*=sz(g)a_z\overline{b_z}D_{g,\sigma z}\,,\nonumber\\
&&e^*=e\,,\ e'{}^*=s\delta_{-1,1}e'+se''_2\,,\ e''_z{}^*=\delta_{z,2}e'+a_z\overline{b_z}\overline{b''(\overline{\sigma z},\sigma^2(-\overline{\sigma z}))}
\,e''_{\sigma^2(-\sigma z)}\,.
\nonumber
 \end{align}
We use in these expressions that $\sigma(2)=-1$, where 2 denotes the element $1+1$ in the
corresponding field $\bbF_{n+1}$; when the characteristic of $\bbF_{n+1}$ is 2, then neither
$-1=1$ nor $2=0$ lie in $\cF$, so the corresponding terms should be ignored. 

Let's begin by solving the half-braiding equations for (a). For later convenience, we'll solve
this for any near-group C$^*$-category.
 Since Hom$(\alpha_{g+h},\alpha_{g+h})$ and Hom$(\alpha_g\rho,
\rho\alpha_g)$ are both one-dimensional, we need to find numbers $c_{g,h},c'_g\in\bbC$ for each
$g,h\in G$, such that $\cE_{\alpha_g}(h)=c_{g,h}$ and $\cE_{\alpha_g}(\rho)=c'_g U_g$.
\eqref{HBb} and \eqref{HBe} give $c_{g,h}$ and $c'_g$, respectively. The remaining three 
equations are then automatically  satisfied. \eqref{HBe} becomes
\begin{align}
&&-s\delta^{-3/2}\sum_h\langle g+k,h\rangle\,\overline{\ddot{(-h\overline{y})}(h)}\,
\overline{a_{\overline{y}}}\,\delta_{-h\overline{y},z'}=\sum_{z,x}\ddot{z}(g)\,\overline{x(k)}\,
b_{z,x}\,\overline{b''_{gz,z',y,x}}\,,\label{sig1a}\\
&&\overline{\ddot{z}(g)}\,c_g^{\prime-1}\,\overline{z(g)}\,\langle g,h-g\rangle\delta_{x,gz}=
W_{xhz}(g)\,,\label{sig1b}\\
&&c_g^{\prime-1}\,\overline{\langle g,g\rangle}\overline{\ddot{z}(g)}\,\overline{z(g)}\,\ddot{x}(g)\,
\delta_{w,gz}\,\delta_{y,gx}=V_{wxyz}(g)
\,,\label{sig1c}\end{align}
where $W,V$ are defined in the proof of Proposition 1.
It is easy to see from Theorems 3 and 4 that \eqref{sig1a} is automatically satisfied for both
classes.

When $n'=n-1$ this condition becomes $c_g'=\overline{{z}(g)}\,\overline{\sigma{
z}(g)}\,\overline{\sigma^2 z(g)}$ is independent of $z$. Using the finite field expression 
$\sigma x=(1-x)^{-1}$ of Proposition 2, ${z}\,\sigma{
z}\,{\sigma^2 z}$ collapses to $-1\in\bbF_q$, so this expression for $c'_g$ simplifies to $c_g'=
(-1)^g$ if $n$ is even, or 
$c'_g=1$ if $n$ is odd. In this case, \eqref{HBd} is automatically satisfied.

Now turn to the proof of (b).
We see from \eqref{mult} or the products collected at the beginning of this proof, 
that $\cA(\sum\alpha)$ is an $n^2$-dimensional  C$^*$-algebra,
so is a direct sum of matrix algebras. The $\pi'_{g}$ obey $\pi'_g{}^*=\pi'_g$, $\pi'_g\pi'_h=\delta_{gh}
\pi'_g$, and $\pi'_gB_{h,k}=\delta_{g,h}B_{h,k}$, so the projections $\pi'_g$ are the diagonal matrix
units. For any $g$, there is an $n$-dimensional space in $\cA(\sum\alpha)$ satisfying
$\pi'_gx=x$, namely the span of $\pi'_g,B_{g,h}$ ($h\ne g$), so $\pi'_g$ must belong to an
$n\times n$ block and hence $\cA(\sum\alpha)$ is $M_{n\times n}(\bbC)$.

To show the $\pi'_g$ are minimal, i.e. $\cA(\sum\alpha)$ is maximal as a simple C$^*$-subalgebra of Tube$\,\Delta$, it suffices
to show that $\pi'_g\pi_k=B_{g,h}\pi_k=\pi'_gC_{k,z}=B_{g,h}C_{k,z}=0$ for all $k,z$ and all 
$g\ne h$. This is clear from the products listed earlier.

The C$^*$-algebras $\cA_{g,g}$ are each isomorphic to $\bbC^{n+1}$, with projections
$\pi_z$ of Theorem 5(a), $\pi'_g$ of Theorem 5(b), and for all $z\in\widehat{G}^*$,
\begin{equation}
p_{g,z}=n^{-1}\sum_h\overline{z(h)}\,A_{g,h}\,.\end{equation}
Together with the $n\times n$ matrix algebra of Theorem 5(b) and the projections $\pi_g$
of Theorem 5(a), these projections $p_{g,z}$ span all of $\sum_{g,h}\cA_{g,h}$. Note that
$p_{g,z}\pi_k=p_{g,z}\pi'_k=0$ and $p_{g,z}^*=p_{g,z}$.


Now turn to (c). Each $\cA(g,w)$ is clearly a 4-dimensional C$^*$-algebra, and is noncommutative
since $C_{g,w}D_{g,\sigma w}\in\sum_h\cA_{g,h}$ and $D_{g,\sigma w}C_{g,w}\in\cA_{\rho,\rho}$
are distinct. Therefore each $\cA(g,w)$ is isomorphic as a C$^*$-algebra to the $2\times 2$
matrix algebra. Each $\cA(g,w)$ is readily seen to be orthogonal to $\pi_g$ and $\cA(\sum\alpha)$.
Note that $C_{g,w}D_{g,\sigma w}$ is a scalar multiple of $p_{g,w}$.
A basis for $\sum_g\cA_{g,\rho}+\sum_g\cA_{\rho,g}+\cA_{\rho,\rho}$ is $\cup_{g,w}\{C_{g,w},D_{g,\sigma w},
C_{g,w}D_{g,\sigma w},D_{g,\sigma w}C_{g,w}\}\cup\{e,e',e''_x\}_{x\in\widehat{G}^*}$.
To verify that $\cA(g,w)$ is maximal as a simple C$^*$-subalgebra of Tube$\,\Delta$, 
it suffices to verify that $C_{g,w}D_{g',w'}=0=D_{g',w'}C_{g,w}$ unless $g=g'$ and $\sigma w=w'$, $C_{g,w}e=C_{g,w}e'=C_{g,w}e''_x=0$ (the other orthogonalities are either trivial or follow from these).
This is elementary.

The diagonal matrix units are $p_{g,w}$ and $n^{-1}\overline{a_w}\,\overline{b_w}D_{g,\sigma w}
C_{g,w}$. From this we obtain the desired quantities.

Now, turn to (d).
We compute\begin{align}
&&e^2=ne\,,\ ee'=ne'=e'e\,,\ ee''_z=ne''_z=e''_ze\,,\ e'e'=n^{-1}s\delta_{-1,1}e+a_2b_2e''_{1/2}\,,
\nonumber\\
&&e'e''_z=n^{-1}\delta_{z,2}e+\overline{b_z}b_{-\sigma z}b''(\overline{\sigma(-{
\sigma z})},\overline{\sigma^2(-\sigma z)})\,e''_{\sigma(-\overline{\sigma z})}\,,\nonumber\\
&&e''_ze'=n^{-1}\delta_{z,2}e+b''(\overline{\sigma^2z},\overline{z})\,e''_{\sigma(-\overline{\sigma z})}\,,\nonumber\\
&&e''_ze''_w=n^{-1}sb_wa_z\overline{b''_{z,\sigma^2z}}\delta_{\sigma w,-\sigma z}e
+a_w\delta_{z,\overline{w}}\overline{b_z}e'+b''_{x,x'}\overline{b''_{z,y}}b''(\overline{\sigma^2z},
y)e''_{\sigma(-\overline{\sigma w\,\sigma z}\,\sigma({wz}))}\,,\nonumber\end{align}
where in the $e_z''e_w''$ equation $x,x',y$ are defined by $\sigma x=-\overline{\sigma w}\,
\sigma z$, $\sigma x'=x\,\sigma z$, and $\sigma y=\overline{\sigma^2z}\,\overline{\sigma w}$. 
Manifestly, the span of $e,e'$ and all $e''_x$ is an $n+1$-dimensional C$^*$-algebra, which we'll
call $\cA(\rho)$. It is immediate now that it is orthogonal to $\pi_g,\cA(\sum\alpha)$ and
all $\cA(g,w)$. Let's identify when it is abelian. First note that
\begin{align}
&&\overline{b_x}\,b_{-\sigma x}\,b''(\overline{\sigma(-\sigma x)},\overline{\sigma^2(-\sigma x)})
=b''(\overline{\sigma^2x},\overline{x})\qquad\mathrm{for}\ x\ne 2\,,\label{abela}\\ 
&&s\,b_w\,a_z\,\overline{b''_{z,\sigma^2z}}=s\,b_z\,a_w\,\overline{b''_{w,\sigma^2w}}\qquad
\mathrm{for}\ \sigma z=-\sigma w\,,\label{abelb}\\
&&a_w\,\overline{b_z}=a_z\,\overline{b_w}\qquad\mathrm{for}\ z=\overline{w}\,.\label{abelc}
\end{align}
Indeed, \eqref{abela} follows by applying \eqref{newb} twice and \eqref{newc} once;
any $b_y$ appearing in these expressions can be replaced with $sa_y$. To see \eqref{abelb},
use \eqref{newc} and the substitution $-\sigma z=\sigma w$; again, $b=sa$ here.
\eqref{abelc} is trivial.

Thus to conclude the argument that $\cA(\rho)$ is commutative, it suffices to verify that
$b''(x,x')\,\overline{b''(z,y')}\,b''(\overline{\sigma^2z},y')$ is invariant under the switch $z
\leftrightarrow w$. 
Using \eqref{bpplcbpp},\eqref{bppbppbpp},\eqref{bpprcbpp} and \eqref{newb}, we obtain
\begin{align}
&&b''_{x,x'}=\overline{b''(w,\sigma z\,\sigma^2w)}\,a_x\,b_x^2\,b_{yz}\,a_{\sigma^2w\,\sigma z}\,
\overline{b_{\sigma^2w\,\sigma z}}^2\,a(\sigma z\,\overline{\sigma w})\,a\overline{b}(\sigma^2w\,\sigma(\sigma^2w\,\sigma z))\,,\\
&&\overline{b''_{z,y}}=\overline{b''(\sigma^2z,\sigma w)}\,a\overline{b}(z)\,\overline{b_w}\,
\overline{b_{yz}}\,b(\sigma^2w\,\sigma z)\,b(\sigma^2z\,\sigma w)\,,\\
&&b''(\overline{\sigma^2z},y)=\overline{b''(z,\sigma^2z\,\sigma w)}\,b_z\,b_{x'}\,b_x\,
\overline{b(\sigma^2z\,\sigma w)}\,,\end{align}
where we write $a\overline{b}(z)$ for $a(z)\,\overline{b(z)}$ etc. We likewise have
\begin{equation}
\overline{b''(\sigma^2w,\sigma z)}=\overline{b''(\sigma^2z,\sigma w)}\,s\,ab(z)\,a\overline{b}(w)\,
a\overline{b}(\sigma^2z\,\sigma w)\,\overline{b(\sigma^2w\,\sigma z)}^2\,.\end{equation}
Since $x,x',y$ depend on $w,z$, they will be affected by the switch $w\leftrightarrow z$:
in particular we find $x$ becomes $\overline{\sigma^2 x}$, $y$ becomes $\overline{\sigma^2w
\,\sigma z}$, and $x'$ is unchanged. Thus $\cA(\rho)$ is commutative iff
\begin{align}
&&a\overline{b}(\sigma^2w\,\sigma z)\,a\overline{b}(\sigma^2w\,\sigma(\sigma^2w\,\sigma z))=
\qquad\qquad\qquad\nonumber\\&&a
\overline{b}(\sigma^2z\,\sigma w)\,a\overline{b}(\sigma^2z\,\sigma(\sigma^2z\,\sigma w))\,
a\overline{b}(\sigma^2z\,\sigma w)\,\overline{b(\sigma^2w\,\sigma z)}^2\,.\end{align}

Consider first the generic case, where $\omega=1$. We find that, provided
$\sigma w\ne -\sigma z$ and $z\ne \overline{w}$, $e''_ze''_w=s\,e''_we''_z$. This means
$\cA(\rho)$ will be commutative provided $s=\omega=1$ or $n\le 3$.

Much more subtle is when $\omega$ is a primitive third root of unity. In this case $-1=\sigma(-1)$ and $s=1$;
$a_x=b_x=1$ except for $b_{-1}=\omega$. Then $\cA(\rho)$ is commutative iff
\begin{equation}
b(\sigma^2w\,\sigma z)\,b(\sigma^2w\,\sigma(\sigma^2w\,\sigma z))=
b(\sigma^2z\,\sigma w)\,b(\sigma^2z\,\sigma(\sigma^2z\,\sigma z))\label{commut1}\end{equation}
for all $w,z$ with $\sigma w\ne -\sigma z$ and $z\ne \overline{w}$. The only possible way
this equation can be violated is if at least one of those $b$'s doesn't equal 1.

Suppose first $\sigma^2z\,\sigma(\sigma^2z\,\sigma w)=-1$. Then hitting both sides with
$\sigma$, we get $-\sigma z=\overline{x}$. This implies $\sigma x'=-1$, i.e. $x'=-1$, and
hence $\sigma w=-\sigma^2x$. Therefore 
\begin{equation}-1=\sigma(-1)=\sigma(\sigma w\,\overline{\sigma^2x})
=\sigma^2w\,\sigma(\sigma x\,\overline{w})=\sigma^2w\,\sigma(-\sigma z\,\overline{\sigma w}
\,\overline{w})\end{equation}
and thus $\sigma^2w\,\sigma(\sigma^2w\,\sigma z)=-1$. Thus in \eqref{commut1},
$b(\sigma^2w\,\sigma(\sigma^2w\,\sigma z))$ and $b(\sigma^2z\,\sigma(\sigma^2z\,\sigma z))$
are always equal.

Finally, suppose $\sigma^2z\,\sigma w=-1$. Then $\sigma^2w\,\sigma z=\overline{x}$ cannot 
equal $-1$. This means that for any pair $x,z$ with $\sigma^2z\,\sigma w=-1$ and
$\sigma z\ne -\sigma x$ and $z\ne \overline{x}$, $e''_ze''_x=\omega e''_xe''_z$.

Consider now the case where $\cA(\rho)$  is commutative, i.e. where 
either $s=\omega=1$,  or $n<7$.
Then the $n+1$ minimal central projections in $\cA_{\rho,\rho}$ are scalar multiples of 
\begin{equation}
\pi(\vec{\zeta})=(n^2+n)^{-1}e+(n+1)^{-1}\zeta_1 e'+(n+1)^{-1}\sum_x\zeta_xe''_x
\label{minprojrho}\end{equation}
for some $\zeta_1,\zeta_x\in\bbC^\times$. These must satisfy
$e'\pi(\vec{\zeta})=\beta_{e'}\pi(\vec{\zeta})$ and $e''_x\pi(\vec{\zeta})=\beta_x\pi(\vec{\zeta})$ for scalars
$\beta_{e'},\beta_x\in\bbC$. This yields the equations
\begin{align}
&&\zeta_2=\overline{\zeta_1}\,,\ \zeta_{\overline{x}}=sa_xb_x\zeta_1\,\overline{\zeta_x}\,,\ 
\zeta_{\sigma(-\overline{\sigma x})}=\overline{b_x}b_{-\sigma x}\zeta_1\,\zeta_x\,b''_{x,\sigma^2x}\ \qquad(x\ne 2)\,,
\nonumber\\
&&\zeta_{\sigma(-\overline{\sigma z\,\sigma x}\,\sigma({xz}))}=\zeta_x\,\zeta_z\,
{b''_{x,x'}}\,\overline{b''_{z,y}}\,{b''_{\overline{\sigma^2z},y}}\qquad 
(z\ne \overline{x},\sigma^2(-{\sigma x}))\,,\label{rhoproj}\end{align}
as well as $\zeta_1^2=s$ when $n+1$ is even. 

First, note we can solve these equations in the special case that $b''$ and $b$ are
identically 1. In this case, identify $G$ with $\bbF_{n+1}^\times$ and take $\zeta_1=\psi(1)$ and 
$\zeta_x=\psi(\sigma x)$ for any of the $n+1$ characters of the additive group
$\bbF_{n+1}^+\cong(\bbZ_p)^k$: 
a little effort shows  these $n+1$ choices of $\zeta$'s all work, and so must exhaust all
solutions.

Recall that the C$^*$-categories with $n'=n-1$
 form a group: $\cC(b_1,b''_2)*\cC(b_2,b''_2)=\cC(b_1b_2,b''_1b''_2)$. Let $\pi(\vec{\zeta}(i))$ be solutions of \eqref{rhoproj} for
$b''_i$ respectively; then $\zeta_x=\zeta(1)_x\zeta(2)_x$ will be a solution for $\cC(b_1b_2,b''_1b''_2)$.
This implies that if you have any particular solution $\vec{\zeta}$ for a given $\cC(b,b'')$,
 all other solutions for that category $\cC(b,b'')$ are obtained by multiplying that particular solution by the solutions $\psi(\sigma x)$ for $\cC(1,1)$. 
 
 Since the sum of the minimal projections of $\cA(\rho)$ must
equal the unit $n^{-1}e$, we know now that the  $\pi(\vec{\zeta})$ given in
\eqref{minprojrho} are indeed the minimal projections (i.e. the coefficient $\lambda^{-1}$ for $e$ 
is correct).

When $s=-1$, both $e,e'$ are central elements. Suppose the centre
is not 2-dimensional. Then there is some $Z=\sum_x c_xe''_x\ne 0$ which commutes with all
$e''_w$. Suppose $c_x\ne 0$; then as long as $n>3$ it will be possible to choose a $z$
not equal to $0,1,x,\overline{x}$. Then $e''_zZ$ and $Ze''_z$ will differ by a sign in at least one
 coefficient (namely that of $e''_{\sigma(\sigma^2x\,\sigma^2z)\,\sigma(\overline{x}\overline{z})}$),
 which will be nonzero because $c_x$ is. This contradiction means that $\cA(\rho)$ will indeed
 be a sum of
 two matrix subalgebras. Note that $\beta e+\gamma e'=(\beta e+\gamma e')^2=\beta^2ne+2\beta\gamma ne'-\gamma^2n^{-1}
 e$ forces $\beta=(2n)^{-1}$ and $\gamma=\pm\i/2$. Thus the identities in the two matrix subalgebras are $1_\pm:=(2n)^{-1}e\pm \i e'/2$, and the two matrix
 subalgebras will be the images $e_\pm(\cA(\rho))=:\cA(\rho)\pm$, and are of equal dimension.
 
 Let's finish off the proof of (e). The C$^*$-algebra $\cA(\rho)_+$ is spanned by
 $1_+$ and $z_+:=(\gamma_z\,e''_z-\i\overline{\gamma_z}e''_{\overline{z}})/2$ where
 $\gamma_z=\sqrt{b''(\sigma z,z)}$ (since $b''(\sigma z,z)=-\overline{b''(\overline{\sigma^2 z},
 \overline{z})}$, we can choose these square-roots so that $\gamma_{\overline{z}}=\i\overline{
 \gamma_z}$, so $z_+=\overline{z}_+$). We compute $z_+^2=1_+$, so each $z_+$ is
 invertible. Provided $w\not\in\{z,\overline{z}\}$, we have
 \begin{equation}
 z_+w_+=\beta_{z,w}\,\left(\frac{wz+1}{w+z}\right)_+\end{equation}
 for some $\beta_{z,w}\in\bbC$, where we use addition in the finite field $\bbF_{n+1}$ 
 to simplify notation. But provided $w\not\in\{z,\overline{z}\}$, we know
 $e_z''e_w''=-e_w''e_z''$ and hence $z_+w_+=-w_+z_+$. Now choose some
 $x\not\in\{z,\overline{z},w,\overline{w},(wz+1)/(w+z),(w+z)/(wz+1)\}$ --- this is possible as
 long as $n>7$. Then
 \begin{equation}
 -(z_+w_+)x_+=x_+(z_+w_+)=-z_+x_+w_+=z_+w_+x_+\,,\end{equation}
 so $z_+x_+w_+=0$, which contradicts invertibility of $z_+,w_+,x_+$. This contradiction
 shows that $s=-1$ requires $n\le 7$. 
 
 When $\omega\ne 1$, say $\omega=e^{\pm 2\pi\i/3}$, we see that every $z\ne -1$ will have 
 precisely one $x$ such that $\sigma^2z\,\sigma x=-1$. This means $e_z''e_x''=e''_xe''_z$
 unless $x\in\{-1/(1+z),-(1+1/z)\}$, using finite field notation for addition as usual. Note that $e''_{-1/(1+z)}
 \in\bbC e'e''_z$ and $e''_{-(1+1/z)}\in\bbC e''_{-1}e_z$.
 The centre of $\cA(\rho)$ manifestly contains $e,e',e''_{-1}$, and
exactly as in the $s=-1$ argument, we see that the centre cannot be more than 3-dimensional.
 Therefore $\cA(\rho)$ here must be a sum of 3 matrix subalgebras. We compute the corresponding
 identities as before, obtaining $1_\zeta=(3n)^{-1}e+e^{\mp 2\pi\i/9}\zeta e'/3+e^{\mp 2\pi\i/9}
 \zeta^2e''_{-1}/3$ for each third root $\zeta$ of unity. Thus $\cA(\rho)$ will be a sum of the
 matrix subalgebras $\cA(\rho)_\zeta:=1_\zeta\cA(\rho)=\mathrm{span}\{1_\zeta,1_\zeta e''_z\}$.
 Now, $1_\zeta e''_x\in \bbC 1_\zeta e''_z$ iff $x\in\{z,-1/(1+z),-(1+1/z)\}$, so we see that
 $1_\zeta e''_z$ and $1_\zeta e''_x$ will always commute, and thus each subalgebra
 $\cA(\rho)_\zeta$ will be commutative. Since they are also simple, each must be $\bbC$,
 and we have that $n+1=3$. This concludes the proof of (f).
 \qquad
 \textit{QED to Theorem 5}\medskip

We already knew (from Proposition 5) that $s=1$ for $n\ne 1,3,7$ and $\omega=1$ for $n\ne 2$,
but we wanted to derive it directly to demonstrate the effectiveness of the structure of
the tube algebra in constraining the sets of solutions. 
Let us now give particular solutions for $\zeta$ in the known cases where $a,b,b''$ are not all identically 
1.  When $n=1$ and $s=1$, there are no parameters $\zeta$. When $n=2$ and 
$\omega=e^{\pm 2\pi\i/3}$, take $\zeta_1=e^{\mp 2\pi\i/9}=\overline{\zeta_{-1}}$.
When $n=3$ and $s=-1$, take $\zeta_1=\zeta_{\omega^2}=\i$ and $\zeta_\omega=1$, where
we identify $\bbF_4$ with $\bbZ_2[\omega]$ for $\omega^3=1$ and identify $G$ with $\bbF_4^\times$. 

\medskip
\noindent\textbf{Corollary 6.} \textit{Fix any finite field $\bbF=\bbF_q$ for $q=n+1$, and identify
$G$ with $\bbF^\times$. Here is the modular data of the double of any system 
covered by Theorem 3. The global dimension is $\lambda=n^2+n$. The primaries come in 4 families, parametrised as follows:} 

\begin{itemize} 

\item $g\in G$;

\item \textit{the symbol} $\Sigma$;

\item  $w+h$ \textit{for
$w\in\widehat{G}^*$ and} $h\in G$;

\item \textit{either $\rho^\psi$ for $\psi\in\widehat{\bbF_{n+1}}$ (when $s=\omega=1$ or $n\le 3$), or
$\rho^{s_1}$ for $s_1\in\{\pm \}$ (when  $s=-1$ and $n=7$).}

\end{itemize}

\noindent \textit{Then the $T$ and $S$ matrices are given in block form by}
\begin{align}
&&T=\mathrm{diag}\left(1;1;\overline{w(h)};\overline{\zeta_1\psi(1)}\right)\,,\qquad\qquad\qquad\qquad\qquad\\
&&S=\frac{1}{n+1}{\scriptsize \left(\begin{matrix}n^{-1}&1&w(g)(n+1)n^{-1}&{1}\\
{1}&n&0&-1\\ w(g)(n+1)n^{-1}&0&(n+1)n^{-1}w'(h)w(h')&0\\
1&-1&0&\sum_{x\in\bbF^\times}\overline{\zeta_x}^2\overline{\psi(
\sigma x)}\overline{\psi'(({\sigma x})^{-1})}\end{matrix}\right)}\,.\end{align}
\smallskip
\textit{except for $n=7$ when $s=-1$, when}
\begin{align}
&&T=\mathrm{diag}\left(1;1;\overline{w(h)};\pm\i\right)\,,\qquad\qquad\qquad\\
&&S=\frac{1}{8}{\left(\begin{matrix}7^{-1}&1&\frac{8}{7}w(g)&2\\
{1}&7&0&-2\\ \frac{8}{7}w(g)&0&\frac{8}{7}w'(h)w(h')&0\\
2&-2&0&-4s_1s_1'\end{matrix}\right)}\,.\end{align}
\smallskip

The primary labelled $g$ corresponds to the half-braiding with $\sigma={\alpha_g}$ in Theorem 5(a); $\Sigma$ 
corresponds to $\sigma=\sum\alpha_g$ in 5(b);  $w+g$ corresponds to the half-braiding
$\cE_{\rho+g}^w$ of 5(c); $\rho_j$ or $\rho,\rho'$ or $\rho,\rho',\rho''$ correspond to 
the half-braidings with $\sigma=\rho$ in 5(d). The proof of Corollary 6 is an elementary
calculation based on the matrix entries listed in Theorem 5, as well as the formulae of
Section 4.1. The most difficult is the bottom-right block in the $S$ matrix. Consider
$n+1$ odd (so $s=1$). Then
\begin{align}
&&S_{\rho^\psi,\rho^{\psi'}}=\frac{n^2}{\lambda}S_0^*\left(\overline{\zeta_1\psi(1)}\sum_g
(-1)^g\rho(S_g)\rho(S_g)^*+\sum_z\overline{\zeta_z\psi(\sigma z)}\rho(T_{\overline{\sigma^2z}})
\rho(T_z)^*\right)\nonumber\\
&&\times\left(\overline{\zeta_1\psi'(1)}\sum_h
(-1)^h\rho(S_h)\rho(S_h)^*+\sum_w\overline{\zeta_w\psi'(\sigma w)}\rho(T_{\overline{\sigma^2w}})S_0
\rho(T_w)^*\right)\nonumber\\
 &&=n^{-1}\left(\overline{\zeta_1\psi(1)}T_2^*T^*_{1/2}+\sum_z\overline{\zeta_z\phi(\sigma z)}
 b_z\overline{b''(\sigma^2(-\overline{\sigma z}),\sigma(-\overline{\sigma z}))}T^*_{\sigma(-\overline{\sigma z})}
 T^*_{\sigma^2(-\overline{\sigma z})}\right)\nonumber\\
&&\times\left(\overline{\zeta_1\psi'(1)}T_{1/2}T_{2}+\sum_w\overline{\zeta_w\phi'(\sigma w)}
 \overline{b_w}\overline{b''(\sigma^2(-{\sigma w}),\sigma(-{\sigma w}))}T_{\sigma^2(-{\sigma w})}
 T_{\sigma(-{\sigma w})}\right)\,,\nonumber\end{align}
which simplifies down to the given expression.
 
 Note that when $s=\omega=1$, this recovers the modular data for the double of the
 affine group Aff$_1(\bbF_{n+1})$ (see e.g. \cite{DV3,CGR} for the general theory of finite
 group modular data and its twists by cocycles in $H^3(G;\bbT)$). Recall that the primaries of the (untwisted) double of a 
 finite group are pairs $(g,\psi)$ of a conjugacy class representative, and an irrep of 
 the centraliser of $g$ in the full group. The primaries of the double of $\cC(1,1)$ and the 
 double of Aff$_1(\bbF_{n+1})\cong \bbF_q\sdprod\bbF^\times_q$ match up quite nicely as follows:
 what we call $g\in G$ corresponds to the pair $(e,\psi)$ where $e$ is the identity and
 $\psi$ is a 1-dimensional representation of Aff$_1(\bbF_{n+1})$;
 $\Sigma$ corresponds to $(e,\rho)$ where $\rho$ is the $n$-dimensional irrep;
 $w+h$ corresponds to conjugacy class $(w,0)$ and irrep $g\in\widehat{\bbF^\times_q}$;
 $\rho^\psi$ corresponds to conjugacy class (1,1) and irrep $\psi$ of the centraliser
 $\bbF^+_q$.

Note that in each case the modular data is inequivalent for the two systems at
$n=1$, the three at $n=2$, the two at $n=3$ and at $n=7$. This then verifies that for $n=1,2,3,7$,
the solution $b=b''=1$ corresponds to Mod(Aff$_1(\bbF_q)$) (for the other $n$, this is clear
by the uniqueness in Proposition 5. This inequivalence of the modular data also means that those systems
are not Morita equivalent, i.e. there cannot exist a subfactor $N\subset M$ for which the
principal even sectors form say the $s=+1$ system at $n=7$ and the dual principal sectors form 
the $s=-1$ system at $n=7$.

 The modular data for the 3 systems with $G=\bbZ_2$ was also computed in
 Section 4 of \cite{iz3}, where it was remarked that the modular data corresponds to that of
 the double of $S_3\cong\mathrm{Aff}_1(\bbZ_2)$ and its twists by order-3 cocycles in 
 $H^3(S_3;\bbT)\cong\bbZ_6$. The order-3 twist arises because all that the twist is allowed to
affect is the primary corresponding to conjugacy class (1,1) and (projective) irreps $\psi$ of
its centraliser $\bbZ_3$. In other words the cocycle must be nontrivial on $\bbZ_3<S_3$,
and coboundary on the $\bbZ_2$ subgroups.

It appears that $s=-1$ for $n=1,3,7$ likewise corresponds to twisted modular data.
This is clear for $n=1$: $H^3(\bbZ_2;\bbT)\cong\bbZ_2$ and $s=\pm 1$
corresponds to $\pm 1$-twisted data for $\bbZ_2\cong\mathrm{Aff}_1(\bbF_2)$.
$\mathrm{Aff}_1(\bbF_4)$ is isomorphic to the alternating group $A_4$, and the natural
restriction of $H^3(A_4;\bbT)$ to the subgroup $\bbZ_2^2<A_4$ is $\bbZ_2^3$. The $s=-1$
modular data appears to agree with the twist by some cocycle in that $\bbZ_2^3$, although
we haven't yet fleshed out the details.
For $n=7$, something analogous will hold, except now we want a twist which, when
restricted to $\bbZ_2^3<\bbZ_2^3\sdprod\bbZ_7$, is not `CT' in the sense of \cite{CGR} (because
the number of primaries for $s=-1$, $n=7$ is less than for $s=1$, $n=7$).
Such cocycles do indeed exist here.

\subsection{The tube algebra in the second class}

Consider now the near-group systems of type $G+n'$ where $n'\in n\bbZ$. It is certainly 
expected that only finitely many $n'$ will work for a given $n$; in fact to our knowledge
all known near-group categories for abelian $G$ have $n'\in\{0,n-1,n\}$. Last subsection, we used the existence and properties
of the tube algebra for $n'=n-1$ to prove that $s=\omega=1$ except for $n=1,2,3,7$. Likewise,
we expect that the existence of the tube algebra should constrain the possible values $n'\in
n\bbZ$. 

As a preliminary step towards working out the tube algebra structure here for $n'>n$, consider the
half-braidings for $\sigma=\alpha_g$, when $n$ divides $n'$. We find (following the
proof of Theorem 5(a)) that 
$\cE_{\alpha_g}$ exists here iff, whenever $\widetilde{z}=
\widetilde{w}=1$,  we have $\dot{z}(g)=\dot{w}(g)$. Indeed,
when $n|n'$, equations \eqref{sig1b} and \eqref{sig1c} are both satisfied iff 
\begin{equation} c'_g=\epsilon_{\langle,\rangle}(g)\,\overline{\dot{z}(g)}\,\overline{\widetilde{z}(g)}
\end{equation}
is independent of $z$. In this case, the only condition from \eqref{HBd} is $c'_g{}^2=
\overline{\langle g,g\rangle}$. Because $\dot{gz}=\mu^{-g}\dot{z}$ while $\widetilde{gz}=\mu^g
\widetilde{z}$, it suffices to consider  $z$ with trivial $\widetilde{z}$. Moreover, because $\overline{\dot{z}}
=\dot{\overline{z}}$, if $\overline{\dot{z}(g)}\,\overline{\widetilde{z}(g)}$ is independent of $z$
then it must be real and hence in $\pm 1$. Thus $c'_g{}^2=\epsilon_{\langle,\rangle}(g)^2=\overline{\langle g,g
\rangle}$ and \eqref{HBd} is automatic.

This condition (independence of $z$) is automatically true when $n=n'$ (the case considered in \cite{iz3}) because then
$\widetilde{z}$ uniquely determines $z$. At this time, we don't know whether it is also
true when  
$n'>n$ --- for all we know, $n'>n$ is never realised.

When the half-braiding $\cE_{\alpha_g}$ exists, it is unique and defined by
$\mathcal{E}_{\alpha_g}(\alpha_h)=\langle g,h\rangle$ and $\mathcal{E}_{\alpha_g}(\rho)=
\epsilon_{\langle,\rangle}(g)\dot{z}(g)\widetilde{z}(g) U_g$, for any choice of $z$.
This then allows us to compute the corresponding parts of the $S$ and $T$ matrices.
In particular (assuming all half-braidings $\cE_{\alpha_g}$ exist), we have
\begin{align}
&&T_{\alpha_g,\alpha_g}=\cE_{\alpha_g}(\alpha_g)=\langle g,g\rangle\,,\\
&&S_{\alpha_g,\alpha_h}=\lambda^{-1}\cE_{\alpha_g}(\alpha_h)^*\cE_{\alpha_h}(\alpha_g)^*
=\lambda^{-1}\overline{\langle g,h\rangle}^2\,,\end{align}
where $\lambda$ is given in \eqref{globdim}.

In the remainder of this subsection we turn to $n'=0$ and $n'=n$.
When $n'=0$, \eqref{HBe} no longer applies and we have two 
half-braidings for $\alpha_g$ given by choosing either sign in $c_g'$. The remaining half-braidings
for $n'=0$ are  for $\sigma=\rho$, with
precisely $2n$ half-braidings, and precisely one each  for $\sigma=\alpha_g+\alpha_h$ for each
$g\ne h$. In this Tambara-Yamagami case, as analysed in Section 3 of \cite{iz3},
Tube$\,\Delta$ is isomorphic as a C$^*$-algebra to 
$\bbC^{4n}\oplus \left(M_{2\times 2}\right)^{n(n-1)/2}$, and  elementary expressions for the modular 
data fall out directly.

Thanks to Proposition 6, the case $n=n'$ reduces to that studied in \cite{iz3}, and so its tube
algebra is fully analysed in Section 6 of \cite{iz3}.
We find there is a unique half-braiding and simple summand $\bbC$ in Tube$\,\Delta$
for each $\sigma=\alpha_g$, while each $\sigma=\rho+\alpha_g$ corresponds to a unique
summand $M_{2\times 2}$ and half-braiding. $\sigma=\rho+\alpha_g+\alpha_h$ ($g\ne h$)
gives a unique braiding and summand $M_{3\times 3}$. Finally, there are exactly
$n(n+3)/2$ half-braidings with $\sigma=\rho$, and each contributes a $\bbC$ to Tube$\,\Delta$.
Thus 
\begin{equation}
\mathrm{Tube}\,\Delta_{n'=n}\cong \bbC^{n(n+5)/2}\oplus\left(M_{2\times 2}\right)^n\oplus
\left(M_{3\times 3}\right)^{n(n-1)/2}\,.\end{equation}

The modular data for $n'=n$ is described in \cite{iz3} as follows. 
First, 
find all  functions  $\xi:G\rightarrow \bbT$
and $\omega\in\bbT$, $\tau\in G$ such that
\begin{align}
&&\sum_g \xi(g)=\sqrt{n}\,\omega^2\,a(\tau)\,c^3-n\delta^{-1}\,,\label{xile1}\\
&&\overline{c}\sum_kb(g+k)\,\xi(k)=\omega^2\,c^3\,a(\tau)\,\overline{\xi(g+\tau)}-\sqrt{n}\delta^{-1}\,,\label{xile2}\\
&&\xi(\tau-g)=\omega\,c^4\,a(g)\,a(\tau-g)\,\overline{\xi(g)}\,,\label{xile3}\\
&&\sum_k\xi(k)\,b(k-g)\,b(k-h)=c^{-2}\,b(g+h-\tau)\,\xi(g)\,\xi(h)\,\overline{a(g-h)}-c^2\delta^{-1}\,.
\end{align}
There will be a total of $n(n+3)/2$ such triples $(\omega_j,\tau_j,\xi_j)$.

The $n(n+3)$ primaries fall into four classes:\smallskip

\noindent\textit{1.} $n$ primaries, denoted $\mathfrak{a}_g$, $g\in G$;

\noindent\textit{2.} $n$ primaries, denoted $\mathfrak{b}_h$ for  $h\in G$;

\noindent\textit{3.} $n(n-1)/2$ primaries, denoted $\mathfrak{c}_{k,l}=\mathfrak{c}_{l,k}$ for 
$k,l\in G$, $k\ne l$;

\noindent\textit{4.} $n(n+3)/2$ primaries, denoted $\mathfrak{d}_{j}$, corresponding to the
triples $(\omega_j,\tau_j,\xi_j)$.\smallskip

We can write the $S$ and $T$ matrices in block form as
\begin{align}
&&T=\mathrm{diag}(\langle g,g\rangle;\langle h,h\rangle;\langle k,l\rangle;\omega_j)\,,\qquad\qquad\qquad\qquad\\
&&S=\frac{1}{\lambda}{\scriptsize \left(\begin{matrix}\langle g,g'\rangle^{-2}&(\delta+1)\langle g,h'\rangle^{-2}
&(\delta+2)\overline{\langle g,k'+l'\rangle}&\delta\langle g,\tau_{j'}\rangle\\
(\delta+1)\langle h,g'\rangle^{-2}&\langle h,h'\rangle^{-2}&(\delta+2)\overline{\langle h,k'+l'\rangle}
& -\delta\langle h,\tau_{j'}\rangle\\ (\delta+2)\overline{\langle k+l,g'\rangle}&(\delta+2)\overline{
\langle k+l,h'\rangle}&(\delta+2)(\overline{\langle k,k'\rangle\langle l,l'\rangle}+\overline{
\langle k,l'\rangle\langle l,k'\rangle})&0\\ \delta\langle\tau_j,g'\rangle&-\delta\langle \tau_j,h'\rangle
&0&S_{j,j'}\end{matrix}\right)}\,,\end{align}
where
\begin{equation}\label{Sjjprime}
S_{j,j'}=\omega_j\omega_{j'}\sum_g\langle \tau_j+\tau_{j'}+g,g\rangle+\delta \omega_j\omega_{j'}
c^6a(\tau_j)a(\tau_{j'})n^{-1}\sum_{g,h}\overline{\xi_j(g)\xi_{j'}(h)\langle\tau_j-\tau_{j'}+h-g,h-g
\rangle}\,.\end{equation}
This is all perfectly simple, except for the $n(n+3)/2\times n(n+3)/2$ block $S_{j,j'}$.

\subsection{The modular data for the double of $G+n$ when $n$ is odd}

The point of this subsection is to compute the mysterious part \eqref{Sjjprime}. 
We will show that, rather unexpectedly, $S_{j,j'}$ is built up from a quadratic form on
an abelian group  of order $n+4$.

\medskip\noindent\textbf{Definition 3.} \textit{Let $G$ be any finite abelian group.
By a} nondegenerate quadratic form \textit{$Q$ on $G$ we mean a map $Q:G\rightarrow 
\bbQ/\bbZ$ such that $Q(-g)\equiv Q(g)$ (mod 1) for all $g\in G$, and $\langle,\rangle_Q:G\times G\rightarrow
\bbT$ defined by $\langle g,h\rangle_Q=e^{2\pi\i\,(Q(g+h)-Q(g)-Q(h))}$ is a nondegenerate 
symmetric pairing in the sense of Definition 1.} \medskip

For example, when $G=\bbZ_n$ for $n$ odd,  these are precisely 
  $Q(g)={mg^2}/n$ for any integer $m$ coprime to $n$. More generally, for $|G|$ odd, the
  nondegenerate quadratic forms and nondegenerate symmetric pairings are in natural
  bijection. In such a case, we can always write $G$ as $\bbZ_{n_1}\times\cdots\times\bbZ_{n_k}$
  where $Q$ restricted to each $\bbZ_{n_i}$ is nondegenerate and $\langle \bbZ_{n_i},\bbZ_{n_j}
  \rangle=1$ for $i\ne j$. When $|G|$ is even, things are more complicated but $G$ will have
precisely  $|G/2G|$ nondegenerate quadratic forms for each nondegenerate symmetric
pairing $\langle,\rangle$.  The map $a:G\rightarrow \bbT$ of
  Corollary 5 is the exponential of  a nondegenerate  quadratic form.

Given a nondegenerate quadratic form $Q$ and $a\in\bbZ$, define the \textit{Gauss sum} 
$$\alpha_Q(a)=\frac{1}{\sqrt{|G|}}\sum_{k\in G}\exp(2\pi\i\, a \,Q(k))\,.$$
For example,  note from \eqref{5.1} that the quantity $c^3$ of Corollary 5 is a Gauss sum.
Provided $aQ$ is nondegenerate, $\alpha_Q(a)$ will be a root of unity (this is a consequence
of Proposition 7(a) below). All Gauss sums needed in this paper can be computed from the
classical Gauss sums, corresponding to $G=\bbZ_n$ and   $Q(g)=mg^2/n$, which equal
$$\left\{\begin{matrix}
\left(\frac{am}{n}\right)&\mathrm{for}\ n\equiv_41\\  \i\left(\frac{am}{n}\right)&
\mathrm{for}\ n\equiv_43\\  
0 &\mathrm{for}\ n\equiv_4 2\\
(1\pm \i)\left(\frac{n}{am}\right)&\mathrm{when}\  n\equiv_40\ \mathrm{and}
\ a\equiv_4 \pm m\end{matrix}\right.\,,$$
where $\left(\frac{a}{b}\right)$ is the Jacobi symbol. For $n$ even, these classical Gauss
sums are not modulus 1, because $mg^2/n$ is degenerate in $\bbZ_n$.

\medskip\noindent\textbf{Proposition 7(a)} \textit{Let $Q$ be a nondegenerate quadratic form
on any abelian group $G$. Define  matrices 
\begin{equation}S^Q_{g,h}=\frac{\alpha}{\sqrt{|G|}}\overline{\langle g,h\rangle}_Q\ ,\qquad\qquad
T^Q_{g,h}=\beta\,\delta_{g,h}\exp(2\pi \i\, Q(g))\,,\end{equation}
for any $\alpha,\beta\in\bbC$. Then $S^Q,T^Q$ define modular data iff $\alpha=\pm 1$ and
$\beta^3=\alpha\alpha_Q(1)$. In this case, the identity is $\mathfrak{a}_0$.}\medskip

\noindent\textbf{(b)} \textit{Let $G,G'$ be abelian groups of odd order $n$ and $n+4$ respectively.
Choose any nondegenerate quadratic forms $Q$ and $Q'$ on them, 
 and write $\langle g,h\rangle=\langle g,h\rangle^{(n+1)/2}_Q$ and $\langle \beta,\gamma\rangle'
 =\langle \beta,\gamma\rangle_{Q'}^{(n+5)/2}$ for all $g,h\in G$, $\beta,\gamma\in G'$, so
 $Q(g)=\langle g,g\rangle$ and $Q'(\gamma)=\langle \gamma,\gamma\rangle'$. 
Let $\Phi$ consist of the following $n(n+3)=n+n+n(n-1)/2+n(n+3)/2$ elements: $\mathfrak{a}_g$
 $\forall g\in G$; $\mathfrak{b}_h$  $\forall h\in G$;  $\mathfrak{c}_{k,l}=\mathfrak{c}_{l,k}$  
$\forall k,l\in G$ with $k\ne l$; and  $\mathfrak{d}_{m,\gamma}=\mathfrak{d}_{m,-\gamma}$  
$\forall m\in G,\gamma \in G'$,  $\gamma\ne 0$. Define
\begin{align}
&&T^{Q,Q'}=\mathrm{diag}(\langle g,g\rangle;\langle h,h\rangle;\langle k,l\rangle;\langle m,m\rangle\,\langle \gamma,\gamma\rangle')\,,\qquad\qquad\qquad\qquad\qquad\qquad\cr
&&S^{Q,Q'}=\frac{1}{\lambda}\scriptsize{\left(\begin{matrix}\overline{\langle g,g'\rangle}{}^2&
(\delta+1)\overline{\langle g,h'\rangle}{}^2
&(\delta+2)\overline{\langle g,k'+l'\rangle}&\delta\overline{\langle g,{m'}\rangle}{}^2\\
(\delta+1)\overline{\langle h,g'\rangle}{}^2&\overline{\langle h,h'\rangle}{}^2&(\delta+2)\overline{\langle h,k'+l'\rangle}
& -\delta\overline{\langle h,m'\rangle}{}^2\\ (\delta+2)\overline{\langle k+l,g'\rangle}&(\delta+2)\overline{
\langle k+l,h'\rangle}&(\delta+2)\left(\overline{\langle k,k'\rangle}\overline{\langle l,l'\rangle}+\overline{
\langle k,l'\rangle}\overline{\langle l,k'\rangle}\right)&0\\ \delta\overline{\langle m,g'\rangle}^2&-\delta\overline{\langle m,h'\rangle}{}^2
&0&-\delta\overline{\langle m,m'\rangle}\left(\langle \gamma,\gamma'\rangle'+\overline{\langle \gamma,\gamma'\rangle'}\right)
\end{matrix}\right)}\nonumber\end{align}
where $\lambda$ is given in \eqref{globdim}.
Then these define modular data iff $\alpha_{Q}(1)\,\alpha_{Q'}(1)=-1$. The identity is
$\mathfrak{a}_0$.}\medskip

 The straightforward proof is by direct 
calculation: $S^2=C$, $S^*=CS$, $ST^*S=TS^*T$, and Verlinde's formula \eqref{verl}.
In (a), $\alpha^2=1$ arises from the requirement that $S^2$ be a permutation matrix.
The conditions $\alpha^3\beta^3=\alpha_Q(-1)$ and $\alpha_Q(-1)\alpha_{Q'}(-1)=-1$
for (a) and (b) respectively both come from $ST^*S=TS^*T$. 
 We find that the fusion coefficients of part (a) are
$N_{g,h}^k=\delta_{k,g+h}$, every primary is a simple-current, and 
charge-conjugation  $C$ acts by $-1$. In (b), charge-conjugation sends $\mathfrak{a}_g\mapsto \mathfrak{a}_{-g}$, 
$\mathfrak{b}_h\mapsto \mathfrak{b}_{-h}$, $\mathfrak{c}_{k,l}\mapsto \mathfrak{c}_{-k,-l}$,
$\mathfrak{d}_{m,\gamma}\mapsto\mathfrak{d}_{-m,\gamma}$.
The nonzero fusion coefficients there are
\begin{align}
N_{\mathfrak{a}_g,\mathfrak{a}_{h},\mathfrak{a}_{k}}=&\,
N_{\mathfrak{a}_g,\mathfrak{b}_{h},\mathfrak{b}_{k}}=
N_{\mathfrak{b}_g,\mathfrak{b}_{h},\mathfrak{b}_{k}}=
N_{\mathfrak{b}_g,\mathfrak{b}_{h},\mathfrak{d}_{k,\gamma}}=\delta(g+h+k)\,;\cr
N_{\mathfrak{a}_g,\mathfrak{c}_{h,k},\mathfrak{c}_{h',k'}}=&\,\delta(g+h+h')\delta(g+k+k')+
\delta(g+h+k')\delta(g+h'+k)\in\{0,1\}\,;\cr
N_{\mathfrak{a}_g,\mathfrak{d}_{h,\beta},\mathfrak{d}_{k,\gamma}}=&\,\delta(g+h+k)\delta_{\beta,\gamma}\,;\cr
N_{\mathfrak{b}_g,\mathfrak{b}_{h},\mathfrak{c}_{k,k'}}=&\,
N_{\mathfrak{b}_g,\mathfrak{c}_{k,k'},\mathfrak{d}_{h,\gamma}}=\delta(2g+2h+k+k')\,;\cr
N_{\mathfrak{b}_g,\mathfrak{c}_{h,k},\mathfrak{c}_{h',k'}}=&\,\delta(2g+h+k+h'+k')+
\delta(g+h+h')\delta(g+k+k')\cr&\qquad\qquad\qquad\qquad\qquad
+\delta(g+h+k')\delta(g+h'+k)\in\{0,1,2\}\,;\cr
N_{\mathfrak{b}_g,\mathfrak{d}_{h,\beta},\mathfrak{d}_{k,\gamma}}=&\,\delta(g+h+k)(1-
\delta_{\beta,\gamma}) \,;\cr
N_{\mathfrak{c}_{g,h},\mathfrak{c}_{g',h'},\mathfrak{d}_{k,\gamma}}=&\,\delta(g+h+g'+h'+2k)\,;\cr
N_{\mathfrak{c}_{g,h},\mathfrak{d}_{k,\beta},\mathfrak{d}_{k',\gamma}}=&\,\delta(g+h+2k+2k')\,;\cr
N_{\mathfrak{c}_{g,h},\mathfrak{c}_{g',h'},\mathfrak{c}_{g'',h''}}=&\,\delta(g+h+g'+h'+g''+h'')
(1+\delta(g+g'+g'')+\delta(g+h'+g'')\cr&\qquad\qquad\qquad\qquad+\delta(g+g'+h'')+\delta(g+h'+h'')) \in\{0,1,2\}\,;\cr
N_{\mathfrak{d}_{g,\gamma},\mathfrak{d}_{g',\gamma'},\mathfrak{d}_{g'',\gamma''}}=&\, 
\delta(g+g'+g'')(1-\delta(\gamma+\gamma'+\gamma'')-\delta(\gamma-\gamma'+\gamma'')
-\delta(\gamma+\gamma'-\gamma'')\cr&\quad\qquad\qquad\qquad-\delta(\gamma-\gamma'-\gamma''))\in\{0,1\}\,,
\nonumber\end{align}
where we write $\delta(g)=1$ or 0 depending
on whether or not $g=0$.

We'll let $\cM\cD_{G,G'}(Q,Q')$ denote the modular data of (b).
Of course, associating SL$_2(\bbZ)$ representations to quadratic forms is an old story. See
for instance \cite{NW}, who
study these in  similar generality  (though their $G$ are $p$-groups, and they require
$\beta=1$), and call these Weil representations. To our knowledge, Proposition 7(b)
is completely new, but what is more important is its relation to near-group doubles:

\medskip\noindent\textbf{Conjecture 2.} \textit{When $|G|=n$ is odd, the modular data for $G+n$ is
$\cM\cD_{G,G'}(Q,Q')$, where $Q$ is the nondegenerate quadratic form on $G$
corresponding to $\langle,\rangle$, and $Q'$ is a nondegenerate quadratic form on
some abelian group $G'$ of order $n+4$.}\medskip

This is true for all groups $G$ of odd order $\le 13$ (recall in Table 2). The easiest way
to verify modular datum are equivalent is to first identify their $T$ matrices (straightforward
since they must have finite order) and then compare the floating point values of the $S$
matrices --- see Section 4.1 of \cite{EGh} for details. In all cases in Table 2, $G'
\cong \bbZ_{n+4}$, except for the first entry for $G=\bbZ_5$ when $G'=\bbZ_3\times\bbZ_3$. 
The quadratic form
$Q'$ is then identified in Table 2 by the integer $m'$ in the $Q'$ column. For the first entry of
$G=\bbZ_5$, $Q'(\gamma_1,\gamma_2)=(\gamma_1^2+\gamma_2^2)/3$.

Even for $n=3$, this is vastly simpler than the modular data as it appears in Example A.1 of \cite{iz3}.
In fact we have no direct proof that they are equal for $n=3$ --- our proof that they are the same
is that they both yield nonnegative integer fusions, they have the same $T$ matrices,
and their $S$ matrices are numerically close. The simplicity of this modular data
$\cM\cD_{Q,Q'}(G,G')$ supports our claim that the doubles of these near-group categories
$G+n$ should not be regarded as exotic. We would expect that these doubles are realised
by rational conformal nets of factors, and by rational vertex operator algebras.




The quantum-dimensions $S_{x,0}/S_{0,0}$ are  $1,\delta+1,\delta+2,\delta$
for primaries of type $\mathfrak{a},\mathfrak{b},\mathfrak{c},\mathfrak{d}$ respectively.
We see from the above that the $\mathfrak{a}_g$ are {simple-currents}, and obey the
fusions $\mathfrak{a}_g*\mathfrak{a}_h=\mathfrak{a}_{g+h}$, $\mathfrak{a}_g*\mathfrak{b}_h=
\mathfrak{b}_{g+h}$, $\mathfrak{a}_g*\mathfrak{c}_{h,k}=\mathfrak{c}_{g+h,g+k}$, and 
$\mathfrak{a}_g*\mathfrak{d}_{h,\gamma}=\mathfrak{d}_{g+h,\gamma}$. They form a
group isomorphic to $G$, and act without fixed-points. They supply the ultimate explanation
for the $G$-action of Proposition 6.7 of \cite{iz3}. The phases $\varphi_g(x)$ defined by
$S_{\mathfrak{a}_gx,y}=\varphi_g(y)S_{x,y}$ are $\langle g,h\rangle$ for $\mathfrak{a}_h,\mathfrak{b}_h,
\mathfrak{d}_{h,\gamma}$, and $\langle g,k+l\rangle$ for $\mathfrak{c}_{k,l}$.

The Galois symmetry is useful in understanding the modular invariants. For $\ell\in\bbZ_{n(n+4)}^\times$,
$\mathfrak{a}_g^\ell=\mathfrak{a}_{\ell g}$ or $\mathfrak{b}_{\ell g}$ depending on whether or
not the Jacobi symbol $\left({\ell\over n(n+4)}\right)$ equals 1. Similarly, 
$\mathfrak{b}_g^\ell=\mathfrak{b}_{\ell g}$ or $\mathfrak{a}_{\ell g}$ respectively. Finally,
$\mathfrak{c}_{g,h}^\ell=\mathfrak{c}_{\ell g,\ell h}$ and $\mathfrak{d}_{g,\gamma}^\ell=
\mathfrak{d}_{\ell g,\ell \gamma}$. All parities $\epsilon_\ell(x)=+1$ except for
$\epsilon_\ell(\mathfrak{d}_{g,\gamma})=  \left({\ell\over n(n+4)}\right)$. The requirement of
a coherent Galois symmetry is what led us to the simplified modular data given above.

A \textit{modular invariant} is a matrix $Z$ with nonnegative integer entries (often formally written as a generator function
$\cZ=\sum_{a,b}Z_{a,b}ch_a\overline{ch_b}$), with $Z_{0,0}=1$, which commutes with the modular data $S,T$. It is called \textit{type I} if $\cZ$  can be written as a sum of squares.
There are exactly 3 type I modular invariants when both $n$ and $n+4$ are prime (e.g. for $n=3$):
\begin{align}&&\mathcal{Z}_1=\sum_g|\mathfrak{a}_g|^2+\sum_g|\mathfrak{b}_g|^2+\sum_{g,h}|\mathfrak{c}_{g,h}
|^2+\sum_{g,\gamma}|\mathfrak{d}_{g,\gamma}|^2\,,\qquad\qquad\qquad\nonumber\\\nonumber
&&\mathcal{Z}_2=\sum_g|\mathfrak{a}_g+\mathfrak{b}_g|^2+2\sum_{g,h}|\mathfrak{c}_{g,h}|^2\,,
\qquad\mathcal{Z}_3=|\mathfrak{a}_0+\mathfrak{b}_0+\sum_{g\ne 0}\mathfrak{c}_{g,0}|^2\,.
\nonumber\end{align}
The most important modular invariants are the \textit{monomial} ones, of form $\cZ=|\Sigma Z_{a,0}
ch_a|^2$,
as explained in Section 1.3 of \cite{EGh}, as they give a canonical endomorphism
$\theta$ as a sum of sectors, and
can be used to recover the original  system from its double. This is M\"uger's forgetful functor \cite{m1}.
For example, there are exactly 3 monomial modular
invariants for the modular data of the double of the  Haagerup subfactor; these should correspond
bijectively to
the three systems found in \cite{GrSn} which are Morita equivalent to the principal even
Haagerup system (see their Theorem 1.1), as each of those must correspond to a monomial
modular invariant. 

 We see however that for the $\cM\cD_{G,G'}(Q,Q')$ modular data there
 is  only one generic monomial modular invariant, namely $\mathcal{Z}_3$. This
suggests that Grossman-Snyder perhaps isn't as interesting here as it was for the Haagerup
(at least not for general $n$). On the other hand, 
recall our comments earlier that the type $G+n$ systems with $n=n'=\nu^2$ may be related to
the Haagerup-Izumi system for groups of order $\nu$. Consider first $G=\bbZ_n$ for $n=\nu^2$  a perfect square
and write $H=\nu G\cong\bbZ_\nu$; then
$\cM\cD_{G,G'}(Q,Q')$ has at least one other monomial invariant, namely
\begin{equation}
\cZ_4=|\sum_{h\in H}\mathfrak{a}_h+\sum_{h\in H}\mathfrak{b}_h+2\sum_{h<h'\in H}\mathfrak{c}_{h,h'}|^2\,.\end{equation}
Alternatively, when  $G=H_1\times H_2$ where each $H_i\cong\bbZ_\nu$, and $\langle
H_1,H_2\rangle_Q=1$ (which can always be arranged), another monomial invariant is
\begin{equation}
\cZ_4'=|\mathfrak{a}_0+\mathfrak{b}_0+\sum_{h\in H_1,h'\in H_2}\mathfrak{c}_{(h,0),(0,h')}|^2\,.\end{equation}
We would expect that systems of type $\bbZ_{\nu^2}+\nu^2$ or $\bbZ_\nu\times\bbZ_\nu+\nu^2$
should have nontrivial quantum subgroups in the sense of \cite{GrSn}.

It isn't difficult to see why  $n+4$ arises here, i.e. why it can't be replaced by some other
positive integer $n'$. In particular, after some work, 
the nonzero  fusions of the form $N_{\mathfrak{b},\mathfrak{b},\mathfrak{b}}$ are found to
be $4/(n'-n)$;,and the $ST^*S={T}S^*{T}$ calculation requires the product of
Gauss sums for $G$ and  $G'$ to be $-1$, which forces $4|(n'-n)$.

When $G$ has even order, the situation is similar but (as always with $n$ even) somewhat
messier; we will provide its modular data elsewhere. Again we have $n$ simple-currents
(the $\mathfrak{a}_g$), but for each $g\in G$ of order 2, $\mathfrak{a}_g$ now has 
$n/2$ fixed-points, which complicates things. The $T$ entries for the first several
even $G$ are provided by the pairs $(m',m'')$  in the $Q'$ column of Table 2, and from this the
$S$ matrix follows quickly from the equations of the last subsection. In particular,
$$T_{\mathfrak{d}_{g,\gamma},\mathfrak{d}_{g,\gamma}}=\left\{\begin{matrix}
\xi_n^{\langle g,g\rangle}\,\xi_{n+4}^{m'\gamma^2}&\mathrm{if}\ \gamma+n/2\ \mathrm{is\ odd}\\
\xi_n^{\langle m,m-1\rangle}\,\xi_{n(n+4)}^{m''\,(1+n\gamma)^2}&\mathrm{if}\ \gamma+n/2\ \mathrm{is\ even}\end{matrix}\right.$$
where $\tau_\gamma=0,1$ for $\gamma+n/2$ odd respectively even. Here $1\le \gamma\le(n+4)/2$
and $g\in G$ except for $\gamma=(n+4)/2$ when $g\in G/2$.

Recall our observation in Section 3 of \cite{EGh} that the modular data of the double of the Haagerup-Izumi series
at $G=\bbZ_{n}$ resembles that of the affine algebra so$(n^2+4)^{(1)}$ at level 2.
The analogous statement here is that the modular datum of the double of  type $\bbZ_n+n$ near-group systems 
resemble that of  the affine algebra 
so$(n+4)^{(1)}$ at level 2. In particular, for an appropriate choice of $Q'$
(corresponding to  $m'=(n+3)/2$), this recovers $T_{\mathfrak{d}_{0,\gamma},\mathfrak{d}_{0,
\gamma}}$ and $S_{\mathfrak{d}_{0,\gamma},\mathfrak{d}_{0,\gamma'}}$.  
This could hint at ways to construct the corresponding vertex operator
algebra.  



\newcommand\biba[7]   {\bibitem{#1} {#2:} {\sl #3.} {\rm #4} {\bf #5,}
                    {#6 } {#7}}
                    \newcommand\bibx[4]   {\bibitem{#1} {#2:} {\sl #3} {\rm #4}}

\def\ASENS            {Ann. Sci. \'Ec. Norm. Sup.}
\def\AM   {Acta Math.}
   \def\AnM              {Ann. Math.}
   \def\CMP              {Commun.\ Math.\ Phys.}
   \def\IJM              {Internat.\ J. Math.}
   \def\JAMS             {J. Amer. Math. Soc.}
\def\JFA              {J.\ Funct.\ Anal.}
\def\JMP              {J.\ Math.\ Phys.}
\def\JRA              {J. Reine Angew. Math.}
\def\JPAA             {J.\ Pure Appl.\ Algebra.}
\def\JSP              {J.\ Stat.\ Physics}
\def\LMP              {Lett.\ Math.\ Phys.}
\def\RMP              {Rev.\ Math.\ Phys.}
\def\RNM              {Res.\ Notes\ Math.}
\def\RIMS             {Publ.\ RIMS.\ Kyoto Univ.}
\def\Inv              {Invent.\ Math.}
\def\npbp             {Nucl.\ Phys.\ {\bf B} (Proc.\ Suppl.)}
\def\nupb             {Nucl.\ Phys.\ {\bf B}}
\def\nup              {Nucl.\ Phys. }
\def\nupp             {Nucl.\ Phys.\ (Proc.\ Suppl.) }
\def\adma             {Adv.\ Math.}
\def\coma             {Con\-temp.\ Math.}
\def\PAMS             {Proc. Amer. Math. Soc.}
\def\PJM              {Pacific J. Math.}
\def\ijmp             {Int.\ J.\ Mod.\ Phys.\ {\bf A}}
\def\jpa              {J.\ Phys.\ {\bf A}}
\def\PLB              {Phys.\ Lett.\ {\bf B}}
\def\RIMS             {Publ.\ RIMS, Kyoto Univ.}
\def\Top               {Topology}
\def\TAMS             {Trans.\ Amer.\ Math.\ Soc.}
\def\Duke              {Duke Math.\ J.}
\def\K                 {K-theory}
\def\JOP               {J.\ Oper.\ Theory}
\def\JKT               {J.\ Knot Theory and its Ramifications}

\vspace{0.2cm}\addtolength{\baselineskip}{-2pt}
\begin{footnotesize}
\noindent{\it Acknowledgement.}

The authors thank the Erwin-Schr\"odinger-Institute,  Cardiff School of 
Mathematics, and
Swansea University Dept of Computer Science,  for generous
hospitality while researching this
paper. They thank 
 Masaki Izumi for sharing with us an early draft of \cite{iz4} and informing
us of \cite{EGO}, and Eric Rowell for informing
us of \cite{Sie}. Their research was supported in part by EPSRC,   
EU-NCG Research Training Network: MRTN-CT-2006 031962,
and NSERC.

\end{footnotesize}

\end{document}